\documentclass[11pt]{article}
\usepackage{xcolor}
\usepackage{amsmath,amssymb,amsfonts,amsthm,amscd,mathrsfs}
\usepackage[T1]{fontenc}
\usepackage{lmodern}
\usepackage{graphicx,float,url}
\usepackage[hidelinks]{hyperref}
\usepackage{cite}
\allowdisplaybreaks[4]

\newcommand{\e}{\eps}
\newcommand{\p}{\partial}
\newcommand{\cD}{\mathcal D}
\newcommand{\cA}{\mathcal A}
\newcommand{\cL}{\mathcal L}
\newcommand{\cJ}{\mathcal J}
\newcommand{\n}{\ensuremath{\nonumber}}
\newcommand{\eps}{\ensuremath{\varepsilon}}
\newcommand{\bigO}{\mathcal{O}}
\newtheorem{theorem}{Theorem}[section]
\newtheorem{lemma}[theorem]{Lemma}
\newtheorem{proposition}[theorem]{Proposition}
\newtheorem{corollary}[theorem]{Corollary}
\newtheorem{remark}[theorem]{Remark}

\numberwithin{equation}{section}
\newcommand{\beq}{\begin{equation}\begin{aligned}}
		\newcommand{\eeq}{\end{aligned}\end{equation}}
\title{Validity of three-dimensional axisymmetric steady Prandtl expansion}
\author{Yan Guo, Hailiang Li, Xiangyu Ma, Mingyue Zhang}
\date{}

\begin{document}
\maketitle

\begin{abstract}
We study steady axisymmetric Navier--Stokes flows with swirl in the circular cylinder $\mathcal C_L=(0,L)\times B_2(0,1)$ with the no-slip condition on the stationary lateral wall.  For the leading Euler velocity $(u_e^0,w_e^0,\theta_e^0)=(1,0,r K(r))$, where $K(r)=H(r^2)$ is a prescribed smooth radial function, we construct and justify a Prandtl expansion as $\e \to 0$ for sufficiently small fixed $L$. The expansion accounts for both axial and azimuthal boundary layers, including the case of a Blasius leading axial profile with nonzero swirl. Weighted estimates at the symmetry axis and
coupled remainder estimates yield a uniform velocity error bound.
\end{abstract}
\section{Introduction}\label{sec:introduction}

Consider steady incompressible axisymmetric Navier--Stokes flows
with swirl in the circular cylinder
$\mathcal C_L=(0,L)\times B_2(0,1)$. In cylindrical coordinates,
the velocity and pressure have the form
\[
 \mathbf U=U(z,r)\mathbf e_z+V(z,r)\mathbf e_r
               +\Theta(z,r)\mathbf e_\theta,\qquad P=P(z,r),
\]
where $U,V,\Theta$ are the axial, radial, and azimuthal velocities.
For viscosity $\e>0$, the system reads:
\begin{equation}\label{eq:NS-axisymmetric-zr}
 \mathbf U\cdot\nabla\mathbf U+\nabla P=\e\Delta\mathbf U,
 \qquad \nabla\cdot\mathbf U=0,
\end{equation}
with the no-slip condition $(U,V,\Theta)|_{r=1}=(0,0,0)$ at the
stationary lateral wall. As viscosity $\e\to0$, the inviscid flow satisfies
the impermeability condition, but its axial and azimuthal velocities
need not vanish at the wall. Following Prandtl \cite{prandtl1905uber},
we introduce boundary layers of thickness $O(\sqrt\e)$ in these two
tangential components to satisfy the no-slip condition.

For three-dimensional axisymmetric steady flows with swirl, the
azimuthal velocity couples to the meridional flow through the
centrifugal force. We construct Euler--Prandtl profiles and estimate the coupled remainders throughout the cylinder, including near the symmetry axis.

Introduce the variables and velocities
\begin{equation}\label{eq:change-of-variables}
\begin{aligned}
 &x=z,\qquad Y=1-r,\qquad y=\frac{Y}{\sqrt\e},\\
 &W(x,Y)=-V(x,1-Y),\qquad S=\frac{\Theta}{r}.
\end{aligned}
\end{equation}
Here $W$ is the inward radial velocity and $S$ is the angular velocity.
We retain the notation $U,\Theta,P$ in $(x,Y)$.
For the outer flow, let $K(r)=H(r^2)$, where
$H\in C^\infty([0,1])$ is prescribed and $K$ satisfies
$K(1)=\kappa$ for a fixed constant $\kappa\ge0$.
We take the leading Euler state to be
\beq \label{ue0}
 u_e^0=1,\qquad w_e^0=0,\qquad
 S_e^0=K(r),\qquad \theta_e^0=rK(r),\qquad
 p_e^0=p_*+\int_0^r tK(t)^2\,dt.
\eeq

Fix $N\ge8$. For sufficiently small fixed $L$, and
then for $\e$ sufficiently small depending on $L$, we seek a
solution of the form
\begin{equation}\label{eq:sw-remainder-ansatz}
\begin{aligned}
 &U=\sum_{i=0}^N\e^{\frac{i}{2}}\bigl(u_e^i(x,Y)+u_p^i(x,y)\bigr)
       +\e^{\frac{7}{2}}u^\e(x,Y)=:u_s+\e^{\frac{7}{2}}u^\e,\\
 &W=\sum_{i=0}^N\e^{\frac{i}{2}}\bigl(w_e^i(x,Y)+\sqrt\e\,w_p^i(x,y)\bigr)
       +\e^{\frac{7}{2}}w^\e(x,Y)=:w_s+\e^{\frac{7}{2}}w^\e,\\
 &S=\sum_{i=0}^N\e^{\frac{i}{2}}\bigl(S_e^i(x,Y)+S_p^i(x,y)\bigr)
       +\e^{\frac{7}{2}}s^\e(x,Y)=:s_s+\e^{\frac{7}{2}}s^\e,\\
 &P=\sum_{i=0}^N\e^{\frac{i}{2}}\bigl(p_e^i(x,Y)+\sqrt\e\,p_p^{i+1}(x,y)\bigr)
       +\e^{\frac{7}{2}}p^\e(x,Y)=:p_s+\e^{\frac{7}{2}}p^\e.
\end{aligned}
\end{equation}
Here $[u_e^i,w_e^i,S_e^i]$ are the Euler profiles in the variables
$(x,Y)$, and $[u_p^i,w_p^i,S_p^i]$ are the Prandtl profiles in the
variables $(x,y)$. For simplicity, we take $u_e^0=1$, but any positive
constant also works.

We construct the Euler and Prandtl profiles in
Section~\ref{sec:profiles-construction}. The profile data and the
untruncated profiles are independent of $\e$ at every order, and the
finite profile estimates are uniform in $L$. The Prandtl profiles are
constructed for $y\in\mathbb R_+$ and then truncated as in
Section~\ref{subsec:profile-truncation}. The resulting approximate velocity, given by \eqref{eq:profile-composite},
satisfies the continuity equation exactly and vanishes at the lateral wall.

With these profiles fixed, we now state our main result on the validity
of the expansion \eqref{eq:sw-remainder-ansatz}.

\begin{theorem}[Main Result]
\label{thm:main}
Assume that the profile data satisfy the hypotheses of
Section~\ref{subsec:swirl-leading}. There is
$\e_0(L)\in(0,L^{\frac{1}{8}}]$ such that, for $0<\e\le\e_0(L)$,
the coupled remainder problem \eqref{eq:sw-coupled-remainder} for
these fixed profiles has a unique solution in the closed unit ball
of $\mathscr Z_\e$ defined in Section~\ref{subsec:overview}.
The velocity and pressure reconstructed from \eqref{eq:sw-remainder-ansatz}
satisfy \eqref{eq:NS-xY} with the end conditions
\eqref{eq:sw-physical-end-data} and no slip at the lateral wall.
The pressure is unique up to an additive constant. The Navier-Stokes solutions satisfy
\begin{equation}\label{M2}
\begin{aligned}
 &\|U-u_s\|_\infty+\|W-w_s\|_\infty\le C_{L,K}\e^{\frac{11}{4}},\\
 &\|\Theta-rs_s\|_\infty\le C_{L,K}\e^3.
\end{aligned}
\end{equation}
\end{theorem}

\subsection{Notation}\label{subsec:notation}

Let $\Omega=(0,L)\times(0,1)$ in the wall coordinates $(x,Y)$.
On $\Omega$, write $(f,g):=\int_\Omega f g\,dx\,dY$ and
$\|f\|:=\|f\|_{L^2(\Omega)}$. Boundary norms carry a subscript
indicating the side. We use $\nabla=(\partial_x,\partial_Y)$ and
$\langle t\rangle=(1+t^2)^{\frac{1}{2}}$.

The radial weights used below are $w=r^{\frac32}$ and
$\Lambda(r)=\left(1+\log\frac1r\right)^{\frac12}$, and we denote
the radial differential operators by
$\cL=\partial_Y^2-\frac3r\partial_Y$ and
$\bar\cL=\partial_Y^2+\frac3r\partial_Y$.
Let $B_n=B_n(0,1)$ be the unit ball in $\mathbb R^n$. For $F(x,Y)$, write
$\widetilde F(x,r)=F(x,1-r)$ and define its radial lift by
$F^\sharp(x,\xi)=F(x,1-|\xi|)$ on
$\Omega^\sharp=(0,L)\times B_4$.
Polar integration gives
$\|F^\sharp\|_{L^2(\Omega^\sharp)}^2=|\mathbb S^3|\|Fw\|^2$,
and the radial Laplacian satisfies
$(\cL F)^\sharp=\Delta_\xi F^\sharp$.
The weighted Sobolev spaces and norms are defined by
\[
\begin{aligned}
 &\mathcal H^k(w)=\{F:F^\sharp\in H^k(\Omega^\sharp)\},\\
 &\|F\|_{\dot H^k(w)}
   =|\mathbb S^3|^{-\frac{1}{2}}\|\nabla_{x,\xi}^kF^\sharp\|_{L^2(\Omega^\sharp)},\qquad
 \|F\|_{H^k(w)}^2=\sum_{j=0}^k\|F\|_{\dot H^j(w)}^2.
\end{aligned}
\]
Here $\nabla_{x,\xi}^k$ is the full Cartesian derivative tensor
with its Euclidean norm. Unless a slice $B_4$ is specified, $L^2$ norms of lifted functions
are taken over $\Omega^\sharp$. For unlifted functions, derivatives
and $W^{k,\infty}$ norms are taken in $(x,Y)$.

In the profile construction, we use
$\Omega_p=(0,L)\times\mathbb R_+$ for the domain of the inner
variables. A tilde on a Prandtl profile denotes the raw profile on
$\Omega_p$, independent of $\e$. The corresponding untilded profile
includes the final truncation.
In formulas for the profile equations, we adopt the convention that
any profile with a negative index is zero, and that a sum whose upper
limit is smaller than its lower limit is an empty sum and equals zero.

Fix a smooth nonincreasing cutoff $\chi:[0,\infty)\to[0,1]$ equal to one on $[0,1/4]$ and zero on $[1/2,\infty)$.
For the boundary data, fix $\varphi\in C_c^\infty([0,\infty))$ with
$\varphi(0)=1$ and $\int_0^\infty\varphi(y)\,dy=0$.

In the a priori estimates, $a\lesssim b$ means
$a\le Cb$, where $C$ is independent of small $L,\e$. The constant
may depend on the fixed profile data and finite orders specified in
Section~\ref{sec:profiles-construction}. We write $a\sim b$ when
$a\lesssim b\lesssim a$. When dependence on $L$ is allowed, we
write $C_L$ or $C_{L,K}$. Constants in the preliminary regularity, trace, and
norm-equivalence estimates may depend on both $L$ and $\e$, as
indicated in those estimates. We use $o_L(1)$ for a quantity whose
absolute value is bounded by $CL^\alpha$ for some $\alpha>0$, with
$C$ independent of small $L,\e$.

\subsection{Overview of the proof}\label{subsec:overview}

In $(x,Y)$ coordinates, the Navier--Stokes equations are
\begin{equation}\label{eq:NS-xY}
\left\{
\begin{aligned}
 &UU_x+WU_Y+P_x
 =\e\left(U_{xx}+U_{YY}-\frac1rU_Y\right),\\
 &UW_x+WW_Y+\frac{\Theta^2}{r}+P_Y
 =\e\left(W_{xx}+W_{YY}-\frac1rW_Y-\frac{W}{r^2}\right),\\
 &U\Theta_x+W\Theta_Y-\frac Wr\Theta
 =\e\left(\Theta_{xx}+\Theta_{YY}
                  -\frac1r\Theta_Y-\frac{\Theta}{r^2}\right),\\
 &U_x+W_Y-\frac Wr=0,\\
 &(U,W,\Theta)|_{Y=0}=(0,0,0).
\end{aligned}
\right.
\end{equation}
At the two ends, we impose
\begin{equation}\label{eq:sw-physical-end-data}
\begin{aligned}
 &(U,W_x,\Theta)|_{x=0}
   =(u_s,w_{s,x},rs_s)|_{x=0},\\
 &(W,W_{xx},\Theta_x)|_{x=L}
   =(w_s,w_{s,xx},r(s_s)_x)|_{x=L}.
\end{aligned}
\end{equation}
The profiles determining these end data are constructed in
Section~\ref{sec:profiles-construction}.
Substituting \eqref{eq:sw-remainder-ansatz} into \eqref{eq:NS-xY} gives
\begin{equation}\label{reminder uw}
\left\{
\begin{aligned}
 &u_su_x^\e+w_su_Y^\e+u^\e u_{s,x}+w^\e u_{s,Y}+p_x^\e
 -\e\left(u_{xx}^\e+u_{YY}^\e-\frac1r u_Y^\e\right)\\
 &\qquad=-\e^{-\frac72}E_1
          -\e^{\frac72}(u^\e u_x^\e+w^\e u_Y^\e),\\
 &u_sw_x^\e+w_sw_Y^\e+u^\e w_{s,x}+w^\e w_{s,Y}
       +2rs_s s^\e+p_Y^\e-\e\left(w_{xx}^\e+w_{YY}^\e
                        -\frac1r w_Y^\e-\frac{w^\e}{r^2}\right)\\
 &\qquad=-\e^{-\frac72}E_2
          -\e^{\frac72}(u^\e w_x^\e+w^\e w_Y^\e+r(s^\e)^2),\\
 &(u^\e)_x+(w^\e)_Y-\frac{w^\e}{r}=0,\\
 &u^\e|_{Y=0}=w^\e|_{Y=0}=0.
\end{aligned}
\right.
\end{equation}
where
\begin{equation}\label{eq:physical-profile-residuals}
\begin{aligned}
 &E_1={}u_su_{s,x}+w_su_{s,Y}+p_{s,x}
 -\e\left(u_{s,xx}+u_{s,YY}-\frac1r u_{s,Y}\right),\\
 &E_2={}u_sw_{s,x}+w_sw_{s,Y}+rs_s^2+p_{s,Y}-\e\left(w_{s,xx}+w_{s,YY}
                     -\frac1r w_{s,Y}-\frac{w_s}{r^2}\right).
\end{aligned}
\end{equation}

Define the stream function $\Psi$ by $\Psi_Y=ru^\e$,
$\Psi_x=-rw^\e$, and $\Psi(x,0)=0$. Set $G=-\frac{\Psi}{r^2}$ and $q=\frac{G}{u_s}$.
Then $u^\e=2G-rG_Y$ and $w^\e=rG_x$.
Taking the meridional curl of the preceding equations and including
the scalar equation, we obtain
\begin{equation}\label{eq:sw-coupled-remainder}
\left\{
\begin{aligned}
 &\partial_xR[q]-\e\cA^2G+\cJ[G]
                         +2\partial_x(s_s s^\e)\\
 &\qquad=-\e^{-\frac72}\mathfrak R_\beta
 -\e^{\frac72}\bigl(
      (2G-rG_Y)(\cA G)_x+rG_x(\cA G)_Y
                         +\partial_x((s^\e)^2)\bigr),\\
 &u_ss_x^\e+w_ss_Y^\e-\frac{2w_s}{r}s^\e-\e\cA s^\e+(2G-rG_Y)(s_s)_x
                  +\bigl(r(s_s)_Y-2s_s\bigr)G_x\\
 &\qquad=-\e^{-\frac72}\mathfrak R_S
 -\e^{\frac72}\bigl(
       (2G-rG_Y)s_x^\e+rG_xs_Y^\e-2G_xs^\e\bigr),\\
 &G|_{x=0}=G_{xx}|_{x=0}=s^\e|_{x=0}=0,\\
 &G_x|_{x=L}=G_{xxx}|_{x=L}=s_x^\e|_{x=L}=0,\\
 &G|_{Y=0}=G_Y|_{Y=0}=s^\e|_{Y=0}=0,
\end{aligned}
\right.
\end{equation}
where $\cA=\partial_x^2+\partial_Y^2-\frac3r\partial_Y$,
$\beta_s=\frac{w_{s,x}-u_{s,Y}}{r}$,
\[
 R[q]
 =\partial_x(u_s^2q_x)
   +\frac1{r^3}\partial_Y(r^3u_s^2q_Y)
 =u_s\cA G-G\cA u_s,
\]
\begin{equation}\label{eq:J-def}
\begin{aligned}
 &\cJ[G]={}-u_{s,x}\cA G+G_x\cA u_s
            +G(\cA u_s)_x+w_s(\cA G)_Y+(2G-rG_Y)(\beta_s)_x+rG_x(\beta_s)_Y,
\end{aligned}
\end{equation}
and
\[
\begin{aligned}
 &\mathfrak R_\beta
 =\frac{\partial_xE_2-\partial_YE_1}{r},\\
 &\mathfrak R_S
 =u_s(s_s)_x+w_s(s_s)_Y-\frac{2w_s}{r}s_s-\e\cA s_s.
\end{aligned}
\]

Since $u^\e=2G-rG_Y$ and $w^\e=rG_x$, the conditions
$G=G_{xx}=0$ at $x=0$ give $u^\e=w_x^\e=0$ there,
and $G_x=G_{xxx}=0$ at $x=L$ give
$w^\e=w_{xx}^\e=0$ there.  Together with the scalar boundary
conditions, these give \eqref{eq:sw-physical-end-data} exactly.
The wall conditions $G=G_Y=s^\e=0$ give
$u^\e=w^\e=s^\e=0$ at $Y=0$.

Section~\ref{sec:mixed-quotient-derivation} justifies the division by
$r^2$ at the axis for the principal stream-function problem.
Lemma~\ref{lem:main-principal-solvability} constructs $G$ directly.
The subsequent weighted estimates hold throughout the cylinder.
For related stream-function regularity, see
\cite{liu2009characterization,nowakowski2022weighted}.

\begin{remark}[The equations at the axis]\label{rem:solution-axis}
Radial Hardy, \eqref{equiv-norm},
\eqref{eq:G-Linfty-H4eps}--\eqref{eq:grad-G-Linfty-H4eps},
and the swirl estimate in Section~\ref{sec:nonlinear} give
$\mathbf U\in H^1(\mathcal C_L)\cap L^\infty(\mathcal C_L)$.
Equation~\eqref{reminder uw} gives $\nabla p^\e w\in L^2(\Omega)$.
Normalize $p^\e$ to have zero mean on the annular cylinder
$\frac{1}{2}<r<1$. Radial Hardy then gives $p^\e\in L^2(\mathcal C_L)$,
hence $P\in L^2(\mathcal C_L)$.

For $\boldsymbol\phi\in C_c^\infty(\mathcal C_L;\mathbb R^3)$,
test the momentum equation away from the axis with
$\bigl(1-\chi(\frac r\eta)\bigr)\boldsymbol\phi$,
where $0<\eta<1$. The cutoff derivative is supported in
$\frac\eta4<r<\frac\eta2$, and
\[
 \left\|\nabla\chi\left(\frac r\eta\right)
       \right\|_{L^2(\mathcal C_L)}^2
 \lesssim\frac{L}{\eta^2}\int_{\frac\eta4}^{\frac\eta2}r\,dr
 \lesssim L.
\]
The viscous and pressure cutoff errors vanish by Cauchy--Schwarz,
since the $L^2$ norms of $\nabla\mathbf U$ and $P$ on
$\{r<\frac\eta2\}$ tend to zero as $\eta\to0$. The convection term converges
directly because $\mathbf U\cdot\nabla\mathbf U\in L^2(\mathcal C_L)$.
The divergence equation holds almost everywhere, since
$\mathbf U\in H^1(\mathcal C_L)$ and the axis has measure zero.
Thus \eqref{eq:NS-axisymmetric-zr} holds in
$\mathcal D'(\mathcal C_L)$.
\end{remark}
In comparison with \cite{guo2023validity}, the principal new
analytic difficulty lies in the possible singular behavior at
the axis $r=0$. We devote Section~\ref{sec:axis} to investigating
such a singularity via Hardy-type inequalities and five-dimensional
elliptic regularity. In contrast to the planar problem, integration by parts in the
remainder estimates produces boundary terms at the axis.
The axis is not a physical boundary, and we prescribe no independent
data there. To eliminate these terms, we use a weight that vanishes
at the axis together with bounds on the growth of $G$ and its
derivatives. The divergence form
$\cL=r^{-3}\partial_Y(r^3\partial_Y)$ gives the natural choice
$w^2=r^3$. The following proposition provides the required growth
estimates.

\begin{proposition}[Near-axis estimates for $G$]
\label{prop:near-axis-H4-blowup}
Let $G\in\mathcal H^4(w)$ satisfy $\|G\|_{H^4(w)}\le M$.
For $\frac12<Y<1$, we have
\begin{equation}\label{eq:near-axis}
\begin{aligned}
 &\|G_{xY}(Y)\|_2\lesssim M,\\
 &\|\cL G(Y)\|_2+\|\nabla_{x,Y}^2G(Y)\|_2
       +\left\|\frac{G_Y(Y)}r\right\|_2
       \lesssim M\Lambda(r),\\
 &\|\nabla_{x,Y}\cL G(Y)\|_2
       +\|\nabla_{x,Y}^3G(Y)\|_2
       \lesssim Mr^{-1},
\end{aligned}
\end{equation}
where $\|F(Y)\|_2:=\|F(\cdot,Y)\|_{L^2(0,L)}$.
\end{proposition}
With the weight $r^3$, the boundary terms at radius $r$
are bounded by $Cr\|G\|_{H^4(w)}^2$ and vanish as $r\to0$.

We study \eqref{eq:sw-coupled-remainder} through two
auxiliary problems. For the main equation we impose the mixed stream-function boundary
conditions of \cite{guo2026prandtl} and consider
\begin{equation}\label{eq:main-equation}
\left\{
\begin{aligned}
 &\partial_xR[q]-\e\cA^2G+\cJ[G]=F_G,\\
 &G|_{x=0}=G_{xx}|_{x=0}=0,\\
 &G_x|_{x=L}=G_{xxx}|_{x=L}=0,\\
 &G|_{Y=0}=G_Y|_{Y=0}=0,
\end{aligned}
\right.
\end{equation}
where
\[
\begin{aligned}
 &F_G={}-\e^{-\frac72}\mathfrak R_\beta
          +F_G^{\mathrm L}[s^\e]+F_G^{\mathrm{NL}}[G,s^\e],\\
 &F_G^{\mathrm L}[s^\e]={}-2\partial_x(s_s s^\e),\\
 &F_G^{\mathrm{NL}}[G,s^\e]={}-\e^{\frac72}\bigl(
      (2G-rG_Y)(\cA G)_x+rG_x(\cA G)_Y
                         +\partial_x((s^\e)^2)\bigr).
\end{aligned}
\]

The swirl equation is
\begin{equation}\label{eq:swirl-equation}
\left\{
\begin{aligned}
 &u_ss_x^\e+w_ss_Y^\e-\frac{2w_s}{r}s^\e-\e\cA s^\e
       =F_s,\\
 &s^\e|_{x=0}=0,\qquad s_x^\e|_{x=L}=0,\\
 &s^\e|_{Y=0}=0,
\end{aligned}
\right.
\end{equation}
where
\[
\begin{aligned}
 &F_s={}-\e^{-\frac72}\mathfrak R_S
          +F_s^{\mathrm L}[G]+F_s^{\mathrm{NL}}[G,s^\e],\\
 &F_s^{\mathrm L}[G]={}-(2G-rG_Y)(s_s)_x
          -\bigl(r(s_s)_Y-2s_s\bigr)G_x,\\
 &F_s^{\mathrm{NL}}[G,s^\e]={}-\e^{\frac72}\bigl(
      (2G-rG_Y)s_x^\e+rG_xs_Y^\e-2G_xs^\e\bigr).
\end{aligned}
\]

Define the spaces
\[
\begin{aligned}
 &\mathcal H_0^4(w)=\{G\in\mathcal H^4(w):
          G\text{ satisfies the boundary conditions in }\eqref{eq:main-equation}\},\\
 &\mathscr V_s^2=\{h\in\mathcal H^2(w):
          h|_{x=0}=h_x|_{x=L}=h|_{Y=0}=0\},
\end{aligned}
\]

For $G\in\mathcal H_0^4(w)$, set
\[
 \mathcal H_4(G)
 :=\e\|\{G_{xx},\sqrt2G_{xY},\cL G\}w\|+\e^2\|\{G_{xxxx},G_{xxxY},\cL G_{xx},
                   (\cL G_x)_Y,\cL^2G\}w\|.
\]

We equip $\mathscr X=\mathcal H_0^4(w)$ with the norm
\begin{equation*}
\begin{aligned}
&\|G\|_{\mathscr X}:=\cD[q],\\
&\cD [q]
:=\cD_2[q] +\cD_3[q]+ \cD_4[G]\\
&\cD_2[q]
:=
        \sqrt\e\Big(
        \left\|u_sq_{xx}w\right \| +\left\|u_sq_{xY}w\right \| 
\Big)\\
        &\cD_3[q]:=\cD_{3, x}[q]+\cD_{3,Y}[q]+|q|_{\p,2},\\
        &\cD_{3, x}[q]
        :=
        \e \Big(
        \left\|\sqrt{u_s}q_{xxx}w\right \| 
        +\left\|\sqrt{u_s}q_{xxY}w \right\| 
        +\left\|\sqrt{u_s}\cL q_xw  \right\| 
        \Big)\\
                &\cD_{3,Y}[q]
        :=\e^\frac32\left\|\sqrt{u_s}(\cL q)_Yw\right\|\\
                        &\cD_{4}[G]
        :=
        \e^2\Big(
        \left\|G_{xxxx}w  \right\| 
        +\left\|G_{xxxY}w  \right\| 
        +\left\|\cL G_{xx}w  \right\|+\left\|\sqrt{u_s}(\cL G_{x})_Y w \right\| 
        +\left\|\cL^2 G w \right\|
        \Big)+| G |_{\p, 3} \\
        &|q|_{\p,2}:=\|\sqrt\e u_s q_{xY}w\|_{x=0}
        +\|\e^{\frac{3}{4}}q_{xY}w\|_{Y=0}
        +\|\sqrt\e u_s q_{xx}w\|_{x=L}
        +\|\e^{\frac{5}{4}}\cL q\,w\|_{Y=0},\\
        &|G|_{\p,3}:=\|\e^{\frac{3}{2}}\sqrt{u_s}G_{xxx}w\|_{x=0}
        +\|\e^{\frac{3}{2}}u_s\cL G_xw\|_{x=0}
        +\|\e^{\frac{3}{2}}u_sG_{xxY}w\|_{x=L}.
\end{aligned}
\end{equation*}

For the swirl remainder, define
\begin{equation}\label{eq:sw-scalar-norms}
\begin{aligned}
 &\mathcal S_\e[s^\e]^2
 =\e\|\sqrt{u_s}s_x^\e w\|^2+\e^2\|\nabla s^\e\,w\|^2+\e^3\bigl(\|s_{xx}^\e w\|^2+2\|s_{xY}^\e w\|^2+\|\cL s^\e\,w\|^2\bigr) +|s^\e|_{\p,1}^2,\\
 &|s^\e|_{\p,1}^2
 :=\e^2\|s_x^\e w\|_{x=0}^2+\e^2\|s_Y^\e w\|_{x=L}^2
       +\e\|\sqrt{u_s}s^\e w\|_{x=L}^2.
\end{aligned}
\end{equation}
Equip $\mathscr Z_\e=\mathscr X\times\mathscr V_s^2$ with the norm
$\|(G,s^\e)\|_{\mathscr Z_\e}=\cD[q]+\mathcal S_\e[s^\e]$.

In Section~\ref{sec:profiles-construction}, we construct the Euler
and raw Prandtl profiles from the leading entrance data. At each inner order, we solve the meridional equations before the
scalar swirl equation. The radial momentum equation determines the
inner pressure. We truncate the meridional velocity profiles through
their stream function, preserving continuity and no slip.
Theorem~\ref{thm:approximate-profiles} gives the bounds
$u_s\sim \frac{Y}{\sqrt\e}$ for $Y\le\sqrt\e$ and $u_s\sim1$ otherwise,
together with the derivative and quotient estimates
\eqref{eq:usbound}--\eqref{eq:profile-quotient-bounds} used in the
remainder equations. Using the profile equations to cancel like powers of $\sqrt\e$ and
the decay of the inner profiles, we obtain, for $N\ge8$,
\[
 \sqrt\e\|\mathfrak R_\beta w\|+\|\mathfrak R_Sw\|
 \lesssim\e^{\frac{9}{2}}.
\]

In Sections~\ref{sec:preliminary} and~\ref{sec:axis}, we prove
the preliminary inequalities and the radial regularity and recovery
estimates.
The recovery estimates control $\mathcal H_4(G)$ by $\cD[q]$,
and the five-dimensional Agmon inequalities then give
\[
 \|G\|_\infty\lesssim_L\e^{-\frac{5}{4}}\cD[q],\qquad
 \|\nabla G\|_\infty\lesssim_L\e^{-\frac{7}{4}}\cD[q].
\]
By \eqref{eq:energy-norm-equivalence},
$\mathscr Z_\e$ is complete for fixed $\e,L$.

In Section~\ref{sec:coupled-estimates}, the estimates
\eqref{qxx}, \eqref{Lq-est}, \eqref{qx-est}, \eqref{trace 1}, and
\eqref{trace 2}, together with the main equation, give
\[
\begin{aligned}
 &\cD_3^2[q]
 \lesssim L^{-\frac{3}{8}}\cD_2^2[q]
       +\sqrt L\cD_4^2[G]+\e L^{\frac{1}{8}}\|F_Gw\|^2,\\
 &\cD_2^2[q]
 \lesssim\sqrt L\bigl(\cD_3^2[q]+\cD_4^2[G]\bigr)
       +\e L^{\frac{3}{2}}\|F_Gw\|^2,\\
 &\cD_4^2[G]
 \lesssim\cD_2^2[q]+\cD_3^2[q]+\e^2\|F_Gw\|^2.
\end{aligned}
\]
Combining these estimates yields
\[
 L^{-\frac{3}{8}}\cD_2^2[q]+\cD_3^2[q]
 \lesssim\e L^{\frac{1}{8}}\|F_Gw\|^2.
\]
Our main and scalar swirl estimates are
\[
 \cD[q]\lesssim\sqrt\e L^{\frac{1}{16}}\|F_Gw\|,\qquad\mathcal S_\e[s^\e]\lesssim\|F_sw\|.
\]
For the coupled equations, with the nonlinear sources prescribed,
\eqref{eq:main-feedback} and \eqref{eq:sw-scalar-feedback} give
\[
\begin{aligned}
 &\cD[q]
 \lesssim L^{\frac{1}{16}}\Bigl(
       \e^{-3}\|\mathfrak R_\beta w\|
       +\sqrt\e\|F_G^{\mathrm{NL}}w\|
       +\mathcal S_\e[s^\e]\Bigr),\\
 &\mathcal S_\e[s^\e]
 \lesssim \e^{-\frac{7}{2}}\|\mathfrak R_Sw\|
       +\|F_s^{\mathrm{NL}}w\|+\cD[q].
\end{aligned}
\]
The factor $L^{\frac{1}{16}}$ allows us to absorb the coupling and obtain
\eqref{eq:sw-linear-final}.

In Section~\ref{sec:nonlinear}, we estimate the nonlinear terms by
\[
\begin{aligned}
 &\sqrt\e\|F_G^{\mathrm{NL}}[G,s^\e]w\|
 \lesssim_{L,K}\e^{\frac{3}{4}}\cD[q]^2
                       +\e\mathcal S_\e[s^\e]^2,\\
 &\|F_s^{\mathrm{NL}}[G,s^\e]w\|
 \lesssim_{L,K}\e^{\frac{3}{4}}\cD[q]\mathcal S_\e[s^\e].
\end{aligned}
\]
The profile residual bound and the coupled estimate yield
\[
 \cD[q]+\mathcal S_\e[s^\e]
 \lesssim_{L,K}\e
       +\e^{\frac{3}{4}}\bigl(\cD[q]+\mathcal S_\e[s^\e]\bigr)^2.
\]
The corresponding difference estimates give a contraction on the
closed unit ball of $\mathscr Z_\e$. At the fixed point,
$\cD[q]+\mathcal S_\e[s^\e]\lesssim_{L,K}\e$, and the velocity
recovery estimates give \eqref{M2}.

\subsection{Other works}\label{subsec:other-works}

Guo and Nguyen \cite{guo2017moving} justified steady Prandtl
expansions over a moving plate on a sufficiently short tangential
interval. The motion of the plate enters their positivity estimate
through the nonzero tangential velocity at the boundary.
Iyer \cite{Iyer2019steady} extended the moving-wall analysis to
nonshear outer Euler flows. For a sufficiently small velocity
mismatch between the plate and the outer flow, the works
\cite{iyer2019globalI,iyer2019globalII,iyer2020globalIII}
established global-in-$x$ expansions.

The classical setting of a stationary wall was treated by Guo
and Iyer \cite{guo2023validity}. Their result justifies the
expansion on a short tangential interval for a class of scaled
Prandtl profiles, including the classical Blasius profiles.
Their proof combines a quotient formulation of the differentiated
vorticity equation with estimates for the entrance traces.
See also Gao and Zhang \cite{gao2023steady} for nonshear
outer Euler flows. In the tangentially periodic half-plane,
G\'erard-Varet and Maekawa \cite{gerardvaret2019sobolev} obtained
Sobolev stability of Prandtl-type shear flows for sufficiently
small periods. Their hypotheses require positivity away from
the wall and positive wall shear, but impose neither global
monotonicity nor concavity.

Global-in-$x$ validity at a stationary wall was established by
Iyer and Masmoudi \cite{iyermasmoudi2026global} for boundary-layer
data near Blasius. Their analysis uses the Euler--Prandtl profiles
constructed in \cite{iyermasmoudi2021profiles}, with decay estimates
for both the profiles and the Navier--Stokes remainder.
For an outer Euler shear flow in a bounded channel, Guo and Yang
\cite{guo2024small} constructed steady Navier--Stokes solutions
with an additional $O(\e)$ streamwise scale in the transverse
velocity near a channel end. In their construction, this scale
coexists with the classical $O(\sqrt\e)$ Prandtl layers at the
lateral walls.

For the steady Prandtl equation itself, Oleinik initiated the
existence theory presented in \cite{oleinik1999mathematical}.
For a class of compatible entrance profiles that are positive away
from the wall and have positive wall shear, this theory gives
local-in-$x$ solutions, which extend globally under a favorable
pressure gradient. The higher regularity needed to construct
successive terms in the expansion was established by Guo and Iyer
\cite{guo2021regularity} in the absence of an external pressure
gradient. Their estimates for a differentiated quotient equation
retain the decay of the entrance data and permit the construction
of Euler--Prandtl expansions to arbitrary finite order.
Under a favorable pressure gradient, Wang and Zhang
\cite{wangzhang2021regularity} proved smoothness of Oleinik's
solutions up to the wall at every positive downstream distance.

For zero pressure gradient, the downstream convergence of steady
Prandtl solutions to Blasius has been studied under several
assumptions on the entrance data. Iyer \cite{iyer2020blasius}
proved asymptotic stability for small localized perturbations,
and Wang and Zhang \cite{wangzhang2023asymptotic} treated more
general entrance data with exponential decay. Sharp convergence
rates were obtained by Jia, Lei and Yuan \cite{jia2025sharp}.
Under an adverse pressure gradient, Dalibard and Masmoudi
\cite{dalibard2019separation} proved boundary-layer separation
and the Goldstein singularity for a class of entrance data.

Beyond flat boundaries, Iyer \cite{iyer2017steady} justified
steady Prandtl expansions over a rotating disk. Boundary layers
in circular pipes have also been studied for time-dependent
flows with symmetry. Mazzucato and Taylor
\cite{mazzucato2011pipe} considered the vanishing viscosity
limit for circular pipe flows independent of the axial variable,
including estimates for their boundary layers. Gie et al.\
\cite{gie2019symmetric} constructed boundary-layer correctors
and obtained convergence rates for symmetric channel and pipe
flows, with annular cross sections in the pipe case.

For the time-dependent Prandtl equations, the distinction between
monotone and nonmonotone profiles appears in the results of Guo
and Nguyen \cite[Theorems~1.2, 1.4 and 5.2]{guo2011note}.
They proved weak well-posedness in Oleinik's monotone class,
whereas suitable nonmonotone shear flows lead to failure of
Lipschitz Sobolev well-posedness and of uniform Sobolev estimates
for the Navier--Stokes expansion.

\section{Construction of the profiles}
\label{sec:profiles-construction}

We construct the profiles in \eqref{eq:sw-remainder-ansatz}.
For the prescribed $N$, choose finite derivative and weight orders
$Q$ sufficient for the remainder estimates and take $M=M(N,Q)$
sufficiently large. Fix a bound for $\|H\|_{C^{4M}([0,1])}$.
The higher Prandtl entrance data and Euler end traces are chosen
compatibly on $[0,L]$.
For the velocity profiles, we follow \cite[Section~2]{gao2023steady}
and \cite[Appendix~A.3]{Iyer2017nonshear}. At each inner order,
we then solve the scalar swirl equation. Truncating the stream function gives
\eqref{eq:profile-composite}.

\allowdisplaybreaks[1]

\subsection{The leading Prandtl profiles}
\label{subsec:swirl-leading}
Prescribe smooth leading Prandtl entrance velocities $U_{\mathrm{in}}$ and
$S_{\mathrm{in}}$ for the axial and normalized swirl components, with
\begin{equation}\label{eq:sw-profile-input-data}
\begin{aligned}
 &U_{\mathrm{in}}(0)=0,\qquad U_{\mathrm{in}}(y)>0\quad(y>0),
 \qquad U_{\mathrm{in}}'(0)>0,\qquad
 U_{\mathrm{in}}(y)\to1,\\
 &S_{\mathrm{in}}(0)=0,\qquad S_{\mathrm{in}}(y)\to\kappa,\\
 &\|e^{\sigma_0y}\partial_y^j
       (U_{\mathrm{in}}-1,S_{\mathrm{in}}-\kappa)\|_{L_y^\infty}
 \le C_j
 \quad\text{for $0\le j\le4M$}.
\end{aligned}
\end{equation}
Here $\sigma_0>0$ is fixed.
The entrance data are also assumed to satisfy the corner compatibility
conditions through axial-jet order $M$. For the axial and normal
velocities, these are the conditions of
\cite[Definition~2.1]{guo2021regularity}. For the swirl component, compute the axial entrance derivatives of
orders $1,\ldots,M$ recursively from the third equation in
\eqref{eq:original-method-leading-Prandtl}. Each must extend smoothly
to $y=0$ and vanish there.

For example, one may choose $U_{\mathrm{in}}$ so that the
leading axial velocity in the boundary layer is a downstream translate
of the classical Blasius profile.
The swirl entrance profile need not be proportional to the axial
profile and may be prescribed independently, subject to these
decay and compatibility conditions.

For the leading axial and normal
profiles, we write $\bar u=1+\widetilde u_p^0$ and
$\bar w=\widetilde w_p^0-\widetilde w_p^0(x,0)$.

The leading swirl matches $\kappa=K(1)$. At the wall, the axial
equation and its $y$ derivative require
$U_{\mathrm{in}}''(0)=U_{\mathrm{in}}'''(0)=0$.
The leading profiles solve
\begin{equation}\label{eq:original-method-leading-Prandtl}
\left\{
\begin{aligned}
 &\bar u\widetilde u_{p,x}^0+
       \bar w\widetilde u_{p,y}^0-\widetilde u_{p,yy}^0=0,\\
 &\widetilde u_{p,x}^0+\widetilde w_{p,y}^0=0,\\
 &(\bar u\partial_x+\bar w\partial_y-\partial_y^2)
                                      \widetilde S_p^0=0,\\
 &(\widetilde u_p^0,\widetilde S_p^0)(0,y)
       =\bigl(U_{\mathrm{in}}(y)-1,
                      S_{\mathrm{in}}(y)-\kappa\bigr),\\
 &(\widetilde u_p^0,\widetilde S_p^0)(x,0)=(-1,-\kappa),\\
 &\lim_{y\to\infty}(\widetilde u_p^0,\widetilde w_p^0,
                    \widetilde S_p^0)(x,y)=(0,0,0).
\end{aligned}
\right.
\end{equation}
The normal velocity is recovered by
\begin{equation}\label{eq:original-method-leading-tail}
 \widetilde w_p^0(x,y)
       =\int_y^\infty\widetilde u_{p,x}^0(x,z)\,dz.
\end{equation}

\begin{proposition}[Leading Prandtl profiles]
\label{thm:leading-prandtl-one-wall}
The problem \eqref{eq:original-method-leading-Prandtl} has a unique
solution satisfying the bounds below, with velocity in the Oleinik class.
For all sufficiently small $L$, uniformly on $[0,L]$,
$\bar u\sim\min\{y,1\}$ and $\bar u_y(x,0)\gtrsim1$.
For the orders allowed by the data,
\begin{equation}\label{eq:sw-leading-profile-bounds}
 \|e^{\sigma y}\partial_x^a\partial_y^b
       (\widetilde u_p^0,\widetilde w_p^0,\widetilde S_p^0)
                                      \|_{L^\infty(\Omega_p)}\lesssim1.
\end{equation}
Here $0<\sigma<\sigma_0$ is fixed.
\end{proposition}

\begin{proof}
Existence, positivity, and the finite derivative bounds for the velocity
follow from Oleinik's theorem and
\cite[Theorem~1.2]{guo2021regularity}, with the decaying normal profile
given by \eqref{eq:original-method-leading-tail}.
The exponential-weight estimates of
\cite[Theorems~1.2 and 1.3]{guo2021regularity} apply at any fixed
positive rate after the dilation $x=a^2X$, $y=az$, with the normal
velocity multiplied by $a$. This preserves the Prandtl equations
and their linearization, with rescaled forcing, and replaces $\sigma$
by $a\sigma$. Choose a fixed $a$ sufficiently large and scale back.
These bounds give
$|\partial_x^k\bar u|\lesssim\bar u$ and
$|\partial_x^k\bar w|\lesssim\min\{y^2,1\}$.

For the scalar equation, let
\[
 (\bar u\partial_x+\bar w\partial_y-\partial_y^2)v=F,
 \qquad v(x,0)=v(x,\infty)=0.
\]
A wall cutoff and integration by parts give, for $\eta>0$,
\[
\begin{aligned}
 &\|e^{\sigma y}v\|_{L_y^2}^2
 \le\eta\|e^{\sigma y}v_y\|_{L_y^2}^2
       +C_\eta\|\sqrt{\bar u}\,e^{\sigma y}v\|_{L_y^2}^2,\\
 &\|e^{\sigma y}v_y\|_{L_y^2}^2
 \le\eta\|e^{\sigma y}v_{yy}\|_{L_y^2}^2
       +C_\eta\|\sqrt{\bar u}\,e^{\sigma y}v\|_{L_y^2}^2.
\end{aligned}
\]
Apply $\partial_x^j$ and test with $e^{2\sigma y}\partial_x^jv$.
Using $\bar u_x+\bar w_y=0$, the preceding bounds, and Gronwall's
inequality, we obtain by induction on $j$
\begin{equation}\label{eq:sw-Prandtl-scalar-energy}
\begin{aligned}
 &\sum_{j=0}^J\bigg[
  \sup_{0\le x\le L}\|\sqrt{\bar u}\,e^{\sigma y}
                  \partial_x^jv(x)\|_{L_y^2}^2
  +\|e^{\sigma y}\partial_y\partial_x^jv\|_{L^2(\Omega_p)}^2\bigg]\\
 &\lesssim\sum_{j=0}^J\bigg[
  \|\sqrt{U_{\mathrm{in}}}\,e^{\sigma y}
                  \partial_x^jv(0)\|_{L_y^2}^2
  +\|e^{\sigma y}\partial_x^jF\|_{L^2(\Omega_p)}^2\bigg].
\end{aligned}
\end{equation}
One further axial derivative and the fundamental theorem of calculus give
\[
 \|e^{\sigma y}\partial_y\partial_x^jv(x)\|_2^2
 \le\|e^{\sigma y}\partial_y\partial_x^jv(0)\|_2^2
 +2\|e^{\sigma y}\partial_y\partial_x^jv\|_{L^2(\Omega_p)}
    \|e^{\sigma y}\partial_y\partial_x^{j+1}v\|_{L^2(\Omega_p)}.
\]
The equation gives $v_{yy}$ and higher normal derivatives, and
Sobolev embedding in $y$ yields the mixed bounds, uniformly for
$0<L\le L_0$. The argument also applies to polynomial weights.
For existence, subtract the compatible axial Taylor polynomial $T$
through an order higher than those to be estimated. The remaining
forcing $F-(\bar u\partial_x+\bar w\partial_y-\partial_y^2)T$
has vanishing entrance jets to the required order. In this zero-data
problem replace $\bar u$ by $\bar u+\rho$, $\rho>0$, and solve by
weighted Galerkin approximation. The preceding energy estimate,
with $\sqrt{\bar u+\rho}$ in place of $\sqrt{\bar u}$, is uniform in
$\rho$. The differentiated entrance data are zero. Weak compactness
as $\rho\to0$, followed by the normal derivative recovery above,
gives the solution.
Estimate~\eqref{eq:sw-Prandtl-scalar-energy} gives uniqueness.
For wall data $h(x)$, subtract $h(x)\varphi(y)$ from the unknown
and its image under the operator from $F$.

Applying these estimates to
$v=\widetilde S_p^0-\kappa\widetilde u_p^0$, $F=0$, with entrance
value $S_{\mathrm{in}}-\kappa U_{\mathrm{in}}$ proves the scalar
bounds for the prescribed $S_{\mathrm{in}}$.
\end{proof}

\subsection{Coefficient equations}
\label{subsec:sw-coefficient-equations}

At each higher order we impose
\begin{equation}\label{eq:profile-matching}
 \widetilde u_p^i(x,0)=-u_e^i(x,0),\qquad
 \widetilde S_p^i(x,0)=-S_e^i(x,0),\qquad
 w_e^i(x,0)=-\widetilde w_p^{i-1}(x,0).
\end{equation}
The inner profiles vanish at infinity. We match powers of
$\sqrt\e$ before truncation and before adding the remainders.
After subtracting the outer equations from \eqref{eq:NS-xY} and
writing the swirl equation as $US_x+WS_Y-\frac{2WS}{r}=\e\cA S$, we expand
the remaining outer factors at $Y=0$. For example,
\[
 \frac1{1-\sqrt{\e}\,y}
   \sim\sum_{j\geq0}\e^{\frac{j}{2}}y^j,\qquad
 u_e^i(x,\sqrt{\e}\,y)
   \sim\sum_{j\geq0}\frac{\e^{\frac{j}{2}}y^j}{j!}
                   \partial_Y^ju_e^i(x,0).
\]
At each order, continuity gives
\begin{equation}\label{eq:inner-continuity-source}
 \widetilde u_{p,x}^i+\widetilde w_{p,y}^i=h_p^i,\qquad
 h_p^i=\sum_{j=0}^{i-1}y^{i-1-j}\widetilde w_p^j .
\end{equation}

For $i\ge1$, the coefficient of $\e^{\frac{i-1}{2}}$ in the radial
equation determines $\partial_y\widetilde p_p^i$, with
$\widetilde p_p^i(x,\infty)=0$. The right-hand side involves outer
profiles through $i-1$, inner meridional profiles through $i-2$,
and inner swirl profiles through $i-1$. Here,
\[
 \partial_y\widetilde p_p^1=-\bigl[(\kappa+\widetilde S_p^0)^2-\kappa^2\bigr].
\]
The coefficients of $\e^{\frac{i}{2}}$ in the axial and swirl equations
determine $(\widetilde u_p^i,\widetilde w_p^i)$ and $\widetilde S_p^i$,
respectively, after replacing $w_e^{i+1}(x,0)$ by
$-\widetilde w_p^i(x,0)$. We construct the profiles in the order
\[
 (u_e^i,w_e^i,S_e^i,p_e^i),\qquad
 \widetilde p_p^i,\qquad (\widetilde u_p^i,\widetilde w_p^i),\qquad \widetilde S_p^i .
\]
The sources, including the centrifugal and inner pressure terms,
are sums of decaying inner factors multiplied by polynomials in $y$
and known outer wall derivatives.

\subsection{Euler correctors}
\label{subsec:swirl-Euler}
Given the preceding profiles, we solve
\begin{equation}\label{eq:sw-Euler-i}
\left\{\begin{aligned}
 &(u_e^i)_x+(p_e^i)_x=F_{e,1}^i,\\
 &(w_e^i)_x+(p_e^i)_Y+2rKS_e^i=F_{e,2}^i,\\
 &(S_e^i)_x-\frac{2K+rK'}r w_e^i=F_{e,S}^i,\\
 &(u_e^i)_x+(w_e^i)_Y-\frac{w_e^i}{r}=0,\\
 &w_e^i(x,0)=-\widetilde w_p^{i-1}(x,0),
\end{aligned}\right.
\end{equation}
where $K'=\partial_rK$ and
\begin{equation}\label{eq:sw-Euler-forcing}
\begin{aligned}
 &F_{e,1}^i
 =\left(\partial_x^2+\partial_Y^2-\frac1r\partial_Y\right)u_e^{i-2}
   -\sum_{\substack{j+k=i\\j,k\ge1}}
        \bigl(u_e^j(u_e^k)_x+w_e^j(u_e^k)_Y\bigr),\\
 &F_{e,2}^i
 =\left(\partial_x^2+\partial_Y^2-\frac1r\partial_Y
                                  -\frac1{r^2}\right)w_e^{i-2}
   -\sum_{\substack{j+k=i\\j,k\ge1}}
        \bigl(u_e^j(w_e^k)_x+w_e^j(w_e^k)_Y+rS_e^jS_e^k\bigr),\\
 &F_{e,S}^i
 =\cA S_e^{i-2}
   -\sum_{\substack{j+k=i\\j,k\ge1}}
        \left(u_e^j(S_e^k)_x+w_e^j(S_e^k)_Y
                                      -\frac{2w_e^j}{r}S_e^k\right).
\end{aligned}
\end{equation}
Here $\cA K=8H'(r^2)+4r^2H''(r^2)$.
Set
\[
 u_e^i=2\Gamma_e^i-r(\Gamma_e^i)_Y,\qquad
 w_e^i=r(\Gamma_e^i)_x.
\]
Taking the curl and using the swirl equation gives
\begin{equation}\label{eq:sw-Euler-derivative}
 \bigl(\cA+2K(2K+rK')\bigr)\partial_x\Gamma_e^i
 =\frac{\partial_xF_{e,2}^i-\partial_YF_{e,1}^i}{r}
                                      -2K F_{e,S}^i.
\end{equation}
Choose the integration data by
\begin{equation}\label{eq:sw-Euler-potential}
\begin{aligned}
 &\bigl(\cA+2K(2K+rK')\bigr)\Gamma_e^i
 =\int_0^x
 \left(\frac{\partial_xF_{e,2}^i-\partial_YF_{e,1}^i}{r}
                            -2KF_{e,S}^i\right)(t,Y)\,dt,\\
 &\Gamma_e^i(x,0)=-\int_0^x\widetilde w_p^{i-1}(t,0)\,dt,
\end{aligned}
\end{equation}
and recover
\[
 S_e^i=(2K+rK')\Gamma_e^i+\int_0^xF_{e,S}^i(t,Y)\,dt.
\]
With $\Gamma_e^0=\frac{1}{2}$, continuity gives
\[
 \frac{\partial_xF_{e,2}^i-\partial_YF_{e,1}^i}{r}
 =\cA^2\Gamma_e^{i-2}
 -\sum_{\substack{j+k=i\\j,k\ge1}}
 \left[u_e^j\partial_x(\cA\Gamma_e^k)
       +w_e^j\partial_Y(\cA\Gamma_e^k)
       +\partial_x(S_e^jS_e^k)\right].
\]
The forcing has a regular radial lift, since for a regular radial scalar $F$,
\[
 (w_e^j\partial_YF)^\sharp
 =-(\partial_x\Gamma_e^j)^\sharp\xi\cdot\nabla_\xi F^\sharp.
\]

\begin{proposition}[Euler profiles]\label{prop:sw-Euler-finite-cylinder}
Let $0<L\le1$ satisfy
\begin{equation}\label{eq:sw-Euler-small-length}
 \frac{L^2}{\pi^2}
       \bigl\|[2K(2K+rK')]_+\bigr\|_\infty\le\frac12.
\end{equation}
For $m\ge0$, take $M=3m+6$. Assume that the
preceding $\Gamma_e^j,S_e^j$ have regular radial lifts, the lifted
right-hand side of \eqref{eq:sw-Euler-derivative} belongs to
$C_x^MH_\xi^M$, and
$\widetilde w_p^{i-1}(\cdot,0)\in C^M([0,L])$.
There are axial Dirichlet data determined by these quantities,
independent of $\e$, for which \eqref{eq:sw-Euler-potential} has a
unique solution, regular at the axis, satisfying
\begin{equation}\label{eq:sw-Euler-profile-estimate}
\begin{aligned}
 &\|(\Gamma_e^i)^\sharp\|_{C_x^{m+1}H_\xi^{m+2}}
 \lesssim
 \left\|\left(
 \frac{\partial_xF_{e,2}^i-\partial_YF_{e,1}^i}{r}-2KF_{e,S}^i\right)^\sharp\right\|_{C_x^MH_\xi^M}+\|\widetilde w_p^{i-1}(\cdot,0)\|_{C^M([0,L])}.
\end{aligned}
\end{equation}
Moreover,
\begin{equation}\label{eq:sw-Euler-velocity-estimate}
\begin{aligned}
 &\|\bigl(u_e^i,\frac{w_e^i}{r}\bigr)^\sharp\|_{C_x^mH_\xi^{m+1}}
 \lesssim\|(\Gamma_e^i)^\sharp\|_{C_x^{m+1}H_\xi^{m+2}},\\
 &\|(S_e^i)^\sharp\|_{C_x^mH_\xi^{m+1}}
 \lesssim\|(\Gamma_e^i)^\sharp\|_{C_x^{m+1}H_\xi^{m+2}}
             +\|(F_{e,S}^i)^\sharp\|_{C_x^mH_\xi^{m+1}}.
\end{aligned}
\end{equation}
The constants depend on the indicated finite norms and orders,
and are independent of $L,\e$. Normalize $p_e^i(0,1)=0$.
\end{proposition}

\begin{proof}
Choose a wall lifting $B$ with
\[
 B(x,0)=-\int_0^x\widetilde w_p^{i-1}(t,0)\,dt,
 \qquad \partial_Y^{2j+1}B(x,0)=0\quad(0\le j\le m+1).
\]
Determine its even normal jets by
\begin{equation}\label{eq:sw-Euler-lifting-jets}
\begin{aligned}
 &\left.(\cL+2K(2K+rK'))^{j+1}B\right|_{Y=0}\\
 &\quad=\left.(\cL+2K(2K+rK'))^j
       \int_0^x\left(
       \frac{\partial_xF_{e,2}^i-\partial_YF_{e,1}^i}{r}
                              -2KF_{e,S}^i\right)(t,Y)\,dt\right|_{Y=0}\\
 &\qquad\quad-\partial_x^2
       \left.(\cL+2K(2K+rK'))^jB\right|_{Y=0},
       \qquad 0\le j\le m+1.
\end{aligned}
\end{equation}
The coefficient of $\partial_Y^{2j+2}B(x,0)$ is one. Realize these jets by
\[
 B(x,Y)=\chi(Y)
          \sum_{\ell=0}^{2m+4}\frac{Y^\ell}{\ell!}
                                  \partial_Y^\ell B(x,0).
\]
The lifting vanishes near the axis. Its required norms are bounded
by the right-hand side of \eqref{eq:sw-Euler-profile-estimate}, by
finite differentiation, the trace estimates on the fixed ball, and
$\|\int_0^x a(t)\,dt\|_\infty\le L\|a\|_\infty$.
Set
\begin{equation}\label{eq:sw-Euler-end-traces}
 \Gamma_e^i(0,Y)=B(0,Y),\qquad
 \Gamma_e^i(L,Y)=B(L,Y).
\end{equation}
Then $\|B(L,\cdot)-B(0,\cdot)\|_{C^k}
\le L\|B_x\|_{C_x^0C^k}$.

Let $T=\Delta_\xi+2K(2K+rK')$ be the Dirichlet operator on $B_4$.
The lifted function $v=(\Gamma_e^i-B)^\sharp$ satisfies
\[
 (\partial_x^2+T)v=f,
 \qquad v|_{\partial((0,L)\times B_4)}=0.
\]
Here $f$ is the radial lift of the right-hand side of
\eqref{eq:sw-Euler-potential} after subtracting $(\cA+2K(2K+rK'))B$.
Equation~\eqref{eq:sw-Euler-lifting-jets} gives
\[
 T^j\partial_x^qf|_{\partial B_4}=0
 \qquad(0\le j\le m+1)
\]
for every retained $q$. Thus $\partial_x^qf\in D(T^{m+2})$
with uniformly bounded full graph norms.
For $q\le m+1$, these norms use at most $3m+6$ axial derivatives
of the wall datum, and $3m+4$ axial and $2m+4$ radial derivatives
of the right-hand side of \eqref{eq:sw-Euler-derivative}.

For $\lambda\in\operatorname{spec}T$,
$\lambda\le\|[2K(2K+rK')]_+\|_\infty$.
The Dirichlet Green kernel of $\partial_x^2+\lambda$ on $(0,L)$ satisfies
\[
 |G_\lambda(x,t)|\le CL,\qquad
 |\partial_xG_\lambda(x,t)|\le C.
\]
For $\lambda=-a^2<0$,
\[
 G_\lambda(x,t)=
 -\frac{\sinh(a\min\{x,t\})\sinh(a(L-\max\{x,t\}))}
              {a\sinh(aL)},
\]
whose first derivative has absolute value at most one.
For $\lambda=a^2>0$, replace $\sinh$ by $\sin$ and use
$aL\le\frac{\pi}{\sqrt2}$ from \eqref{eq:sw-Euler-small-length}.
The case $\lambda=0$ follows by a limit.
Taking the supremum over the spectrum and applying Minkowski's
inequality, we obtain, for $0\le j\le m+2$,
\begin{equation}\label{eq:sw-Euler-Green-estimate}
 \|T^jv\|_{C_x^0L_\xi^2}\le CL^2\|T^jf\|_{C_x^0L_\xi^2},\qquad\|T^jv_x\|_{C_x^0L_\xi^2}\le CL\|T^jf\|_{C_x^0L_\xi^2}.
\end{equation}
Since $T$ is independent of $x$,
\[
 \partial_x^{q+2}v=\partial_x^qf-T\partial_x^qv.
\]
The inequality
$\lceil\frac{m+2}{2}\rceil+\lfloor\frac{m+1}{2}\rfloor\le m+2$
and Dirichlet elliptic regularity on the fixed ball give
$C_x^{m+1}H_\xi^{m+2}$ bounds. Adding $B$ proves
\eqref{eq:sw-Euler-profile-estimate}.

Existence follows from the Green representation. For uniqueness,
Poincar\'e's inequality in the axial variable and
\eqref{eq:sw-Euler-small-length} give
\[
 \int_{(0,L)\times B_4}
   \bigl(|\partial_xv|^2+|\nabla_\xi v|^2
                        -2K(2K+rK')|v|^2\bigr)
 \ge\tfrac12\|\partial_xv\|_2^2+\|\nabla_\xi v\|_2^2,
\]
which proves uniqueness. Rotational invariance and regularity on
the fixed ball give regularity at the axis. The reconstructed profiles
satisfy \eqref{eq:sw-Euler-velocity-estimate}, continuity, the swirl
equation, and
\[
 \cA\Gamma_e^i+2KS_e^i
   =\int_0^x
       \frac{\partial_xF_{e,2}^i-\partial_YF_{e,1}^i}{r}(t,Y)\,dt.
\]
We recover the pressure from
\[
 (p_e^i)_x=F_{e,1}^i-(u_e^i)_x,\qquad
 (p_e^i)_Y=F_{e,2}^i-(w_e^i)_x-2rKS_e^i.
\]
The preceding identity shows that the $Y$-derivative of the first
right-hand side equals the $x$-derivative of the second. The condition
$p_e^i(0,1)=0$ therefore fixes $p_e^i$ uniquely.
\end{proof}

The entrance values and axial velocity are
\begin{equation}\label{eq:sw-Euler-inlet-recovery}
\begin{aligned}
 &u_e^i(0,Y)=2B(0,Y)-r B_Y(0,Y),\\
 &S_e^i(0,Y)=(2K+rK')B(0,Y),\\
 &u_e^i(x,Y)=u_e^i(0,Y)
       -\int_0^x\left((w_e^i)_Y-\frac{w_e^i}{r}\right)(t,Y)\,dt.
\end{aligned}
\end{equation}

\subsection{The inner profiles}
\label{subsec:inner-profiles}

For $1\le i\le N$, the inner equations are
\begin{equation}\label{eq:inner-fixed-system}
\begin{aligned}
 &\bar u\,\widetilde u_{p,x}^i+\bar u_x\widetilde u_p^i+\bar w\,\widetilde u_{p,y}^i
   +\bar u_y\bigl(\widetilde w_p^i-\widetilde w_p^i(x,0)\bigr)-\widetilde u_{p,yy}^i=f_p^i,\\
 &\widetilde u_{p,x}^i+\widetilde w_{p,y}^i=h_p^i,\\
 &\bar u\,\widetilde S_{p,x}^i+\bar w\,\widetilde S_{p,y}^i-\widetilde S_{p,yy}^i=g_p^i,\\
 &\widetilde u_p^i(x,0)=-u_e^i(x,0),\qquad \widetilde S_p^i(x,0)=-S_e^i(x,0),\\
 &(\widetilde u_p^i,\widetilde S_p^i)(0,y)
                 =(U_{p,0}^i,S_{p,0}^i)(y),\\
 &(\widetilde u_p^i,\widetilde w_p^i,\widetilde S_p^i)(x,y)\to0\qquad(y\to\infty).
\end{aligned}
\end{equation}
Both $f_p^i$ and $h_p^i$ are now known. The pressure contribution
to $f_p^i$ was determined at the preceding step.
Solving the meridional equations determines the remaining terms in
$g_p^i$, including
$-\widetilde u_p^i\widetilde S_{p,x}^0
 -(\widetilde w_p^i-\widetilde w_p^i(x,0))\widetilde S_{p,y}^0$.
The entrance traces satisfy the finite corner conditions below.

\begin{lemma}[Inner profiles]\label{lem:inner-profiles}
Fix finite derivative and weight orders $J,m$. Suppose the known
inner profiles have sufficiently many exponentially weighted
derivatives and the outer profiles have sufficiently many derivatives.
Compatible higher-order entrance traces exist. For each fixed choice
of such traces, \eqref{eq:inner-fixed-system} has a unique solution
on $[0,L]$ satisfying
\begin{equation}\label{eq:inner-uniform-estimate}
 \sum_{a+b\leq J}
 \bigl\|\langle y\rangle^m\partial_x^a\partial_y^b
       (\widetilde u_p^i,\widetilde w_p^i,\widetilde S_p^i)\bigr\|_{L^\infty((0,L)\times\mathbb R_+)}
 \leq C_{i,J,m},\qquad 0<L\leq L_0.
\end{equation}
The profiles and entrance choices are independent of $\e$.
The constant depends only on the finite norms fixed in the
construction and is independent of $\e,L$.
\end{lemma}

\begin{proof}
The first corner conditions for the tangential entrance trace are
\[
\begin{aligned}
 &\widetilde u_p^i(0,0)=-u_e^i(0,0),\qquad 
 \widetilde u_{p,yy}^i(0,0)=-f_p^i(0,0),\\
 &\widetilde u_{p,yyy}^i(0,0)
 =\bar u_y(0,0)h_p^i(0,0)
   -\bar u_{xy}(0,0)u_e^i(0,0)-\partial_yf_p^i(0,0).
\end{aligned}
\]
Momentum and continuity determine the normal entrance trace:
\begin{equation}\label{eq:inner-normal-inlet}
 \frac{\widetilde w_p^i(0,y)-\widetilde w_p^i(0,0)}{\bar u(0,y)}
 =\frac{h_p^i(0,0)+u_{e,x}^i(0,0)}{\bar u_y(0,0)}+\int_0^y
 \frac{\bar u h_p^i+\bar w \widetilde u_{p,y}^i+\bar u_x\widetilde u_p^i
                    -\widetilde u_{p,yy}^i-f_p^i}
      {\bar u^2}(0,z)\,dz .
\end{equation}
The numerator is $O(z^2)$ at zero and decays at infinity.
Choose the additive constant so that $\widetilde w_p^i(0,\infty)=0$.
After $k$ axial differentiations, the wall equations determine
the entrance derivatives of orders $3k+2$ and $3k+3$, with nonzero
coefficients since $\bar u_y(0,0)>0$. The remaining terms contain
only lower entrance jets and known data. Choose the free jets and
realize the prescribed jets by a fixed compactly supported extension.

Subtract $\int_0^y h_p^i(x,z)\,dz$ from
$\widetilde w_p^i-\widetilde w_p^i(x,0)$ and apply a compactly
supported divergence-free wall lifting. The resulting system is the
linearized Prandtl equation with homogeneous wall data. The primitive
of $h_p^i$ is bounded and enters the forcing multiplied by $\bar u_y$.
The estimates in \cite[Theorem~1.3 and Proposition~2.1]{guo2021regularity}
apply to the forced system, with the normal entrance datum given by
\eqref{eq:inner-normal-inlet} and compatible with the tangential trace
as in \cite[Section~2.2]{guo2021regularity}.
A slightly smaller positive exponential weight absorbs the polynomial
factors in the forcing.
This gives the bound for $\widetilde u_p^i$ in
\eqref{eq:inner-uniform-estimate}. The normal profile is recovered by
\[
 \widetilde w_p^i(x,y)=\int_y^\infty
                (\widetilde u_{p,x}^i-h_p^i)(x,z)\,dz .
\]
The tail integral requires one additional axial derivative and a larger
weight. For total output order $J$, apply the cited theorem with
$M_0=J+1$. It requires normal entrance derivatives through
$2J+5$, exponentially weighted source derivatives
$\partial_x^{k+1}\partial_y f$ for $0\le k\le J+1$, and
$\langle y\rangle\partial_x^{J+3}\partial_y f\in L^2$,
where $f$ is the forcing after the liftings.
The equation gives the remaining normal derivatives.
The entrance jets through normal order $3J+9$ ensure compatibility
through axial order $J+2$.

For the scalar entrance trace, impose
\[
\begin{aligned}
 &\widetilde S_p^i(0,0)=-S_e^i(0,0),\qquad 
 \widetilde S_{p,yy}^i(0,0)=-g_p^i(0,0),\\
 &\widetilde S_{p,yyy}^i(0,0)
 =-\bar u_y(0,0)S_{e,x}^i(0,0)-\partial_yg_p^i(0,0).
\end{aligned}
\]
Further axial differentiation determines the higher normal jets.
Realize these finitely many jets by a compactly supported entrance
trace and subtract $-S_e^i(x,0)\varphi(y)$.
Estimate~\eqref{eq:sw-Prandtl-scalar-energy} and its differentiated
forms, with the smaller exponential weight, give the scalar solution
and bounds. The constants and length conditions depend only on the
leading coefficients and finite orders. The source norms occur on
the right-hand sides.
\end{proof}

\subsection{Uniform estimates for the finite hierarchy}

\begin{proposition}\label{prop:finite-profile-hierarchy}
For fixed $N\ge8$ and $Q\ge0$, take $M=M(N,Q)$ sufficiently large
in \eqref{eq:sw-profile-input-data} and in its finite compatibility
conditions. There is $L_0>0$ such that, for every $0<L\le L_0$,
there are raw profiles and compatible higher entrance and Euler
end data, all independent of $\e$, satisfying
\begin{equation}\label{eq:all-profile-estimates}
\begin{aligned}
 &\sum_{i=0}^{N}\sum_{a+b\le Q}
  \left\|\partial_x^a\partial_Y^b
       \left(u_e^i,\frac{w_e^i}{r},S_e^i,p_e^i\right)
                      \right\|_{L^\infty([0,L]\times[0,1])}\\
 &\quad+
 \sum_{i=0}^{N}\sum_{a+b\le Q}
  \left\|\langle y\rangle^Q\partial_x^a\partial_y^b
       (\widetilde u_p^i,\widetilde w_p^i,\widetilde S_p^i)
                      \right\|_{L^\infty([0,L]\times\mathbb R_+)}\\
 &\quad+
 \sum_{i=1}^{N+1}\sum_{a+b\le Q}
  \left\|\langle y\rangle^Q\partial_x^a\partial_y^b
                         \widetilde p_p^i
                      \right\|_{L^\infty([0,L]\times\mathbb R_+)}
 \le C_{N,Q}.
\end{aligned}
\end{equation}
The pressures are normalized as above. The constant is independent
of $\e$ and $L$.
\end{proposition}

\begin{proof}
Starting from Proposition~\ref{thm:leading-prandtl-one-wall}, construct
the successive profiles by Proposition~\ref{prop:sw-Euler-finite-cylinder}
and Lemma~\ref{lem:inner-profiles}. Integrating the radial equation
from infinity gives each inner pressure, including
$\widetilde p_p^{N+1}$ after the last scalar profile.
For Euler output order $m$, Proposition~\ref{prop:sw-Euler-finite-cylinder}
uses data through $3m+6$. The source formula uses at most four further
derivatives of the preceding profiles. The inner equations require the source and entrance orders in
Lemma~\ref{lem:inner-profiles}, one additional axial derivative for
the normal tail, and finitely many corner jets. The scalar estimate
requires one further axial derivative before recovering normal
derivatives. Choose the weights so that the tail integrals for the pressure and
normal velocity converge. Working backwards from $(N,Q)$ then determines
how large $M$ must be to supply the required derivatives and corner
compatibility conditions.
Choose finitely many strictly decreasing positive exponential rates
below $\sigma_0$.
The length bound depends only on the leading coefficients, the
Euler potential and the finite orders. The estimates are therefore
uniform for $0<L\le L_0$.
\end{proof}

\subsection{Truncation and closure}
\label{subsec:profile-truncation}

The continuity equations imply
\begin{equation}
 \partial_x\bigl(\widetilde u_p^i-y\widetilde u_p^{i-1}\bigr)
 +\partial_y\bigl(\widetilde w_p^i-y\widetilde w_p^{i-1}\bigr)=0.
 \label{eq:profile-stream-identity}
\end{equation}
For $i<N$, the wall value $w_e^{i+1}(x,0)=-\widetilde w_p^i(x,0)$
contributes a term independent of $y$ to the normal velocity in the
Prandtl equations.
At order $N$, the final stream term below produces
$\widetilde w_p^N-\widetilde w_p^N(x,0)$ in both the meridional
and scalar equations.

Set
\begin{equation}
\begin{aligned}
 \Psi_p^\e(x,Y)
 =&-\chi(Y)\sum_{i=0}^{N-1}\e^{\frac{i+1}{2}}
       \int_{\frac{Y}{\sqrt\e}}^{\infty}
       \bigl(\widetilde u_p^i(x,z)-z \widetilde u_p^{i-1}(x,z)\bigr)\,dz\\
 &+\chi(Y)\e^{\frac{N+1}{2}}
       \int_0^{\frac{Y}{\sqrt\e}}
       \bigl(\widetilde u_p^N(x,z)-z \widetilde u_p^{N-1}(x,z)\bigr)\,dz .
\end{aligned}
\label{eq:profile-stream-cutoff}
\end{equation}
For each $0\le i\le N$, let
$\e^{\frac{i}{2}}(u_p^i,\sqrt\e\,w_p^i)$ be the velocity field obtained by applying
$r^{-1}(\partial_Y,-\partial_x)$ to the $i$-th summand in
\eqref{eq:profile-stream-cutoff}. Set
$S_p^i(x,y)=\chi(\sqrt\e\,y)\widetilde S_p^i(x,y)$ for $0\le i\le N$,
and $p_p^i(x,y)=\chi(\sqrt\e\,y)\widetilde p_p^i(x,y)$ for $1\le i\le N+1$.
The fields in \eqref{eq:sw-remainder-ansatz} are
\begin{equation}
\begin{aligned}
 &u_s(x,Y)=\sum_{i=0}^N\e^{\frac{i}{2}}u_e^i(x,Y)
                    +\frac{\partial_Y\Psi_p^\e(x,Y)}r,\\
 &w_s(x,Y)=\sum_{i=0}^N\e^{\frac{i}{2}}w_e^i(x,Y)
                    -\frac{\partial_x\Psi_p^\e(x,Y)}r,\\
 &s_s(x,Y)=\sum_{i=0}^N\e^{\frac{i}{2}}S_e^i(x,Y)
                    +\chi(Y)\sum_{i=0}^N\e^{\frac{i}{2}}\widetilde S_p^i(x,y),\\
 &p_s(x,Y)=\sum_{i=0}^N\e^{\frac{i}{2}}p_e^i(x,Y)
                    +\chi(Y)\sum_{i=1}^{N+1}\e^{\frac{i}{2}}\widetilde p_p^i(x,y).
\end{aligned}
\label{eq:profile-composite}
\end{equation}
By \eqref{eq:profile-stream-identity},
\[
\left.\frac{\partial_Y\Psi_p^\e}{r}\right|_{Y=0}
   =\sum_{i=0}^N\e^{\frac{i}{2}}\widetilde u_p^i(x,0),\qquad\left.-\frac{\partial_x\Psi_p^\e}{r}\right|_{Y=0}
   =\sum_{i=0}^{N-1}\e^{\frac{i+1}{2}}\widetilde w_p^i(x,0).
\]
The stream representation and \eqref{eq:profile-matching} give
\begin{equation}\label{eq:original-method-exact-continuity}
 (u_s)_x+(w_s)_Y-\frac{w_s}{r}=0,
 \qquad (u_s,w_s,s_s)|_{Y=0}=(0,0,0).
\end{equation}
The final stream term has zero normal wall trace, and
$\widetilde w_p^N(x,0)$ is determined by its tail integral.
Near the axis the inner terms vanish, leaving the regular Euler profiles.

The resulting end data are those in \eqref{eq:sw-physical-end-data}.

\subsection{Estimates for the approximate solution}
\label{subsec:sw-approximate-estimates}

\begin{proposition}[Approximate solution and residual estimates]
\label{prop:sw-uniform-profiles}\label{prop:profile-residual}
For $N\ge8$, the constructed fields satisfy no slip, exact continuity,
and
\begin{equation}\label{eq:sw-profile-final-residual}
 \|\mathfrak R_\beta w\|\lesssim\e^{\frac{N}{2}+\frac{1}{4}},
 \qquad
 \|\mathfrak R_Sw\|\lesssim\e^{\frac{N+1}{2}}.
\end{equation}
Thus,
\begin{equation}\label{eq:sw-residual-target}
 \sqrt\e\,\|\mathfrak R_\beta w\|+\|\mathfrak R_Sw\|
 \lesssim\e^{\frac{9}{2}}.
\end{equation}
The constants are independent of $\e$ and $0<L\le L_0$.
\end{proposition}

\begin{proof}
Continuity and no slip follow from
\eqref{eq:original-method-exact-continuity}.
The Euler equations cancel all coefficients through $\e^{\frac{N}{2}}$.
For the outer sums, the axial residual has a regular radial lift
and the radial residual is $r$ times a regular radial function.
Their curl divided by $r$ is therefore regular at the axis.
The outer curl and scalar residuals are $O(\e^{\frac{N+1}{2}})$
in the required finite radial norms.

Where $\chi(Y)=1$, the stream contributions telescope to
\[
\frac{\partial_Y\Psi_p^\e}{r}
 =\sum_{i=0}^N\e^{\frac{i}{2}}\widetilde u_p^i
       +\e^{\frac{N+1}{2}}\frac{y\widetilde u_p^N}{r},\qquad-\frac{\partial_x\Psi_p^\e}{r}
 =\sum_{i=0}^N\e^{\frac{i+1}{2}}\widetilde w_p^i
       -\e^{\frac{N+1}{2}}\frac{\widetilde w_p^N(x,0)}r
       +\e^{\frac{N+2}{2}}\frac{y\widetilde w_p^N}{r}.
\]
The terminal constant acting on a leading inner $Y$-derivative has
order $\e^{\frac{N}{2}}$ and is included in the last profile equations through
$\widetilde w_p^N-\widetilde w_p^N(x,0)$.
Its remaining inner contribution gains $\sqrt\e$ from
$r^{-1}-1=\frac{\sqrt\e\,y}{r}$. Products with outer fields have order
$\e^{\frac{N+1}{2}}$.

The pressure coefficients through $\widetilde p_p^{N+1}$ cancel the
radial coefficients through $\e^{\frac{N}{2}}$. After subtracting the outer
equations, the weighted profile bounds give decaying inner momentum
and scalar residuals of order $\e^{\frac{N+1}{2}}$.
The meridional curl takes one physical $Y$-derivative.
On the cutoff transition, the earlier stream tails are bounded by
any required finite power of $\e$ by choosing a sufficiently large
polynomial weight. As $\e\to0$, the last stream integral tends to
\[
 \int_0^\infty
     \bigl(\widetilde u_p^N(x,z)
                   -z\widetilde u_p^{N-1}(x,z)\bigr)\,dz .
\]
This integral is bounded and depends only on $x$. Its contribution
after truncation, the required derivatives of that contribution in
$(x,Y)$, and the resulting term in the transport of $K$ are all
$O(\e^{\frac{N+1}{2}})$.
Thus, for any required fixed $m$,
\[
|\mathfrak R_\beta|
   \lesssim \e^{\frac{N}{2}}\langle \frac{Y}{\sqrt\e}\rangle^{-m}
                         +\e^{\frac{N+1}{2}},\qquad|\mathfrak R_S|
   \lesssim \e^{\frac{N+1}{2}}
                  \bigl(\langle \frac{Y}{\sqrt\e}\rangle^{-m}+1\bigr)
\]
on the inner support, where $r\ge\frac{1}{2}$.
Near the axis, use the outer estimates above.
For $m>\frac{1}{2}$, the change of variables $Y=\sqrt\e\,y$ gives
a factor $\sqrt L\,\e^{\frac{1}{4}}$, and hence
\[
 \|\langle \frac{Y}{\sqrt\e}\rangle^{-m}\|_{L^2((0,L)\times(0,\frac{1}{2}))}
 \lesssim\e^{\frac{1}{4}}.
\]
Since $w=r^{\frac{3}{2}}$,
\[
\begin{aligned}
 &\|\mathfrak R_\beta w\|
 \lesssim\e^{\frac{N}{2}+\frac{1}{4}}+\e^{\frac{N+1}{2}}
 \lesssim\e^{\frac{N}{2}+\frac{1}{4}},\\
 &\|\mathfrak R_Sw\|
 \lesssim\e^{\frac{N+1}{2}+\frac{1}{4}}+\e^{\frac{N+1}{2}}
 \lesssim\e^{\frac{N+1}{2}}.
\end{aligned}
\]
For $N\ge8$, these estimates give \eqref{eq:sw-residual-target}.
\end{proof}

\begin{theorem}[Approximate profiles]\label{thm:approximate-profiles}
Fix a finite set of indices covered by the construction, with $i=0,1$,
$j,k,l\ge0$, and either $m=0$ or $0\le m\le l-\frac12$.
For all sufficiently small $\e$, with the threshold independent of $L$,
the approximate fields satisfy
\begin{equation}\label{eq:usbound}
\begin{aligned}
 &\sup_{i,j,k,l,m}\Bigg[
 \e^{\max\{\frac{i+2j-l}{2},0\}}
 \left\|\frac{Y^l\partial_Y^i\partial_x^k
       \{\cL^j u_s,\bar\cL^j u_s,\cL^j s_s,\bar\cL^j s_s\}}
                         {u_s^m}\right\|_\infty\\
 &+\frac1{\sqrt\e}\left(
       \|Y\partial_x^{k+1}u_s\|_\infty
       +\left\|\partial_x^k\left(\frac{w_s}{r}\right)\right\|_\infty
                         \right)\Bigg]\lesssim1.
\end{aligned}
\end{equation}
Also $u_s\sim \frac{Y}{\sqrt\e}$ for $Y\le\sqrt\e$ and $u_s\sim1$
otherwise. Near the axis, $u_s,\frac{w_s}{r},s_s$ have regular radial lifts.
For the retained axial orders,
\begin{equation}\label{eq:profile-quotient-bounds}
\begin{aligned}
 &|\partial_x^ku_s|\lesssim u_s,\qquad
   \left|\frac{w_s}{r}\right|\lesssim\sqrt\e\,u_s,\\
 &|s_s|+|(s_s)_x|+|Y(r(s_s)_Y-2s_s)|\lesssim u_s,\\
 &\left\|\frac{u_{s,x}}{u_s}\right\|_\infty
 +\sqrt\e\left\|\partial_Y\frac{u_{s,x}}{u_s}\right\|_\infty
 +\sqrt\e\left\|Y\cL\frac{u_{s,x}}{u_s}\right\|_\infty
 \lesssim1.
\end{aligned}
\end{equation}
\end{theorem}

\begin{proof}
Matching and the finite profile bounds give
\[
 |u_s(x,Y)-\bar u(x,\frac{Y}{\sqrt\e})|
       \lesssim\sqrt\e\,\min\{\frac{Y}{\sqrt\e},1\},
\]
and
\[
 |s_s(x,\sqrt\e y)-\kappa-\widetilde S_p^0(x,y)|
       \lesssim\sqrt\e\,y \qquad (y\le1).
\]
The positivity of the leading profile gives the bounds for $u_s$,
with an $\e$ threshold independent of $L$.
For $Y\le\sqrt\e$, each term in the first line of
\eqref{eq:usbound} is bounded by a constant times
\[
 \e^{\max\{\frac{l-i-2j}{2},0\}}y^{\,l-m}\le1.
\]
On $Y\ge\sqrt\e$, use the weighted decay of the inner profiles
and the smooth outer bounds. Both $\cL$ and $\bar\cL$ preserve
smooth functions of $r^2$ near the axis. The cutoff terms lie in
$r\ge\frac{1}{2}$. Their decaying tails and bounded terminal stream constant
satisfy the same estimates.

The only order-one axial term is
$1+\widetilde u_p^0(x,\frac{Y}{\sqrt\e})$, whose $x$-derivatives gain
$\sqrt\e$ after multiplication by $Y$.
All higher corrections and every term in $\frac{w_s}{r}$ contain a factor
$\sqrt\e$, so the weighted profile bounds give the second line of
\eqref{eq:usbound}.
Wall vanishing gives the first and third bounds in
\eqref{eq:profile-quotient-bounds}. On $y\le1$, division at the simple
wall zero gives
\[
 \frac{u_{s,x}(x,Y)}{u_s(x,Y)}
 =\frac{\int_0^1u_{s,xY}(x,tY)\,dt}
        {\int_0^1u_{s,Y}(x,tY)\,dt}.
\]
Multiplying numerator and denominator by $\sqrt\e$ gives a
denominator bounded below and bounded inner derivatives.
Each $Y$-derivative contributes $\e^{-\frac{1}{2}}$. The factor $Y$ in
the last term of \eqref{eq:profile-quotient-bounds} supplies
$\sqrt\e$.
On $Y\ge\sqrt\e$, use $u_s\sim1$, the inner decay and the regular
outer bounds. Finally, continuity and $w_s(x,0)=0$ give
\[
 \frac{w_s(x,Y)}r
 =-\frac1{r^2}\int_0^Y(1-t)u_{s,x}(x,t)\,dt.
\]
For $Y\le\sqrt\e$, its absolute value is bounded by
$C\frac{Y^2}{\sqrt\e}\lesssim\sqrt\e\,u_s$.
For $Y\ge\sqrt\e$, the second line of \eqref{eq:usbound} suffices.
\end{proof}

\section{Preliminary estimates}\label{sec:preliminary}
Recall from Theorem~\ref{thm:approximate-profiles} that $u_s\sim1$ for
$\sqrt\e\le Y\le1$ and $\sqrt\e\,u_s\sim Y$ for
$0\le Y\le\sqrt\e$.
\begin{lemma}\label{lem:weighted-poincare}
Let $f\in\mathcal H^1(w)$ vanish at $x=0$ or $x=L$. Then
$\|fw\|\le L\|f_xw\|$. For fixed $k\ge0$,
$\|u_s^kfw\|\lesssim_k L\|u_s^kf_xw\|$.
\end{lemma}

\begin{lemma}[Weighted Hardy inequality]
Let $f\in\mathcal H^1(w)$ have zero trace at $Y=0$. Then
\begin{equation}\label{eq:hardy-wall}
        \left\|\frac{f}{Y}w\right\|
        \le C\|f_Yw\|.
\end{equation}
\end{lemma}

\begin{proof}
For almost every $x$, the classical Hardy inequality gives
\[
 \int_0^{\frac{1}{2}}\frac{|f|^2}{Y^2}r^3\,dY
 \le4\int_0^{\frac{1}{2}}|f_Y|^2\,dY
 \le32\int_0^{\frac{1}{2}}|f_Y|^2r^3\,dY.
\]
For $\frac{1}{2}<Y<1$, $Y^{-2}\le4$, and
\begin{align*}
|f(Y)|^2
=\left|\int_0^Yf_Y(s)\,ds\right|^2&\le\left(\int_0^Y|f_Y(s)|^2(1-s)^3\,ds\right)
\left(\int_0^Y(1-s)^{-3}\,ds\right)\\
&\lesssim r(Y)^{-2}\int_0^Y|f_Y(s)|^2r(s)^3\,ds.
\end{align*}
Multiplying by $r(Y)^3$ and integrating gives
\begin{align*}
 \int_{\frac{1}{2}}^1r(Y)\int_0^Y|f_Y(s)|^2r(s)^3\,ds\,dY
 &=\int_0^1|f_Y(s)|^2r(s)^3
 \int_{\max\{\frac{1}{2},s\}}^1r(Y)\,dY\,ds\\
 &\lesssim\int_0^1|f_Y(s)|^2r(s)^3\,ds.
\end{align*}
Combining the two regions and integrating in $x$ proves the estimate.
\end{proof}

\begin{lemma}[One-dimensional finite-interval trace]
\label{lem:one-dimensional-trace}
Let $0<a<1$ and $f\in H^1(0,a)$. Then $|f(0)|^2\le a^{-1}\|f\|_{L^2(0,a)}^2+2\|f\|_{L^2(0,a)}\|f_Y\|_{L^2(0,a)}$.
\end{lemma}
\begin{proof}
Integration by parts and Cauchy--Schwarz give
\begin{align*}
 |f(0)|^2
 =-\int_0^a\partial_Y\left[\left(1-\frac Ya\right)f(Y)^2\right]\,dY
 &=\frac1a\int_0^af(Y)^2\,dY
 -2\int_0^a\left(1-\frac Ya\right)f(Y)f_Y(Y)\,dY\\
 &\le\frac1a\|f\|_{L^2(0,a)}^2
 +2\|f\|_{L^2(0,a)}\|f_Y\|_{L^2(0,a)}.
\end{align*}
\end{proof}

\begin{corollary}
Let $0<a<1$ and $H\in\mathcal H^2(w)$. Then
\[
 \|H_Y(\cdot,0)\|_{L_x^2}
 \le C_a\left(\|H_Yw\|
 +\|H_Yw\|^{\frac{1}{2}}
 \bigl(\|\cL H\,w\|+\|H_Yw\|\bigr)^{\frac{1}{2}}\right).
\]
\end{corollary}

\begin{proof}
Since $H_{YY}=\cL H+3r^{-1}H_Y$ and $r\ge1-a$ on $(0,a)$,
\[
 \|H_{YY}\|_{L^2((0,L)\times(0,a))}
 \le C_a\bigl(\|\cL H\,w\|+\|H_Yw\|\bigr).
\]
Apply Lemma~\ref{lem:one-dimensional-trace} at each fixed $x$, then integrate in $x$ to obtain
\[
 \|H_Y(\cdot,0)\|_{L_x^2}^2
 \le C_a\|H_Yw\|^2
 +C_a\|H_Yw\|\bigl(\|\cL H\,w\|+\|H_Yw\|\bigr).
\]
\end{proof}

\section{Radial regularity and preliminary recovery}\label{sec:axis}

\begin{lemma}[Near-axis slice estimates]
\label{lem:near-axis-slice-estimates}

For $\frac{1}{2}<Y<1$, write
$\|F(Y)\|_2:=\|F(\cdot,Y)\|_{L^2(0,L)}$.

\begin{enumerate}
\item[1.] If \(
        Fw,F_Yw\in L^2(\Omega)
\), then
\begin{align}
        &\|F(Y)\|_2
        \lesssim
        \|F w\|+r^{-1}\|F_Y w\|.
        \label{eq:rough-FTC}
\end{align}
\item[2.] If \( Fw,\cL F\,w\in L^2(\Omega)
\), then
\begin{align}
        &\|F_Y(Y)\|_2
        \lesssim
        r^{-1}\|\cL F w\|,
        \label{eq:radial-FY-upper}
\end{align}
and
\begin{align}
        &\|F(Y)\|_2
        \lesssim
        \|F w\|+  \Lambda(r) \|\cL F w\|.
        \label{eq:log-FTC-from-LF}
\end{align}
The radial flux satisfies
\begin{equation}\label{eq:radial-flux-upper}
 r(Y)^3F_Y(\cdot,Y)
 =-\int_Y^1r(s)^3\cL F(\cdot,s)\,ds
 \quad\text{in }L_x^2,\qquad 0<Y<1.
\end{equation}
Hence, for $0<\eta<\frac12$,
\begin{equation}\label{eq:radial-flux-tail-sharp}
 \eta\|F_Y(\cdot,1-\eta)\|_{L_x^2}
 \le\frac12\|\cL F\,w\|_{L^2((0,L)\times(1-\eta,1))}
 \to0\qquad(\eta\to0).
\end{equation}
\end{enumerate}
\end{lemma}

\begin{proof}
To prove \eqref{eq:rough-FTC}, let $Y>\frac12$. Since $\chi(Y)=0$
and $\chi(0)=1$, the fundamental theorem of calculus gives
\[
        F(Y)
        =
        \int_0^Y\partial_s\bigl[(1-\chi(s))F(s)\bigr]\,ds
        =
        -\int_0^Y\chi'(s)F(s)\,ds
        +
        \int_0^Y(1-\chi(s))F_Y(s)\,ds .
\]
Taking \(L_x^2\)-norms,
\begin{align}
        &\|F(Y)\|_2
        \leq
        \int_0^Y|\chi'(s)|\|F(s)\|_2\,ds
        +
        \int_0^Y(1-\chi(s))\|F_Y(s)\|_2\,ds .
        \label{eq:cutoff-FTC-start}
\end{align}
The first term is supported where $\frac14\leq s\leq\frac12$ and $\frac12\leq r(s)\leq\frac34$. Thus
\begin{align}
 \int_0^Y|\chi'(s)|\|F(s)\|_2\,ds
\lesssim
        \left(
        \int_{\frac{1}{4}}^{\frac{1}{2}}\|F(s)\|_2^2\,ds
        \right)^{\frac{1}{2}}
        \lesssim
        \left(
        \int |F|^2r^3\,dY\,dx
        \right)^{\frac{1}{2}}
        =
        \|F w\| .
        \label{eq:cutoff-F-term}
\end{align}
For the second term, Cauchy's inequality gives
\begin{align*}
        \int_0^Y(1-\chi(s))\|F_Y(s)\|_2\,ds
        &\leq
        \int_{\frac{1}{4}}^{Y}\|F_Y(s)\|_2\,ds = 
        \int_{\frac{1}{4}}^{Y}
        r(s)^{-\frac{3}{2}}
        r(s)^{\frac{3}{2}}\|F_Y(s)\|_2\,ds  \\
        &\leq
        \left(
        \int_{\frac{1}{4}}^{Y}r(s)^{-3}\,ds
        \right)^{\frac{1}{2}}
        \left(
        \int_{\frac{1}{4}}^{Y}r(s)^3\|F_Y(s)\|_2^2\,ds
        \right)^{\frac{1}{2}}.
\end{align*}
Since
\[
        \int_{\frac{1}{4}}^{Y}r(s)^{-3}\,ds
        =
        \int_{\frac{1}{4}}^{Y}(1-s)^{-3}\,ds
        \lesssim r^{-2},
        \qquad Y=1-r,
\]
we obtain $\int_0^Y(1-\chi(s))\|F_Y(s)\|_2\,ds \lesssim r^{-1}\|F_Y w\|$. Together with \eqref{eq:cutoff-FTC-start} and \eqref{eq:cutoff-F-term},
this proves \eqref{eq:rough-FTC}.

For \eqref{eq:radial-FY-upper}, use $\partial_Y(r^3F_Y)=r^3\cL F$. Local one-dimensional elliptic regularity gives
$F\in H^2_{\mathrm{loc}}((0,1);L_x^2)$. Since
$r^3\cL F\in L_Y^1L_x^2$, there is a limiting flux
$C_F\in L_x^2$ such that $r^3F_Y(Y)=C_F-\int_Y^1r(s)^3\cL F(s)\,ds$. We claim that $C_F=0$. If $C_F\ne0$, Cauchy's inequality gives
\[
 \left\|\int_Y^1r(s)^3\cL F(s)\,ds\right\|_{L_x^2}
 \le\frac{r^2}{2}\|\cL Fw\|_{L^2((0,L)\times(Y,1))}
 \to0
 \qquad(Y\to1^-).
\]
Pairing the flux identity with $C_F$, we therefore obtain, for
$Y$ sufficiently close to $1$,
\begin{align*}
 \frac{d}{dY}\langle F(\cdot,Y),C_F\rangle_{L_x^2}
 =r^{-3}\left(\|C_F\|_{L_x^2}^2
   -\left\langle\int_Y^1r(s)^3\cL F(s)\,ds,C_F\right\rangle_{L_x^2}\right)\ge\frac12\|C_F\|_{L_x^2}^2r^{-3}.
\end{align*}
Choose $Y_0<1$ in this interval with $F(\cdot,Y_0)\in L_x^2$.
Such a point exists because $Fw\in L^2(\Omega)$ and the weight is bounded
below on every compact subinterval of $(0,1)$. Integration gives
\begin{align*}
 &\langle F(\cdot,Y),C_F\rangle_{L_x^2}
 \ge\langle F(\cdot,Y_0),C_F\rangle_{L_x^2}
       +\frac12\|C_F\|_{L_x^2}^2
                       \int_{Y_0}^Y(1-s)^{-3}\,ds\\
 &=\langle F(\cdot,Y_0),C_F\rangle_{L_x^2}
       +\frac14\|C_F\|_{L_x^2}^2
                       \bigl(r^{-2}-(1-Y_0)^{-2}\bigr).
\end{align*}
For $Y$ sufficiently close to $1$,
\[
 \|F(\cdot,Y)\|_{L_x^2}
 \ge\frac{\langle F(\cdot,Y),C_F\rangle_{L_x^2}}{\|C_F\|_{L_x^2}}
 \ge\frac18\|C_F\|_{L_x^2}r^{-2}.
\]
It follows that
\[
 \int^1\|F(\cdot,Y)\|_{L_x^2}^2r^3\,dY
 \ge\frac{\|C_F\|_{L_x^2}^2}{64}\int^1\frac{dY}{1-Y}
 =\infty,
\]
contrary to $Fw\in L^2(\Omega)$. Thus $C_F=0$, proving
\eqref{eq:radial-flux-upper}.
Another application of Cauchy's inequality yields
\begin{align*}
 \|F_Y(Y)\|_2
 \le r^{-3}\int_Y^1r(s)^3\|\cL F(s)\|_2\,ds
 \le r^{-3}\left(\int_Y^1r(s)^3\,ds\right)^{\frac12}
                         \|\cL Fw\|
 \le\frac1{2r}\|\cL Fw\|.
\end{align*}
This proves \eqref{eq:radial-FY-upper}.

The sharper tail estimate \eqref{eq:radial-flux-tail-sharp} follows by
applying Cauchy's inequality to
\[
 \eta^3F_Y(\cdot,1-\eta)
 =-\int_{1-\eta}^1r(s)^3\cL F(\cdot,s)\,ds.
\]
On every truncated interval, integration by parts yields
\begin{align*}
 \int_0^{1-\eta}\!\|\cL F(\cdot,Y)\|_{L_x^2}^2r^3\,dY
 +3\eta^2\|F_Y(\cdot,1-\eta)\|_{L_x^2}^2
 =&\int_0^{1-\eta}\!\|F_{YY}(\cdot,Y)\|_{L_x^2}^2r^3\,dY\\*
 &+3\int_0^{1-\eta}\!\|F_Y(\cdot,Y)\|_{L_x^2}^2r\,dY
 +3\|F_Y(\cdot,0)\|_{L_x^2}^2.
\end{align*}
The wall trace exists by local one-dimensional elliptic regularity,
since $Fw,\cL F\,w\in L^2(\Omega)$ and the coefficients are regular near $Y=0$.
The tail estimate \eqref{eq:radial-flux-tail-sharp} bounds the
right-hand integrals uniformly. Letting $\eta\to0$ gives
\begin{equation}\label{eq:radial-second-derivative-identity}
 \|\cL F\,w\|^2
 =\|F_{YY}w\|^2
  +3\left\|\frac{F_Y}{r}w\right\|^2
  +3\|F_Y(\cdot,0)\|_{L_x^2}^2.
\end{equation}
Thus,
\begin{equation}\label{equiv-norm}
 \left\|\frac{F_Y}{r}w\right\|+\|F_{YY}w\|
 \lesssim\|\cL Fw\|.
\end{equation}
For \eqref{eq:log-FTC-from-LF}, use the same cutoff.
By \eqref{eq:cutoff-FTC-start} and
\eqref{eq:cutoff-F-term}, $\|F(Y)\|_2 \lesssim \|F w\|+ \int_0^Y(1-\chi(s))\|F_Y(s)\|_2\,ds$. For the last term,
\begin{align*}
        \int_0^Y(1-\chi(s))\|F_Y(s)\|_2\,ds
        &\leq
        \int_{\frac{1}{4}}^{Y}\|F_Y(s)\|_2\,ds  \\
        &\leq
        \left(
        \int_{\frac{1}{4}}^{Y}r(s)^{-1}\,ds
        \right)^{\frac{1}{2}}
        \left(
        \int_{\frac{1}{4}}^{Y}r^3(s)\|\frac{F_Y(s)}{r}\|_2^2\,ds
        \right)^{\frac{1}{2}}\leq
        \Lambda(r)\left\|\frac{F_Y}{r}w \right\|.
\end{align*}
Together with \eqref{equiv-norm}, this proves
\eqref{eq:log-FTC-from-LF}.
\end{proof}
For $F\in\mathcal H^k(w)$ and $k=0,1,2,3,4$,
\begin{equation}\label{eq:radial-sobolev-norms}
 \|F\|_{\dot H^k(w)}^2
 \sim\sum_{\substack{a+b+2j=k\\a,j\ge0,\ b=0,1}}
       \|\partial_x^a\partial_Y^b\cL^jF\,w\|^2.
\end{equation}
For smooth radial lifts, \eqref{eq:radial-flux-upper} gives
\[
 \partial_r^j\left(\frac{\partial_r\widetilde F}{r}\right)(x,r)
 =\int_0^1t^{3+j}
       (\partial_r^j\widetilde{\cL F})(x,tr)\,dt,
 \qquad j=0,1,2.
\]
Radial differentiation, Minkowski's inequality and \eqref{equiv-norm}
give \eqref{eq:radial-sobolev-norms}. Radial smooth approximation
extends it to $\mathcal H^k(w)$.

For $G\in\mathcal H_0^4(w)$ and $0\le a\le j\le3$,
the trace of $\partial_x^a\nabla_\xi^{j-a}G^\sharp$ vanishes at
$x=0$ for even $a$ and at $x=L$ for odd $a$.
Poincar\'e's inequality in $x$ gives
\begin{equation}\label{eq:mixed-Sobolev-Poincare}
 \|G\|_{\dot H^j(w)}\lesssim L\|G\|_{\dot H^{j+1}(w)},
 \qquad j=0,1,2,3.
\end{equation}
Together with \eqref{eq:radial-sobolev-norms}, this yields
\begin{equation}\label{eq:main-H-recovery}
 \mathcal H_4(G)\sim\e\|G\|_{H^2(w)}+\e^2\|G\|_{H^4(w)}.
\end{equation}

\subsection{The stream function transformation at the axis}\label{sec:mixed-quotient-derivation}
We justify the change of variables $G=-\frac{\Psi}{r^2}$ at the axis
for the principal variational solution.
Let
$\Delta_*:=\partial_x^2+\partial_Y^2+r^{-1}\partial_Y=:\partial_x^2+D$.
For $Y<1$, the identity $\Delta_*(r^2G)=r^2\cA G$ puts
\eqref{eq:main-equation} in the form
\[
\begin{aligned}
 u_s(\Delta_*\Psi)_x+w_s(\Delta_*\Psi)_Y
       +\frac{2w_s}{r}\Delta_*\Psi-\e\Delta_*^2\Psi-r\Psi_Y(\beta_s)_x+r\Psi_x(\beta_s)_Y=-r^2F_G.
\end{aligned}
\]
For $\frac{f_\Psi}{\sqrt r}\in L^2(\Omega)$, consider the principal problem
\begin{equation}\label{eq:appA-short-Psi-principal}
\left\{
\begin{aligned}
 &\e\Delta_*^2\Psi=f_\Psi,\\
 &\Psi|_{x=0}=\Psi_{xx}|_{x=0}=0,\\
 &\Psi_x|_{x=L}=\Psi_{xxx}|_{x=L}=0,\\
 &\Psi|_{Y=0}=\Psi_Y|_{Y=0}=0.
\end{aligned}
\right.
\end{equation}
Its variational form is
\begin{equation}\label{eq:appA-short-Psi-weak}
\begin{aligned}
 &\e\mathfrak a_\Psi(\Psi,\Phi)
 =\int_\Omega\frac{f_\Psi\Phi}{r}\,dx\,dY,\\
 &\mathfrak a_\Psi(\Psi,\Phi)
 :=\int_\Omega\frac{D\Psi D\Phi+2\Psi_{xY}\Phi_{xY}
                              +\Psi_{xx}\Phi_{xx}}r\,dx\,dY.
\end{aligned}
\end{equation}
Here $\Phi$ has finite $\mathfrak a_\Psi(\Phi,\Phi)$ and satisfies
$\Phi=\Phi_Y=0$ at $Y=0$, $\Phi|_{x=0}=0$, and $\Phi_x|_{x=L}=0$.
By Poincar\'e's inequality in $x$ and the same modal estimates as in
Lemma~\ref{lem:main-principal-solvability}, we obtain
\begin{equation}\label{eq:appA-short-Psi-graph}
 \left\|\frac{\Psi}{\sqrt r}\right\|
 +\left\|\frac{\Psi_Y}{\sqrt r}\right\|
 +\left\|\frac{D\Psi}{\sqrt r}\right\|
 +\left\|\frac{D^2\Psi}{\sqrt r}\right\|
 \le\frac{C_L}{\e}\left\|\frac{f_\Psi}{\sqrt r}\right\|.
\end{equation}

\begin{lemma}[Axis behavior of the stream function]
\label{lem:appA-short-Psi-axis}
Let $\Psi$ be the variational solution of
\eqref{eq:appA-short-Psi-principal}, defined by
\eqref{eq:appA-short-Psi-weak}, with $\frac{f_{\Psi}}{\sqrt r}\in L^2(\Omega)$.
There exists $\gamma_\Psi\in L^2(0,L)$ such that, for $\frac{1}{2}<Y<1$,
\begin{align}
 &\|\Psi(\cdot,Y)\|_{L_x^2}
 +r\|\Psi_Y(\cdot,Y)\|_{L_x^2}
 \le \frac{C_L}{\e}r^2\Lambda(r)
       \left\|\frac{f_{\Psi}}{\sqrt r}\right\|,
 \label{eq:appA-short-Psi-decay}\\
 &\|\Psi(\cdot,Y)-r^2\gamma_\Psi\|_{L_x^2}
 +r\|\Psi_Y(\cdot,Y)+2r\gamma_\Psi\|_{L_x^2}
 \le \frac{C_L}{\e}r^4\Lambda(r)
       \left\|\frac{f_{\Psi}}{\sqrt r}\right\|.
 \label{eq:appA-short-Psi-expansion}
\end{align}
Thus $\Psi|_{Y=1}=\Psi_Y|_{Y=1}=0$ and
$G(\cdot,Y)\to-\gamma_\Psi$ strongly in $L_x^2$ as $Y\to1$.
\end{lemma}

\begin{proof}

For \eqref{eq:appA-short-Psi-decay}, choose
$Y_0\in(\frac{1}{4},\frac{1}{2})$ such that
$\|\Psi_Y(\cdot,Y_0)\|_{L_x^2}
 \lesssim\|\frac{\Psi_Y}{\sqrt r}\|$.
Since $\partial_Y\left(\frac{\Psi_Y}{r}\right)=\frac{D\Psi}{r}$, integration from $Y_0$ and Cauchy's inequality give

\begin{align*}
 \frac{\|\Psi_Y(\cdot,Y)\|_{L_x^2}}r
 &\le
 \frac{\|\Psi_Y(\cdot,Y_0)\|_{L_x^2}}{1-Y_0}
 +\int_{Y_0}^Y\frac{\|D\Psi(\cdot,s)\|_{L_x^2}}{1-s}\,ds\\*
 &\lesssim
 \left\|\frac{\Psi_Y}{\sqrt r}\right\|
 +\left(\log\frac{1-Y_0}{1-Y}\right)^{\frac{1}{2}}
  \left(
    \int_{Y_0}^Y
    \frac{\|D\Psi(\cdot,s)\|_{L_x^2}^2}{1-s}\,ds
  \right)^{\frac{1}{2}}\le
 \frac{C_L}{\e}\Lambda(r)\left\|\frac{f_{\Psi}}{\sqrt r}\right\|,
\end{align*}

where we used \eqref{eq:appA-short-Psi-graph}. Hence

\begin{equation}\label{eq:appA-short-radial-derivative}
 \|\Psi_Y(\cdot,Y)\|_{L_x^2}
 \le\frac{C_L}{\e}r\Lambda(r)\left\|\frac{f_{\Psi}}{\sqrt r}\right\|.
\end{equation}

Since $(1-Y)\Lambda(1-Y)$ is integrable near $Y=1$,
$\Psi(\cdot,Y)$ is Cauchy in $L_x^2$ as $Y\to1$ and has
a strong trace $\Psi(\cdot,1)$. This trace must vanish.
Otherwise, for some $\eta>0$, $\|\Psi(\cdot,Y)\|_{L_x^2} \ge\frac12\|\Psi(\cdot,1)\|_{L_x^2}$, $1-\eta<Y<1$, and therefore

\[
 \int_{1-\eta}^1
 \frac{\|\Psi(\cdot,Y)\|_{L_x^2}^2}{1-Y}\,dY
 \ge\frac14\|\Psi(\cdot,1)\|_{L_x^2}^2
       \int_{1-\eta}^1\frac{dY}{1-Y}
 =\infty.
\]

This contradicts \eqref{eq:appA-short-Psi-graph}.
Thus $\Psi(\cdot,1)=0$, and
\eqref{eq:appA-short-radial-derivative} gives

\[
\begin{aligned}
\|\Psi(\cdot,Y)\|_{L_x^2}
 \le\int_Y^1\|\Psi_Y(\cdot,s)\|_{L_x^2}\,ds\le\frac{C_L}{\e}r^2\Lambda(r)
          \left\|\frac{f_{\Psi}}{\sqrt r}\right\|.
\end{aligned}
\]

This proves \eqref{eq:appA-short-Psi-decay} and the vanishing axis traces.

To prove \eqref{eq:appA-short-Psi-expansion}, local one-dimensional regularity gives

\[
 \int_{\frac{1}{4}}^{\frac{1}{2}}\|(D\Psi)_Y(\cdot,s)\|_{L_x^2}^2\,ds
 \lesssim
 \left\|\frac{D\Psi}{\sqrt r}\right\|^2
 +\left\|\frac{D^2\Psi}{\sqrt r}\right\|^2.
\]

Choose an interior slice using this estimate and integrate

\[
 \partial_Y\left(\frac{(D\Psi)_Y}{r}\right)=\frac{D^2\Psi}{r}
\]

as above to obtain

\[
 \|(D\Psi)_Y(\cdot,Y)\|_{L_x^2}
 \le\frac{C_L}{\e}r\Lambda(r)\left\|\frac{f_{\Psi}}{\sqrt r}\right\|.
\]

Thus $D\Psi$ also has a strong $L_x^2$ trace at $Y=1$.
Since $\frac{D\Psi}{\sqrt r}\in L^2$, the same argument
shows that this trace is zero. Hence

\[
 \|D\Psi(\cdot,Y)\|_{L_x^2}
 \le\int_Y^1\|(D\Psi)_Y(\cdot,s)\|_{L_x^2}\,ds
 \le\frac{C_L}{\e}r^2\Lambda(r)\left\|\frac{f_{\Psi}}{\sqrt r}\right\|.
\]

Hence $\frac{D\Psi}{r}\in L_Y^1L_x^2$ near $Y=1$.
The identity $\partial_Y(\frac{\Psi_Y}{r})=\frac{D\Psi}{r}$ therefore gives a
finite strong $L_x^2$ limit $-\frac{\Psi_Y(\cdot,Y)}{2r}\to\gamma_\Psi$, with

\begin{align*}
\left\|\frac{\Psi_Y(\cdot,Y)}r+2\gamma_\Psi\right\|_{L_x^2}
 =\left\|
     \int_Y^1\frac{D\Psi(\cdot,s)}{1-s}\,ds
   \right\|_{L_x^2}\le\frac{C_L}{\e}r^2\Lambda(r)
          \left\|\frac{f_{\Psi}}{\sqrt r}\right\|.
\end{align*}

Finally, since $\Psi(\cdot,1)=0$,

\begin{align*}
 \|\Psi(\cdot,Y)-r^2\gamma_\Psi\|_{L_x^2}
 &\le
 \int_Y^1(1-s)
 \left\|\frac{\Psi_Y(\cdot,s)}{1-s}
              +2\gamma_\Psi\right\|_{L_x^2}\,ds\\
 &\le
 \frac{C_L}{\e}\left\|\frac{f_{\Psi}}{\sqrt r}\right\|
 \int_Y^1(1-s)^3\Lambda(1-s)\,ds\le
 \frac{C_L}{\e}r^4\Lambda(r)\left\|\frac{f_{\Psi}}{\sqrt r}\right\|.
\end{align*}

This proves \eqref{eq:appA-short-Psi-expansion}.
Dividing by $r^2$ gives $G(\cdot,Y)\to-\gamma_\Psi$ in $L_x^2$ as $Y\to1$.

\end{proof}

\subsection{Near-axis blow-up estimates}

\begin{lemma}[Five-dimensional Agmon embedding]
Let \(0<\e\leq1\) and $G\in\mathcal H_0^4(w)$. Then, with constants
allowed to depend on $L$,
\begin{equation}
\label{eq:G-Linfty-H4eps}
        \|G\|_{L^\infty(\Omega)}
        \lesssim
        \e^{-\frac54}\mathcal H_4(G),
\end{equation}
and
\begin{equation}
\label{eq:grad-G-Linfty-H4eps}
        \|G_x\|_{L^\infty(\Omega)}
        +
        \|G_Y\|_{L^\infty(\Omega)}
        \lesssim
        \e^{-\frac74}\mathcal H_4(G).
\end{equation}
\end{lemma}

\begin{proof}

Since \(\Omega^\sharp\subset\mathbb R^5\), the five-dimensional Agmon
inequalities give
\begin{align*}
&\|G^\sharp\|_{L^\infty(\Omega^\sharp)}
        \lesssim
        \|G^\sharp\|_{H^2(\Omega^\sharp)}^{\frac34}
        \|G^\sharp\|_{H^4(\Omega^\sharp)}^{\frac14},\\
&\|\nabla_{x,\xi}G^\sharp\|_{L^\infty(\Omega^\sharp)}
        \lesssim
        \|G^\sharp\|_{H^2(\Omega^\sharp)}^{\frac14}
        \|G^\sharp\|_{H^4(\Omega^\sharp)}^{\frac34}.
\end{align*}
By \eqref{eq:main-H-recovery},
\begin{align*}
&\|G^\sharp\|_{L^\infty}
        \lesssim
        \bigl(\e^{-1}\mathcal H_4(G)\bigr)^{\frac34}
        \bigl(\e^{-2}\mathcal H_4(G)\bigr)^{\frac14}
        =
        \e^{-\frac54}\mathcal H_4(G),\\
&\|\nabla_{x,\xi}G^\sharp\|_{L^\infty}
        \lesssim
        \bigl(\e^{-1}\mathcal H_4(G)\bigr)^{\frac14}
        \bigl(\e^{-2}\mathcal H_4(G)\bigr)^{\frac34}
        =
        \e^{-\frac74}\mathcal H_4(G).
\end{align*}
Finally, $\|G^\sharp\|_{L^\infty(\Omega^\sharp)} = \|G\|_{L^\infty(\Omega)}$, and, for \(r>0\), $\partial_xG^\sharp(x,\xi)=G_x(x,1-|\xi|)$, $|\nabla_\xi G^\sharp(x,\xi)|=|G_Y(x,1-|\xi|)|$. Hence
\[
        \|\nabla G\|_{L^\infty(\Omega)}
        \lesssim
        \|\nabla_{x,\xi}G^\sharp\|_{L^\infty(\Omega^\sharp)}
        \lesssim
        \e^{-\frac74}\mathcal H_4(G).
\]
This proves \eqref{eq:G-Linfty-H4eps} and
\eqref{eq:grad-G-Linfty-H4eps}.
\end{proof}

The following wall estimates use the argument
of~\cite[Lemma~3.8]{guo2023validity}.

\begin{lemma}\label{lem:wall-traces}
Let $G\in\mathcal H_0^4(w)$. For
$j=0,1,2$ and $k=0,1,2,3$, we have
\begin{align*}
 &\sup_{0<Y_0\le\frac{1}{2}}
 \Big[\|\nabla^kG\|_{Y=Y_0}+\|\nabla^jq\|_{Y=Y_0}
       +\sqrt{Y_0}\|\nabla^3q\|_{Y=Y_0}\Big]
 \le C_\e,
\end{align*}
where $C_\e\le C_{\e,L}\|G^\sharp\|_{H^4(\Omega^\sharp)}$.
The local trace and Hardy estimates apply since $\frac12\le r\le1$
near $Y=0$ and $\frac{u_s}{Y}$ is smooth and positive there.
\end{lemma}

\begin{proof}[Proof of Proposition~\ref{prop:near-axis-H4-blowup}]

\medskip
\noindent
\textbf{Step 1: estimates for $\cL G$.}
By \eqref{eq:radial-sobolev-norms}, $\cL G\,w\in L^2(\Omega)$, $\cL^2G\,w\in L^2(\Omega)$. Applying Lemma~\ref{lem:near-axis-slice-estimates} with \(F=\cL G\) gives

\begin{align}
        &\|(\cL G)_Y(Y)\|_2
        \lesssim
        Mr^{-1},
        \label{eq:H_Y-rminus1}
        \\
        &\|\cL G(Y)\|_2
        \lesssim
        M\Lambda(r)
        \label{eq:H-log-bound}
\end{align}

Since $(\cL G)_xw\in L^2(\Omega)$, $(\cL G)_{xY}w\in L^2(\Omega)$, the rough estimate \eqref{eq:rough-FTC} gives

\begin{align}
        &\|(\cL G)_x(Y)\|_2
        \lesssim
        Mr^{-1}.
        \label{eq:H_x-rminus1}
\end{align}

It follows that

\begin{equation}
        \|\nabla_{x,Y}\cL G(Y)\|_2
        \leq
        \|(\cL G)_x(Y)\|_2+\|(\cL G)_Y(Y)\|_2
        \lesssim
        Mr^{-1}.
        \label{eq:grad-LG-rminus1}
\end{equation}

\medskip
\noindent
\textbf{Step 2: second derivatives of $G$.}
Since $G_{xx}w,\cL G_{xx}w\in L^2(\Omega)$,
Lemma~\ref{lem:near-axis-slice-estimates} gives
\begin{align}
        &\|G_{xx}(Y)\|_2
        \lesssim
        M\Lambda(r),
        \label{eq:Gxx-log}
        \\
        &\|G_{xxY}(Y)\|_2
        \lesssim
        Mr^{-1}.
        \label{eq:GxxY-rminus1}
\end{align}

Since $\cL G_x=(\cL G)_x$, the flux formula \eqref{eq:radial-flux-upper} gives $r^3G_{xY}(Y) = -\int_Y^1r(s)^3(\cL G)_x(s)\,ds$. Using \eqref{eq:H_x-rminus1},

\begin{align}\label{eq:GxY-bdd}
        \|G_{xY}(Y)\|_2
        \leq r^{-3}\int_Y^1r(s)^3\|(\cL G)_x(s)\|_2\,ds 
        \leq
        Mr^{-3}\int_Y^1r(s)^2\,ds
        \lesssim M.
\end{align}

Finally, $G_{YY}=\cL G+\frac3rG_Y$. The flux formula for \(F=G\) gives $r^3G_Y(Y) = -\int_Y^1r(s)^3\cL G(s)\,ds$. Therefore $\frac{G_Y(Y)}{r} = -r^{-4}\int_Y^1r(s)^3\cL G(s)\,ds$. Using \eqref{eq:H-log-bound},

\begin{align*}
        \left\|\frac{G_Y(Y)}{r}\right\|_2
        \leq
        r^{-4}\int_Y^1r(s)^3\|\cL G(s)\|_2\,ds \leq
        Mr^{-4}\int_Y^1r(s)^3\Lambda(r(s))\,ds
        \lesssim
        M\Lambda(r).
\end{align*}

Hence $\|G_{YY}(Y)\|_2 \lesssim M\Lambda(r)$.

\medskip
\noindent
\textbf{Step 3: third derivatives of $G$.}
Since $G_{xxx}w,G_{xxxY}w\in L^2(\Omega)$,
\eqref{eq:rough-FTC} gives $\|G_{xxx}(Y)\|_2\lesssim Mr^{-1}$.
Equation~\eqref{eq:GxxY-rminus1} gives the same bound for $G_{xxY}$.

Next, $(\cL G)_x=\cL G_x = G_{xYY}-\frac3rG_{xY}$. Therefore $G_{xYY}=(\cL G)_x+\frac3rG_{xY}$. Using \eqref{eq:H_x-rminus1} and \eqref{eq:GxY-bdd}, we get $\|G_{xYY}(Y)\|_2 \lesssim Mr^{-1}$. It remains to estimate \(G_{YYY}\). Differentiating the flux formula $G_Y(Y)=-r(Y)^{-3}\int_Y^1r(s)^3\cL G(s)\,ds$ with respect to $Y$ gives
\[
 G_{YYY}(Y)=(\cL G)_Y(Y)+\frac{12}{r(Y)^5}
       \int_Y^1r(s)^3[\cL G(Y)-\cL G(s)]\,ds.
\]
For $Y\leq s<1$, \eqref{eq:H_Y-rminus1} implies
\[
 \|\cL G(Y)-\cL G(s)\|_2
 \leq\int_Y^s\|(\cL G)_Y(t)\|_2\,dt
 \lesssim M\log\frac{r(Y)}{r(s)}.
\]
Since $\int_Y^1r(s)^3\log\frac{r(Y)}{r(s)}\,ds=\frac{r(Y)^4}{16}$, we obtain $\|G_{YYY}(Y)\|_2\lesssim Mr^{-1}$. Combining the estimates for \(G_{xxx}\), \(G_{xxY}\), \(G_{xYY}\), and
\(G_{YYY}\), we get $\|\nabla_{x,Y}^3G(Y)\|_2 \lesssim Mr^{-1}$. Together with \eqref{eq:grad-LG-rminus1}, this proves
\eqref{eq:near-axis}.
\end{proof}

\begin{corollary}
Assume the hypotheses of Proposition~\ref{prop:near-axis-H4-blowup}.
If $i_n+j_n+2k_n\leq3$ and $j_n\in\{0,1\}$ for $n=1,2$, then
\begin{align}&\label{axis-vanish}
 (\partial_x^{i_1}\partial_Y^{j_1}\cL^{k_1}G,
  \partial_x^{i_2}\partial_Y^{j_2}\cL^{k_2}G\,w^2)_{Y=1}=0.
\end{align}
The pairing denotes the limit at $Y=1-\eta$ as $\eta\to0$.
\end{corollary}

\begin{proof}

By Proposition~\ref{prop:near-axis-H4-blowup}, together with
\eqref{eq:rough-FTC} for the derivatives of order at most one,

\[
 \eta^3
 \|\partial_x^{i_1}\partial_Y^{j_1}\cL^{k_1}G(\cdot,1-\eta)\|_{L_x^2}
 \|\partial_x^{i_2}\partial_Y^{j_2}\cL^{k_2}G(\cdot,1-\eta)\|_{L_x^2}
 \lesssim\eta M^2\to0.
\]

\end{proof}
\begin{lemma}\label{lem:sw-near-axis}
Let $s^\e\in\mathcal H^2(w)$. For $0<\eta<\frac{1}{2}$,
\[
 \Lambda(\eta)^{-1}\|s^\e(\cdot,1-\eta)\|_{L_x^2}
 +\eta\|\{s_Y^\e,s_x^\e\}(\cdot,1-\eta)\|_{L_x^2}
 \lesssim\|s^\e\|_{H^2(w)}.
\]
\end{lemma}
\begin{proof}
Apply part 2 of Lemma~\ref{lem:near-axis-slice-estimates} to $s^\e$,
and \eqref{eq:rough-FTC} to $s_x^\e$, since
$s_x^\e w,s_{xY}^\e w\in L^2(\Omega)$.
\end{proof}

\begin{corollary}\label{cor:sw-axis-flux}
Under the assumptions of Lemma~\ref{lem:sw-near-axis}, at $Y=1-\eta$,
\begin{align*}
 &\eta^3\left|\int_0^L w_s(s^\e)^2\,dx\right|
 \lesssim\eta^4\Lambda(\eta)^2\|s^\e\|_{H^2(w)}^2,\\
 &\eta^3\left|\int_0^L s_Y^\e s^\e\,dx\right|
 \lesssim\eta^2\Lambda(\eta)\|s^\e\|_{H^2(w)}^2,\\
 &\eta^3\left|\int_0^L s_Y^\e s_x^\e\,dx\right|
 \lesssim\eta\|s^\e\|_{H^2(w)}^2.
\end{align*}
All three expressions tend to zero as $\eta\to0$.
\end{corollary}
\begin{proof}
Lemma~\ref{lem:sw-near-axis}, Cauchy's inequality and
$|w_s|\lesssim r$ from \eqref{eq:usbound} give the estimates.
\end{proof}

\begin{lemma}[Green formula]
Let $f,g\in\mathcal H^3(w)$. Then $(\cL f,gw^2) =-(f_Y,g_Yw^2)-(f_Y,gw^2)_{Y=0}$.
The boundary term at $Y=1-\eta$ tends to zero as $\eta\to0$.
\end{lemma}

\begin{proof}
\begin{align*}
&(\cL f,gw^2)
=(\partial_Y(f_Yr^3),g)
=-(f_Y,g_Yr^3)+(f_Yr^3,g)_{Y=1}-(f_Yr^3,g)_{Y=0}.
\end{align*}
For the axis term, the radial slice estimates give
\[
 \eta^3\|f_Y(\cdot,1-\eta)\|_{L_x^2}
       \|g(\cdot,1-\eta)\|_{L_x^2}
 \le C\eta^2\Lambda(\eta)
       \|f\|_{H^2(w)}\|g\|_{H^2(w)}\to0.
\]
Letting $\eta\to0$ proves the formula.
\end{proof}

\subsection{Basic recovery estimates}
\begin{lemma}\label{lem:basic-recovery}
Let $G\in\mathcal H_0^4(w)$ and $w=r^{\frac{3}{2}}$.
We use the profile bounds \eqref{eq:profile-quotient-bounds}.
For $L$ sufficiently small,
\begin{itemize}

\item[1.] The following Poincare type inequalities hold: 
\begin{align}
&\| \p_Y^{j} \p_x^{j_2 }G  w\| \lesssim L \| \p_Y^j \p_x^{j_2 + 1} G  w\| \text{ for } j_2=0,1,2,3,\quad j+j_2\leq3, \notag\\
&\| u_s^k \p_Y^j q  w\| \lesssim_k L \| u_s^k \p_Y^j q_x  w\| \text{ for fixed } k\ge0,\quad 0\leq j\leq2, \label{Poincare.2}\\
&\|u_s^k\p_x^{j_2}q\,w\|\lesssim_k L\|u_s^k\p_x^{j_2+1}q\,w\|
 \text{ for fixed }k\ge0,\quad j_2=1,2.\notag
\end{align}

\item[2.] The sum $\cD_2[q]+\cD_3[q]$ controls the following quantities. 
\begin{equation}    
\begin{aligned} &\label{systemat1}
\| \{\eps\cL  q, \sqrt{\eps}q_{xY}, \sqrt{\eps} q_{xx} \}  w\| + \| \sqrt{\e} \nabla q  w\|  + \| \sqrt{\e}\frac{q_x}{Y}  w\|  \\
&\hspace{5 mm}  + \| \{\eps^{\frac{3}{2}}(\cL G)_Y, \eps \cL G_{x}, \eps G_{xxY}, \eps G_{xxx} \}  w\|  \\
 &\hspace{5 mm} + \| \{\eps \cL G,\sqrt{\eps} G_{xY}, \sqrt{\eps} G_{xx} \}  w\|  + \| \sqrt\eps \nabla G  w\|  \\
&\lesssim \cD_2[q]+\cD_3[q] 
\end{aligned}
\end{equation}

\item[3.] The following quantities can be controlled with a pre-factor of $o_L(1)$:
\begin{align} 
\begin{aligned} &\label{small.guys}
\| \eps \cL G   w\| + \| \sqrt{\eps} \{ \nabla q, \nabla G\}  w\| \lesssim L\bigl(\cD_2[q]+\cD_{3,x}[q]\bigr), \\
&\| \eps \cL q  w\| \lesssim \sqrt{L} \cD_3.
\end{aligned}
\end{align}

\item[4.] Fix $0<\eta\leq1$. The following interpolation estimate holds: 
\begin{align} &\label{interp.1}
\| \sqrt{\eps}\nabla q_x  w\|  \le \eta \cD_{3,x} + C\eta^{-2} \| \sqrt{\eps} u_s \nabla q_x    w\|.
\end{align}

For fixed $\alpha>0$, take $\eta=L^{\frac{\alpha}{2}}$:
\begin{align} &\label{interp.L}
\| \sqrt{\eps}\nabla q_x  w\| \lesssim L^{\frac{\alpha}{2}}\cD_{3,x} + L^{-\alpha}\|\sqrt{\eps} u_s \nabla q_x    w\|.
\end{align}

\item[5.] The following boundary estimates are valid: 
\begin{align}
&\|\sqrt{\eps}q_x w\|_{x=L} + \|\sqrt{\eps} q_x w\|_{x=0} + \|\eps \cL G w\|_{x=L}  \notag\\ 
 &\hspace{2 mm} + \|\sqrt{\eps} q_{xx} w\|_{x=0} + \|\eps \sqrt{u_s} q_{xx} w\|_{x=L} + \|\sqrt{\eps}q_{xY} w\|_{x=L} \lesssim \sqrt L\bigl(\cD_2[q]+\cD_{3,x}[q]\bigr)  \label{b1}\\
&\|\sqrt{\eps}q_Y w\|_{x=L} +\|\sqrt{\eps}G_Y w\|_{x=L}
 +\|\sqrt{\eps}\tfrac{G}{Y} w\|_{x=L}\notag\\
&\qquad+\|\sqrt{\eps}\tfrac{q_x}{Y}w\|_{x=L}
 +\|\eps\sqrt{u_s}\cL q_xw\|_{x=L}
 \lesssim\sqrt L\bigl(\cD_2[q]+\cD_{3,x}[q]\bigr).\label{b3}
\end{align}

\end{itemize}

\end{lemma}

\begin{proof}

The $G$ estimates in \eqref{Poincare.2} follow from the mixed boundary
conditions. Since $|u_{sx}|\lesssim u_s$, the weights $u_s(x,Y)^k$ and
$u_s(0,Y)^k$ are comparable for each fixed $k\ge0$.
The same argument applies to $u_s^k\partial_Y^jq$.
For almost every $Y$, the boundary conditions in \eqref{eq:main-equation}
and $G=u_sq$ give, successively for
$j_2=1,2$,
\[
 \|\partial_x^{j_2}q\|_{L_x^2}
 \le L\|\partial_x^{j_2+1}q\|_{L_x^2}
       +CL\|\partial_x^{j_2}q\|_{L_x^2}.
\]
Absorb the last term, multiply by $u_s(0,Y)^kw(Y)$, and integrate in $Y$.

\textit{ Proof of (\ref{interp.1}),(\ref{interp.L}):}
Write $y=Y/\sqrt\eps$.
For $0<\xi\ll1$, integration by parts gives

\begin{align*}
&\eps^2\| \cL q  \chi(y/\xi)  w\|^2 =  \eps^2( \p_Y \{ Y \}, (\cL q)^2 w^2 \chi(y/\xi)^2) \\
&=  - \eps^\frac52( \frac{2Y}{\sqrt\eps} \cL q, (\cL q)_Y w^2 \chi(y/\xi)^2) + \eps^\frac52( \frac{3Y}{\sqrt\eps} \cL q, \frac{\cL q}{r} w^2 \chi(y/\xi)^2) \\
 &- \eps^\frac52( \frac{2Y}{\sqrt\eps} (\cL q)^2, w^2 \frac{1}{\xi\sqrt{\eps}} \chi(y/\xi)\chi'(y/\xi)) \\
&\lesssim  L \|\eps \sqrt{u_s} \cL q_{x}  \chi(y/\xi)  w\| \| \eps^\frac32 \sqrt{u_s} (\cL q)_{Y}   w\| \\
 &+ L^2\sqrt{\eps}\| \eps\sqrt{u_s} \cL q_{x}  w\|^2 \sup_{Y\le\xi\sqrt\eps/2} |\frac1r| + \frac{L^2}{\xi} \|\eps \sqrt{u_s} \cL q_{x}  w\|^2 \\
&\lesssim  (L + \frac{L^2}{\xi}) \cD_3^2.
\end{align*}

\noindent  Here we used (\ref{Poincare.2}). Hence

\begin{align*} &\n
\eps\| \cL q   w\|  \le \eps\| \cL q  w [1- \chi(y/\xi)] \| + \eps\|  \cL q  w \chi(y/\xi) \| \\ &\n
 \lesssim \frac{L}{\sqrt{\xi}}  \eps\| \sqrt{u_s}  \cL q_x  w\| + (\frac{L}{\sqrt{\xi}} + \sqrt{L}) \cD_3 \\
 &\lesssim \sqrt{L} \cD_3 \text{ for } \xi = L. 
\end{align*}

For $\nabla q_x$, the same argument gives a bound without
the pre-factor of $o_L(1)$. For the far-field component,
$\| \nabla q_x [1-\chi(y/\xi)]  w\| \lesssim \frac{1}{\xi} \| u_s \nabla q_x   w\|$. For the localized component, integration by parts gives 

\begin{align*}
&\eps\| q_{xx}  \chi(y/\xi)  w\|^2 =  - \eps( Y, \p_Y \{ q_{xx}^2 w^2 \chi(y/\xi)^2 \}) \\
&=  - \eps^\frac32( \frac{2Y}{\sqrt{\eps}} q_{xx}, q_{xxY} w^2 \chi(y/\xi)^2 ) + \eps^\frac{3}{2}( 3 \frac{Y}{\sqrt\eps} q_{xx}, \frac{q_{xx}}{r} w^2\chi(y/\xi)^2) \\
 &- 2\eps( \frac{Y}{\xi\sqrt\eps} q_{xx}^2, w^2 \chi(y/\xi)\chi'(y/\xi) ) \\
&\lesssim  \sqrt{\xi}\bigO(\sqrt{\text{LHS}})\|\e \sqrt{u_s} q_{xxY}  w\| +\sqrt\e \xi \sup_{Y\le\xi\sqrt\eps/2} |\frac1r|  \|\sqrt\e q_{xx} w\|^2 \\
 &+ \xi^{-2}\|\sqrt\e u_s q_{xx}w\|^2.  
\end{align*}

Combining these estimates gives 

\begin{equation}\label{eq:wall-local-interpolation}
\|\sqrt\e q_{xx}w\|^2
 \lesssim \xi\|\e\sqrt{u_s}q_{xxY}w\|^2
       +\xi^{-2}\|\sqrt\e u_sq_{xx}w\|^2.
\end{equation}

For \(q_{xY}\),

\begin{align*}
&\eps\| q_{xY}  \chi(y/\xi)  w\|^2 =  - \eps( Y, \p_Y \{ q_{xY}^2 w^2 \chi(y/\xi)^2 \}) \\
&=  - \eps^\frac32( \frac{2Y}{\sqrt{\eps}} q_{xY}, \cL q_x w^2 \chi(y/\xi)^2 ) - \eps^\frac{3}{2}( 3 \frac{Y}{\sqrt\eps} q_{xY}, \frac{q_{xY}}{r} w^2\chi(y/\xi)^2) \\
 &- 2\eps( \frac{Y}{\xi\sqrt\eps} q_{xY}^2, w^2 \chi(y/\xi)\chi'(y/\xi) ) \\
&\lesssim  \sqrt{\xi} \bigO(\sqrt{\text{LHS}})\|\e \sqrt{u_s} \cL q_{x}  w\| + \xi^{-2}\|\sqrt\e u_s q_{xY}w\|^2.  
\end{align*}

Combining this with the far-field estimate gives

\begin{align*} 
&\|\sqrt\e q_{xY}w\|^2
 \lesssim \xi\|\e\sqrt{u_s}\cL q_xw\|^2
       +\xi^{-2}\|\sqrt\e u_sq_{xY}w\|^2.
\end{align*}

Taking square roots and setting $\xi=c\eta^2$ proves
\eqref{interp.1}. Take $\eta=L^{\frac{\alpha}{2}}$ for \eqref{interp.L}.
For every fixed $\beta>0$,

\[
 \|\sqrt\e\nabla q_xw\|^2
 +\|\sqrt\e\nabla q_xw\|\bigl(\cD_2[q]+\cD_3[q]\bigr)
 \leq C_\beta L^{-\beta}\cD_2^2[q]+o_L(1)\cD_3^2[q],
\]

by taking $\eta=L^{\frac{\beta}{5}}$ and applying Young's inequality.
Expanding $G=u_sq$, using the profile bounds, Hardy's inequality and
\eqref{interp.1} with fixed $\eta$, gives

\[
 \|\{\e\cL G_x,\e G_{xxY},\e G_{xxx},
             \sqrt\e G_{xY},\sqrt\e G_{xx}\}w\|
 \lesssim\cD_2[q]+\cD_{3,x}[q].
\]

For the remaining third derivative, use

\[
 (\cL G)_Y=u_s(\cL q)_Y+3u_{sY}\cL q
 +(2u_{sYY}+\bar\cL u_s)q_Y+q(\cL u_s)_Y.
\]

This proves \eqref{systemat1}. Poincare in $x$ gives the first estimate
in \eqref{small.guys}.

\textit{Proof of \eqref{b1} and \eqref{b3}.}
By the fundamental theorem of calculus, $G=u_sq$, and the boundary
conditions in \eqref{eq:main-equation},

\[
 \|\sqrt\e q_xw\|_{x=0,L}
 \lesssim\sqrt L\|\sqrt\e q_{xx}w\|,\qquad
 q_{xx}|_{x=0}=-2\frac{u_{sx}}{u_s}q_x|_{x=0}.
\]

Apply the same argument to $\e\sqrt{u_s}q_{xx}$ for its outlet trace.
Also $\|\e\cL G\,w\|_{x=L}\leq\sqrt L\|\e\cL G_xw\|$.

Since \(q|_{x=0}=G|_{x=0}=0\), the fundamental theorem of calculus gives

\begin{align*}
&\left|\sqrt\e\{q_Y,G_Y,\frac{G}{Y}\}(L,Y)\right|
 \leq\int_0^L\left|\sqrt\e\{q_{xY},G_{xY},\frac{G_x}{Y}\}(x,Y)\right|\,dx\\
 &\leq\sqrt L\|\sqrt\e\{q_{xY},G_{xY},\frac{G_x}{Y}\}(\cdot,Y)\|_{L_x^2}.
\end{align*}

Taking the weighted $L_Y^2$ norm gives these first-order traces.
Since $q_x=-(\frac{u_{sx}}{u_s})q$ at $x=L$, \eqref{eq:profile-quotient-bounds}
and Hardy's inequality on this slice give
\[
 \|\sqrt\e\,(\frac{q_x}{Y})w\|_{x=L}
 \lesssim\|\sqrt\e\,q_Yw\|_{x=L}
 \lesssim\sqrt L\bigl(\cD_2[q]+\cD_{3,x}[q]\bigr).
\]
At $x=L$, the boundary condition in \eqref{eq:main-equation} gives
$q_x=-(u_{sx}/u_s)q$. Differentiate this identity in $Y$ to express
$q_{xY}$ in terms of $q$ and $q_Y$. Applying $\cL$ to the identity gives

\[
 \cL q_x|_{x=L}=-\left[
  \frac{u_{sx}}{u_s}\cL q
  +2\partial_Y\left(\frac{u_{sx}}{u_s}\right)q_Y
  +\cL\left(\frac{u_{sx}}{u_s}\right)q\right]_{x=L}.
\]

Write $q,q_Y,\cL q$ on the right as integrals in $x$ from zero.
The bounds \eqref{eq:profile-quotient-bounds}, together with Hardy's inequality, give
$\|\e\sqrt{u_s}\cL q_xw\|_{x=L}
\lesssim\sqrt L(\cD_2[q]+\cD_{3,x}[q])$.

\end{proof}

\subsection{The main equation}

\begin{proposition}[Solvability of the main equation]
\label{prop:main-solvability}
Assume that, for every source $F_G$ with $F_Gw\in L^2(\Omega)$, solutions
$G\in\mathcal H_0^4(w)$ of \eqref{eq:main-equation} satisfy
\[
 \cD[q]\lesssim\|F_Gw\|.
\]
Then, for each source $F_G$ with $F_Gw\in L^2(\Omega)$, there exists a unique solution
$G\in\mathcal H_0^4(w)$ to
\[
\left\{
\begin{aligned}
 &\partial_xR[q]-\e\cA^2G+\cJ[G]=F_G,\\
 &G|_{x=0}=G_{xx}|_{x=0}=0,\\
 &G_x|_{x=L}=G_{xxx}|_{x=L}=0,\\
 &G|_{Y=0}=G_Y|_{Y=0}=0.
\end{aligned}
\right.
\]
The solution satisfies $\cD[q]\lesssim\|F_Gw\|$.
\end{proposition}

\begin{lemma}\label{lem:main-principal-solvability}
Given $F_G$ with $F_Gw\in L^2(\Omega)$, there is a unique solution
$G\in\mathcal H_0^4(w)$ of $\e\cA^2G=F_G$ with the boundary
conditions in \eqref{eq:main-equation}. It satisfies
\begin{align}&\label{sum.Gxxxx.lap}
 \|\{\e^2\cL G_{xx},\e^2G_{xxxY},\e^2G_{xxxx}\}w\|^2
 \lesssim |(F_G,\e^3G_{xxxx}w^2)|,\\
 &\label{sum.Gxxyy.lap}
 \|\{\e^2(\cL G_x)_Y,\e^2\cL G_{xx},\e^2G_{xxxY}\}
       \sqrt{u_s}w\|^2
 -C\e\bigl(\cD_2^2[q]+\cD_3^2[q]\bigr)\lesssim |(F_G,\e^3u_s\cL G_{xx}w^2)|.
\end{align}
\end{lemma}

\begin{proof}Fix \(F_G^m\in C_c^\infty((0,L)\times(0,1))\) such that
$        \|(F_G-F_G^m)w\|\to0$.
Extend \(F_G^m\) evenly across \(x=L\) to \((0,2L)\), then oddly
across \(x=0\) and periodically with period \(4L\). Denote the resulting
function by \(\widetilde F_G^m\). Its sine series contains only positive odd
modes, by the symmetry about \(x=L\). Set $\lambda_n:=\frac{n\pi}{2L}$, $n=1,3,5,\ldots$. All modal sums below run over these indices. On this periodic domain, we solve $\e\cA^2\bar G^m=\widetilde F_G^m$ and write
\[
        \widetilde F_G^m
        =
        \sum_{n\ {\rm odd}}
        F_n^m(Y)\sin\left(\lambda_n x\right),\qquad
        \bar G^m
        =
        \sum_{n\ {\rm odd}}
        G_n^m(Y)\sin\left(\lambda_n x\right).
\]

Write $G^m:=\bar G^m|_\Omega$.

Since $\cA\bigl(G_n^m(Y)\sin(\lambda_n x)\bigr) = (\cL-\lambda_n^2)G_n^m(Y)\sin(\lambda_n x)$, we obtain the modal ODE
\begin{equation}\label{eq:G-modal-ODE}
        \e\cL^2G_n^m
        -2\e\lambda_n^2\cL G_n^m
        +\e\lambda_n^4G_n^m
        =
        F_n^m .        \qquad n\neq0,\quad n\ {\rm odd}.
\end{equation}

For each fixed \(n\), define
\[
        B_n[G,\phi]
        :=
        \e(\cL G,\cL\phi\,w^2)
        +2\e\lambda_n^2(G_Y,\phi_Yw^2)
        +\e\lambda_n^4(G,\phi w^2)
        :
        V_n\times V_n\to\mathbb R,
\]

where $V_n=\{G:Gw,G_Yw,\cL G\,w\in L^2(0,1),\ G(0)=G_Y(0)=0\}$, with norm $(B_n[G,G])^{\frac{1}{2}}$.

For $n\ne0$, $B_n$ is bounded and coercive on $V_n$. Lax--Milgram gives a
unique \(G_n^m\in V_n\) satisfying
\begin{equation}\label{eq:modal-variational}
        B_n[G_n^m,\phi]
        =
        (F_n^m,\phi w^2),
        \qquad
        \forall \phi\in V_n.
\end{equation}
For each fixed $n$, radial elliptic regularity for the clamped
operator $(\cL-\lambda_n^2)^2$ on $B_4(0,1)$ gives a smooth radial lift
of $G_n^m$, since $F_n^m$ has a smooth radial lift.
Taking $\phi=G_n^m$ in \eqref{eq:modal-variational} gives
\[
 \e\left(
 \|\cL G_n^m w\|
 +\lambda_n\|(G_n^m)_Yw\|
 +\lambda_n^2\|G_n^m w\|
 \right)
 \lesssim
 \lambda_n^{-2}\|F_n^m w\|.
\]
Summing over $n$, Poincar\'e's inequality and
\eqref{eq:radial-sobolev-norms} give $\e\|G^m\|_{H^2(w)}\lesssim\|F_G^m w\|$. Taking $\phi=\e^3\lambda_n^4G_n^m$,

\begin{equation}\label{eq:modal-Gxxxx-bound}
\begin{aligned}
 &\left\|\e^2\lambda_n^2\cL G_n^m w\right\|^2
 +2\left\|\e^2\lambda_n^3(G_n^m)_Yw\right\|^2
 +\left\|\e^2\lambda_n^4G_n^m w\right\|^2=\left(F_n^m,\e^3\lambda_n^4G_n^m w^2\right)
 \lesssim\e^2\|F_n^m w\|^2.
\end{aligned}
\end{equation}

Summing \eqref{eq:modal-Gxxxx-bound} over \(n\) and using Parseval's
identity gives
\begin{equation*}
 \big\|
 \{
   \e^2(\cL G^m)_{xx},
   \e^2G^m_{xxxY},
   \e^2G^m_{xxxx}
 \big\} w\|^2
 \lesssim |(F_G^m,\e^3G^m_{xxxx}w^2)|
 \lesssim\e^2\|F_G^m w\|^2.
\end{equation*}
Equation~\eqref{eq:G-modal-ODE} also gives $\|\e^2\cL^2G^m w\|\lesssim\e\|F_G^m w\|$.
For $(\cL G^m)_{xY}$, integration by parts gives

\begin{equation}
\begin{aligned}&\label{LGmnxY}
\|\e^2\lambda_n(\cL G_n^m)_Y w\|^2=-\e^4\big(\lambda_n^2\cL G_n^m,\cL^2G_n^m w^2\big)
   +\e^4\big[\lambda_n^2r^3(\cL G_n^m)_Y\cL G_n^m\big]_{Y=0}^{Y=1}
 \\
 &\leq\ \|\e^2 \lambda_n^2\cL G_n^m w\|\,\|\e^2 \cL^2G_n^m w\|
+\e^4(\lambda_n^2r^3(\cL G_n^m)_Y\cL G_n^m)_{Y=1}\\
 &+\e^4\big|\lambda_n^{\frac{1}{2}}(\cL G_n^m)_Y(0)\big|
   \big|\lambda_n^{\frac{3}{2}}\cL G_n^m(0)\big|
 \\
 &\lesssim\ \|
 \e^2 \lambda_n^2\cL G_n^m w\|\,\|\e^2\cL^2G_n^m w\|
 \\
 &+\e^4\left(
        \lambda_n^{-\frac12}\|(\lambda_n\cL G^m_n)_Y w\|
        +
        \|(\lambda_n\cL G^m_n)_Yw\|^{\frac{1}{2}}
        \left(
        \|\cL^2G^m_n w\|^{\frac{1}{2}}
        +
        \lambda_n^{-\frac12}\|(\lambda_n\cL G^m_n)_Y w\|^{\frac{1}{2}}
        \right)\right)\\
       &\quad  \times
\left(
        \lambda_n^{-\frac12}\|\lambda_n^2(\cL G^m_n) w\|
        +
        \|\lambda_n^2(\cL G^m_n)w\|^{\frac{1}{2}}
        \|\lambda_n(\cL G^m_n)_Y w\|^{\frac{1}{2}}
        \right)\\
         &\leq C\left(
 \left\|\e^2\lambda_n^2\cL G_n^m w\right\|^2
 +\|\e^2\cL^2G_n^m w\|^2\right)
 +\frac12\left\|\e^2\lambda_n(\cL G_n^m)_Yw\right\|^2
\end{aligned}
\end{equation}

where Lemma~\ref{lem:near-axis-slice-estimates} gives
\begin{align*}
 \left|\e^4\lambda_n^2r^3
          (\cL G_n^m)_Y\cL G_n^m\right|_{Y=1-\eta}\lesssim
 \e^4\lambda_n^2\eta^2\Lambda(\eta)
 \bigl(\|\cL G_n^m w\|^2+\|\cL^2G_n^m w\|^2\bigr).
\end{align*}
Letting $\eta\to0$, the axis contribution vanishes.

We estimate the wall terms using Lemma~\ref{lem:one-dimensional-trace}.

Summing \eqref{LGmnxY} over $n$ gives
\[
 \|\e^2\{G^m_{xxxx},G^m_{xxxY},\cL G^m_{xx},
              (\cL G_x^m)_Y,\cL^2G^m\}w\|
 \lesssim\e\|F_G^m w\|.
\]
Applying \eqref{eq:radial-sobolev-norms} first to finite modal sums
and passing to the full modal sum gives $\e^2\|G^m\|_{\dot H^4(w)}\lesssim\e\|F_G^m w\|$. By \eqref{eq:mixed-Sobolev-Poincare} and the second-order bound,
this proves $\e\|G^m\|_{H^4(w)}\lesssim\|F_G^m w\|$.
The same estimates for differences give convergence in
$\mathcal H_0^4(w)$. Passing to the limit $m\to\infty$ proves
\eqref{sum.Gxxxx.lap} and gives
\begin{equation}\label{eq:main-principal-H4}
 \e\|G\|_{H^4(w)}\lesssim\|F_Gw\|.
\end{equation} Uniqueness follows by taking $F_G=0$ in
 \eqref{eq:modal-variational}. Multiplying
\eqref{eq:G-modal-ODE} by $\lambda_n^j$ gives the higher axial regularity of each
$G^m$. We first integrate by parts on finite modal sums. Their radial lifts
are smooth, and the estimates above give convergence to $G^m$ in
$H^4(w)$.

Set $q^m=\frac{G^m}{u_s}$. Pairing the equation with
$\eps^3 u_s\cL G^m_{xx}$ and integrating by parts in $x$ and $Y$ gives
\begin{align*}
 &\big(\e\mathcal A^2G^m,\eps^3 u_s\cL G^m_{xx}w^2\big)
 =\
 \|\e^2\,\sqrt{u_s}\,(\cL G^m)_{xY} w\|^2
 +2\|\e^2\sqrt{u_s}\,\cL G^m_{xx} w\|^2
 \\
 &+\|\e^2\sqrt{u_s}\,G^m_{xxxY} w\|^2
 +
 \e^4\big(u_{sx}(\cL G^m)_Y,(\cL G^m)_{xY}w^2\big)
 \\
 &+\e^4\big(u_{sxY}(\cL G^m)_Y,(\cL G^m)_xw^2\big)
 +\e^4\big(u_{sY}(\cL G^m)_{xY},(\cL G^m)_xw^2\big)
 \\
 &+\e^4\big(u_{sx}G^m_{xxY},G^m_{xxxY}w^2\big)
 +\e^4\big(u_{sxY}G^m_{xxY},G^m_{xxx}w^2\big)
 \\
 &+\e^4\big(u_{sY}G^m_{xxx},G^m_{xxxY}w^2\big)\\
 &\gtrsim
 \big\|\e^2\sqrt{u_s}
 \{(\cL G^m)_{xY},\cL G^m_{xx},G^m_{xxxY}\}w\big\|^2
 -C\e\bigl(\cD_2^2[q^m]+\cD_3^2[q^m]\bigr).
\end{align*}

The side terms vanish by the mixed boundary conditions, and
the wall terms in the first integration by parts vanish since
$u_s=u_{sx}=G^m_{xxx}=0$ at $Y=0$. The axis terms vanish for the
finite modal sums by smooth radial compatibility.
For the two terms containing $u_{sY}$, integration by parts gives
\begin{align*}
 &\e^4\bigl(u_{sY}(\cL G_x^m)_Y,\cL G_x^m w^2\bigr)
  +\e^4\bigl(u_{sY}G^m_{xxx},G^m_{xxxY}w^2\bigr)\\
 &\qquad=-\frac{\e^4}{2}
       \bigl(\cL u_s,\bigl[(\cL G_x^m)^2+(G^m_{xxx})^2\bigr]w^2\bigr)
       -\frac{\e^4}{2}
       \bigl(u_{sY},(\cL G_x^m)^2w^2\bigr)_{Y=0}.
\end{align*}
Since $\|\cL u_s\|_\infty\lesssim\e^{-1}$, the interior term
is bounded by $C\e(\cD_2^2[q^m]+\cD_3^2[q^m])$ using
\eqref{systemat1}. At the wall, $\cL G_x^m=2u_{sY}q^m_{xY}+2u_{sxY}q^m_Y$, $q^m_Y|_{x=0}=0$. Thus the wall term is bounded by
\[
 C\e^{\frac{5}{2}}\bigl(\|q^m_{xY}w\|_{Y=0}^2
                  +\|q^m_Yw\|_{Y=0}^2\bigr)
 \lesssim\e\,|q^m|_{\partial,2}^2.
\]
We bound the remaining four terms using the profile estimates, \eqref{systemat1},
and Young's inequality. The upper bound in Lemma~\ref{lem:energy-norm-equivalence}
and the wall trace estimates of Lemma~\ref{lem:wall-traces} justify passage
to the limit in these terms for fixed $\varepsilon,L$.
Thus \eqref{sum.Gxxyy.lap} holds for $G$.

\end{proof}

\begin{proof}[Proof of Proposition~\ref{prop:main-solvability}]
Let $PG:=-\e\cA^2G$ with the conditions in
\eqref{eq:main-equation}. By Lemma~\ref{lem:main-principal-solvability},
for every $h\in L^2(w)$ there is a unique $G=P^{-1}h$, and
\[
 \e\|P^{-1}h\|_{H^4(w)}\lesssim\|hw\|.
\]
Write \eqref{eq:main-equation} as
\[
 PG+\mathcal KG=F_G,\qquad
 \mathcal KG:=\partial_xR[q]+\cJ[G].
\]
Setting $h:=PG$, we obtain
\[
 (I+\mathcal KP^{-1})h=F_G.
\]
By $R[q]=u_s\cA G-G\cA u_s$ and \eqref{eq:J-def},
\[
 \mathcal KG=u_s(\cA G)_x+w_s(\cA G)_Y
       +(2G-rG_Y)(\beta_s)_x+rG_x(\beta_s)_Y.
\]
For fixed $\e,L$, the radial regularity and finite derivative bounds
of $\frac{w_s}{r}$ and $\beta_s$ give a bounded map
$\mathcal K:\mathcal H^3(w)\to L^2(w)$.
The embedding $\mathcal H^4(w)\to\mathcal H^3(w)$ is compact by
Rellich's theorem on $\Omega^\sharp$. Thus
$\mathcal KP^{-1}:L^2(w)\to L^2(w)$ is compact.

If $(I+\mathcal KP^{-1})h=0$ and $G=P^{-1}h$, then $G$ solves
the homogeneous main equation. The assumed a priori estimate gives
$\cD[q]=0$, hence $q=0$ by the boundary conditions and $G=h=0$.
The Fredholm alternative therefore gives a unique $h\in L^2(w)$,
and $G=P^{-1}h$ is the unique solution. The bound follows from the assumed estimate. See also
\cite[Proposition~3.2 and Lemma~3.5]{guo2023validity} for this Fredholm argument.
\end{proof}

\subsection{Fourth-order recovery}
\begin{lemma}[Recovery of \((\cL G_x)_Y\)]
\label{lem:recover-LGxY}
Let $G\in\mathcal H^4(w)$ satisfy the $x$-boundary conditions
in \eqref{eq:main-equation}. Then
\begin{align*}
&\e^4\|(\cL G_x)_Yw\|^2
\lesssim
\bigl(\e^2\|\cL G_{xx}w\|\bigr)
\bigl(\e^2\|\cL^2G\,w\|\bigr)
+
\e^2\bigl(\e\|\cL G_xw\|\bigr)^2
+
\e^4\|\sqrt{u_s}(\cL G_x)_Yw\|^2 .
\end{align*}
\end{lemma}

\begin{proof}
Write $\Omega_0=(0,L)\times(0,\tfrac12)$ and $\Omega_1=(0,L)\times(\tfrac12,1)$.
Since $u_s\gtrsim1$ on $\Omega_1$,
\begin{equation}
\label{eq:recover-axis-side}
        \|(\cL G_x)_Yw\|_{L^2(\Omega_1)}
        \lesssim
        \|\sqrt{u_s}(\cL G_x)_Yw\|_{L^2(\Omega_1)} .
\end{equation}

On \(\Omega_0\), \(\frac12\leq r\leq1\) and \(2^{-\frac32}\leq w\leq1\).
The mixed $x$-boundary conditions give the expansion

\begin{align*}
&G(x,Y)=\sum_{n\ {\rm odd}}
        G_n(Y)\sin\left(\frac{n\pi}{2L}x\right),\\
&\cL G
        =
        \sum_{n\ {\rm odd}}
        \cL G_n(Y)\sin\left(\frac{n\pi}{2L}x\right).
\end{align*}

The one-dimensional Gagliardo--Nirenberg inequality on
\((0,\frac12)\) gives
\[
\|g_Y\|_{L^2(0,\frac12)}
        \lesssim
        \|g\|_{L^2(0,\frac12)}^{\frac{1}{2}}
        \|g_{YY}\|_{L^2(0,\frac12)}^{\frac{1}{2}}
        +
        \|g\|_{L^2(0,\frac12)} .
\]
Using $g_{YY}=\cL g+3r^{-1}g_Y$, the bound $r^{-1}\le2$ on this interval, and Young's inequality, we obtain
\[
 \|g_Y\|_{L^2(0,\frac12)}^2
 \lesssim\|g\|_{L^2(0,\frac12)}\|\cL g\|_{L^2(0,\frac12)}
                         +\|g\|_{L^2(0,\frac12)}^2.
\]
We take $g=\cL G_n$, sum over $n$, and apply Parseval's identity and the Cauchy--Schwarz inequality to obtain
\begin{align}
\|(\cL G_x)_Y\|_{L^2(\Omega_0)}^2
&\lesssim
\frac{\pi^2}{8L}\sum_{n\ {\rm odd}}n^2
\|\cL G_n\|_{L^2(0,\frac12)}
\|\cL^2G_n\|_{L^2(0,\frac12)}
\nonumber+
\frac{\pi^2}{8L}\sum_{n\ {\rm odd}}n^2
\|\cL G_n\|_{L^2(0,\frac12)}^2
\nonumber\\
&\lesssim
\|\cL G_{xx}\|_{L^2(\Omega_0)}
\|\cL^2G\|_{L^2(\Omega_0)}
+
\|\cL G_x\|_{L^2(\Omega_0)}^2 .
\label{eq:recover-wall-side}
\end{align}
Since \(2^{-\frac{3}{2}}\leq w\leq1\) on \(\Omega_0\), \eqref{eq:recover-wall-side} implies
\begin{equation}
\label{eq:recover-wall-side-weighted}
        \|(\cL G_x)_Yw\|_{L^2(\Omega_0)}^2
        \lesssim
        \|\cL G_{xx}w\|
        \|\cL^2G\,w\|
        +
        \|\cL G_xw\|^2 .
\end{equation}
Combining \eqref{eq:recover-axis-side} and
\eqref{eq:recover-wall-side-weighted}, and multiplying by \(\e^4\), yields
\[
\e^4\|(\cL G_x)_Yw\|^2
\lesssim
\bigl(\e^2\|\cL G_{xx}w\|\bigr)
\bigl(\e^2\|\cL^2G\,w\|\bigr)
+
\e^4\|\cL G_xw\|^2
+
\e^4\|\sqrt{u_s}(\cL G_x)_Yw\|^2 .
\]

\end{proof}

Lemma~\ref{lem:recover-LGxY}, \eqref{systemat1} and the definitions
of $\mathcal H_4$ and $\cD_4$ give $\mathcal H_4(G)\lesssim\cD_2[q]+\cD_3[q]+\cD_4[G]$.

\begin{lemma}[Equivalence of norms]\label{lem:energy-norm-equivalence}
For fixed $0<L\le L_0$, $0<\e\le\e_0(L)$ and
$(G,s)\in\mathcal H_0^4(w)\times\mathscr V_s^2$,
\begin{equation}\label{eq:energy-norm-equivalence}
 C_{\e,L}^{-1}\bigl(\|G\|_{H^4(w)}+\|s\|_{H^2(w)}\bigr)
 \le\cD[q]+\mathcal S_\e[s]
 \le C_{\e,L}\bigl(\|G\|_{H^4(w)}+\|s\|_{H^2(w)}\bigr).
\end{equation}
Thus $\mathscr Z_\e$ is complete.
\end{lemma}

\begin{proof}
The lower bound for $G$ follows from \eqref{eq:main-H-recovery},
\eqref{systemat1} and Lemma~\ref{lem:recover-LGxY}.
For $s$, use \eqref{equiv-norm}, Poincar\'e's inequality and
\eqref{eq:sw-scalar-norms}.

Near $Y=0$, $\frac{Y}{u_s}$ is smooth for fixed $\e$.
Since $G(x,0)=0$, the identity
\[
 \frac{G(x,Y)}Y=\int_0^1G_Y(x,tY)\,dt
\]
and Minkowski's inequality give $\|q\|_{H^3}\le C_{\e,L}\|G\|_{H^4}$
on a wall neighborhood. Away from the wall, $\frac{1}{u_s}$ has a regular
radial lift. Hence $\|q\|_{H^3(w)}\le C_{\e,L}\|G\|_{H^4(w)}$.
This controls the bulk terms in $\cD_2,\cD_3$. The trace theorem
controls their second-order traces. The same theorem controls the
third-order traces in $\cD_4$ and the first-order traces in
$\mathcal S_\e$. This proves the upper bound.
The boundary conditions define closed subspaces of
$\mathcal H^4(w)$ and $\mathcal H^2(w)$, so the equivalence gives completeness.
\end{proof}

\subsection{The swirl equation}

\begin{proposition}[Solvability of the swirl equation]
\label{prop:swirl-solvability}
Assume that, for every source $F_s$ with $F_sw\in L^2(\Omega)$, solutions
$s^\e\in\mathscr V_s^2$ of \eqref{eq:swirl-equation} satisfy
\[
 \mathcal S_\e[s^\e]\lesssim\|F_sw\|.
\]
Then, for each source $F_s$ with $F_sw\in L^2(\Omega)$, there exists a unique solution
$s^\e\in\mathscr V_s^2$ to
\[
\left\{
\begin{aligned}
 &u_ss_x^\e+w_ss_Y^\e-\frac{2w_s}{r}s^\e-\e\cA s^\e=F_s,\\
 &s^\e|_{x=0}=s_x^\e|_{x=L}=0,\\
 &s^\e|_{Y=0}=0.
\end{aligned}
\right.
\]
The solution satisfies $\mathcal S_\e[s^\e]\lesssim\|F_sw\|$.
\end{proposition}

\begin{lemma}[The principal swirl problem]
\label{lem:swirl-principal-solvability}
Given $F_s$ with $F_sw\in L^2(\Omega)$, there is a unique
$s\in\mathscr V_s^2$ satisfying
\[
\left\{
\begin{aligned}
 &-\e\cA s=F_s,\\
 &s|_{x=0}=s_x|_{x=L}=0,\\
 &s|_{Y=0}=0.
\end{aligned}
\right.
\]
It satisfies
\begin{equation}\label{eq:swirl-principal-H2}
\begin{aligned}
 &\e^2\bigl(\|s_{xx}w\|^2+2\|s_{xY}w\|^2
                         +\|\cL s\,w\|^2\bigr)=\|F_sw\|^2,\\
 &\e\|s\|_{H^2(w)}\lesssim\|F_sw\|.
\end{aligned}
\end{equation}
The constant is independent of $0<L\le1$ and $\e>0$.
Thus $\|\cA s\,w\|\sim\|s\|_{H^2(w)}$ on
$\mathscr V_s^2$.
\end{lemma}

\begin{proof}
Choose $F_s^m\in C_c^\infty(\Omega)$ with
$\|(F_s^m-F_s)w\|\to0$. Use the same extension as in
Lemma~\ref{lem:main-principal-solvability}, and write on $\Omega$
\[
 F_s^m=\sum_{n\ {\rm odd}}F_{s,n}^m(Y)\sin(\lambda_nx),
 \qquad
 s^m=\sum_{n\ {\rm odd}}s_n^m(Y)\sin(\lambda_nx),
 \qquad \lambda_n=\frac{n\pi}{2L},
\]
where $n=1,3,5,\ldots$.
For each $n$, Lax--Milgram on
\[
 V=\{h:hw,h_Yw\in L^2(0,1),\ h(0)=0\}
\]
gives
\[
 \e\bigl((s_n^m)_Y,\phi_Yw^2\bigr)
 +\e\lambda_n^2(s_n^m,\phi w^2)
 =(F_{s,n}^m,\phi w^2),\qquad \phi\in V.
\]
The modal equation implies
\[
 \cL s_n^m=\lambda_n^2s_n^m-\e^{-1}F_{s,n}^m,
 \qquad \cL s_n^m\,w\in L^2(0,1).
\]
By Lemma~\ref{lem:near-axis-slice-estimates},
\begin{align*}
 &\bigl|\eta^3(s_n^m)_Y(1-\eta)s_n^m(1-\eta)\bigr|\\
 &\quad\lesssim
 \eta^2\|\cL s_n^m\,w\|\|s_n^m w\|
 +\eta\|\cL s_n^m\,w\|\|(s_n^m)_Yw\|
 \to0.
\end{align*}
Since $s_n^m(0)=0$, integration by parts gives
\[
 (\cL s_n^m,s_n^m w^2)=-\|(s_n^m)_Yw\|^2.
\]
Hence,
\begin{align*}
 &\|F_{s,n}^m w\|^2
 =\e^2\|\cL s_n^m\,w\|^2
 -2\e^2\lambda_n^2(\cL s_n^m,s_n^m w^2)
 +\e^2\lambda_n^4\|s_n^m w\|^2\\
 &=\e^2\bigl(\|\cL s_n^m\,w\|^2
 +2\lambda_n^2\|(s_n^m)_Yw\|^2
 +\lambda_n^4\|s_n^m w\|^2\bigr).
\end{align*}
Summing over $n$ and using Parseval's identity gives
\[
 \e^2\bigl(\|s^m_{xx}w\|^2+2\|s^m_{xY}w\|^2
                  +\|\cL s^m\,w\|^2\bigr)=\|F_s^m w\|^2.
\]
By \eqref{equiv-norm} and Poincar\'e's inequality,
\[
 \e\|s^m\|_{H^2(w)}\lesssim\|F_s^m w\|.
\]
Applying the same estimate to differences gives convergence in $H^2(w)$,
first as the number of modes tends to infinity and then as $m\to\infty$.
Passing to the limit proves \eqref{eq:swirl-principal-H2}.
The homogeneous energy identity gives uniqueness.
\end{proof}

\begin{proof}[Proof of Proposition~\ref{prop:swirl-solvability}]
Let $P_s s:=-\e\cA s$ with the boundary conditions in
\eqref{eq:swirl-equation}. By Lemma~\ref{lem:swirl-principal-solvability},
for every $h\in L^2(w)$ there is a unique $s=P_s^{-1}h$, and
\[
 \e\|P_s^{-1}h\|_{H^2(w)}\lesssim\|hw\|.
\]
Write \eqref{eq:swirl-equation} as
\[
 P_s s^\e+\mathcal K_s s^\e=F_s,\qquad
 \mathcal K_s s:=u_ss_x+w_ss_Y-\frac{2w_s}{r}s.
\]
Setting $h:=P_s s^\e$, we obtain
\[
 (I+\mathcal K_sP_s^{-1})h=F_s.
\]
For fixed $\e,L$, the profile bounds give a bounded map
$\mathcal K_s:\mathcal H^1(w)\to L^2(w)$.
Since $\mathcal H^2(w)\to\mathcal H^1(w)$ is compact on
$\Omega^\sharp$, $\mathcal K_sP_s^{-1}$ is compact on $L^2(w)$.

If $(I+\mathcal K_sP_s^{-1})h=0$ and $s^\e=P_s^{-1}h$, the
assumed a priori estimate gives $\mathcal S_\e[s^\e]=0$.
Hence $\nabla s^\e=0$, and the inlet condition gives $s^\e=h=0$.
The Fredholm alternative gives a unique $h\in L^2(w)$ and thus a
unique $s^\e=P_s^{-1}h$. The bound follows from the assumed estimate.
\end{proof}

\section{Linear remainder estimates}\label{sec:coupled-estimates}
We prove a priori estimates for
\eqref{eq:main-equation} and \eqref{eq:swirl-equation} with
$F_Gw,F_sw\in L^2(\Omega)$ and the stated boundary conditions.
Propositions~\ref{prop:main-solvability} and
\ref{prop:swirl-solvability} then give solvability. The coupling terms
are absorbed for sufficiently small $L$.

\subsection{Quotient estimates}

\begin{lemma}\label{lem:qxx}
Let $G$ solve \eqref{eq:main-equation}. For $0<\beta<\frac12$, we have
		\begin{equation}\label{qxx}
			\begin{aligned}
				&\varepsilon^2(\|\sqrt{u_s}q_{xxx}w\|^2+\| \sqrt{u_s} q_{xxY}w\|^2
				+\|\sqrt{u_s}\cL q_{x}w\|^2)\\
                &\quad+\varepsilon\|u_sq_{xY}w\|^2_{x=0}+\varepsilon\|u_sq_{xx}w\|^2_{x=L}+\varepsilon^{\frac{3}{2}}\|q_{xY}w\|_{Y=0}^2\\
				&\lesssim C_\beta L^{-\beta}\cD_2^2[q]
 +o_L(1)\cD_3^2[q]+\sqrt L\cD_4^2[G]
 +|(F_G,\varepsilon q_{xx}w^2)|.
			\end{aligned}
		\end{equation}
	\end{lemma}

\begin{proof}

\subsubsection*{\normalfont \textit{Step 1: Rayleigh terms}}

\begin{align*} \n
(-\p_x R[q], \e q_{xx}w^2) &\le  -c\Big(\|\sqrt\e u_s q_{xY} w\|_{x = 0}^2 + \|\sqrt\e u_s q_{xx}w\|_{x = L}^2\Big) +CL\bigl(\cD_2^2[q]+\cD_{3,x}^2[q]\bigr)
 + C\| \sqrt{\e}\nabla q_{x}w  \|^2.
\end{align*}

\begin{align*}
\begin{aligned}
-( \p_{xx} \{ u_s^2 q_x \},\e q_{xx}w^2) &= -\e (u_s^2 q_{xxx} + 4 u_s u_{sx} q_{xx} + 2[u_s u_{sxx} + u_{sx}^2] q_x, q_{xx}w^2)\\
 &\le - \frac{1}{2}\| \sqrt{\eps} u_s q_{xx} w \|_{x = L}^2
 + \frac{1}{2}\| \sqrt{\eps} u_s q_{xx} w \|_{x = 0}^2\\
&\quad + C\| \sqrt{\eps}u_s q_{xx} w \|^2
 +CL\bigl(\cD_2^2[q]+\cD_{3,x}^2[q]\bigr) \\
 &\le - \frac{1}{2}\| \sqrt{\eps} u_s q_{xx} w \|_{x = L}^2
 + C\| \sqrt{\eps}u_s q_{xx} w \|^2
 +CL\bigl(\cD_2^2[q]+\cD_{3,x}^2[q]\bigr),
\end{aligned}
\end{align*}

We use $|\partial_x^j u_s|\lesssim u_s$,
\eqref{b1} to control the inlet boundary term, and Poincar\'e's
inequality for the lower-order interior terms.

\begin{align*}
-( \frac{1}{r^3}\p_{xY} \{ u_s^2 q_{Y}r^3 \}, \e q_{xx} w^2) &= \e( \p_x \{ u_s^2 q_Y \}, q_{xxY} w^2)-\e(\p_{x} \{ u_s^2 q_{Y} \}, q_{xx}w^2)|_{Y=0}^{Y=1}\\
 &= \e( 2 u_s u_{sx} q_Y, q_{xxY} w^2) + \e( u_s^2 q_{xY}, q_{xxY} w^2) \\ 
  &= - \e( 2 u_s u_{sx}, q_{xY}^2 w^2) - \e( 2 \p_x \{ u_s u_{sx} \} q_{Y}, q_{xY} w^2) \\ 
&\quad - \e(  u_s u_{sx}, q_{xY}^2 w^2) +2 \e(u_s u_{sx} q_Y, q_{xY} w^2)_{x = L} \\
&\quad + \frac{1}{2}\|\sqrt{\e}u_s q_{xY} w\|_{x = L}^2 - \frac{1}{2}\|\sqrt{\e}u_s q_{xY} w\|_{x = 0}^2 \\
 &\le - \frac{1}{2}\|\sqrt{\e}u_s q_{xY} w\|_{x = 0}^2
 + C\|\sqrt{\e}u_s q_{xY} w\|_{x = L}^2
 +C\|\sqrt{\e}u_s q_{xY}w\|^2\\
 &\le - \frac{1}{2}\|\sqrt{\e}u_s q_{xY} w\|_{x = 0}^2
 + CL\bigl(\cD_2^2[q]+\cD_{3,x}^2[q]\bigr)
 +C\|\sqrt{\e}u_sq_{xY}w\|^2,
\end{align*}

where we used \eqref{b1} and Poincar\'e's inequality. The boundary
term vanishes at the wall by Theorem~\ref{thm:approximate-profiles}
and at the axis by \eqref{eq:near-axis}.
\subsubsection*{\normalfont \textit{Step 2: $\e\mathcal A^2$ Terms}}

\begin{align*} \n
(\e\mathcal{A}^2 G, \e q_{xx}w^2) &\le  - c\|\e \sqrt{u_s} \{ \cL q_{x}, q_{xxY},  q_{xxx} \} w \|^2 - c\e^\frac32 \|q_{xY}w\|^2_{Y=0}\\
&\quad +C_\beta L^{-\beta}\cD_2^2[q]+o_L(1) \bigl(\cD_2^2[q]+\cD_3^2[q]\bigr) + C\sqrt{L} \cD_4^2[G].
\end{align*}

Integrating the $G_{xxxx}$ term by parts in $x$ gives 

\begin{align}\label{Gxxxxqxx}
( \eps^2 G_{xxxx}, q_{xx}w^2) =  - (\eps^2 [ u_{sxxx}q + 3 u_{sxx}q_x + 3 u_{sx}q_{xx}  + u_s q_{xxx} ], q_{xxx} w^2) + (\eps^2 G_{xxx}, q_{xx}w^2)^{x = L}_{x = 0}.
\end{align}

The outlet boundary term vanishes because $G_{xxx}|_{x=L}=0$.
By \eqref{Poincare.2} with $k=\frac12$, the first three terms satisfy
\begin{align*}
&|(\ref{Gxxxxqxx}.{1, 2, 3})| \lesssim L \| \sqrt{u_s} \eps q_{xxx}  w \|^2 + \sqrt{\eps} L \| \sqrt{\eps} q_{xx}  w \|^2. 
\end{align*}
$(\ref{Gxxxxqxx}.{4})$ is favorable. For the boundary term, use
$G=u_sq$ and $G=G_{xx}=0$ at $x=0$ from \eqref{eq:main-equation}.

\begin{align*}
(\ref{Gxxxxqxx}.5) &= (2 \eps^2 G_{xxx}, \frac{u_{sx}}{u_s}  q_x w^2)_{x = 0} \\
 &\lesssim  \| \eps^{\frac{3}{2}} G_{xxx}  w \Big( \frac{u_{sx}}{u_s} \Big) \|_{x = 0} \| \sqrt{\eps} q_x  w \|_{x = 0} \\
 &\lesssim  \| \eps^{\frac{3}{2}} G_{xxx}  w ( \frac{u_{sx}}{u_s} ) \|_{x = 0} \sqrt{L} \bigl(\cD_2[q]+\cD_3[q]\bigr) \\
 &\lesssim  \sqrt{L}\cD_{4}^{\frac12}[G] \bigl(\cD_2[q]+\cD_{3,x}[q]\bigr)^{\frac{3}{2}}\\
 &\lesssim  \sqrt{L}\cD_{4}^2[G] +\sqrt{L} \bigl(\cD_2^2[q]+\cD_{3,x}^2[q]\bigr).
\end{align*}

We use \eqref{systemat1}, \eqref{b1}, and the trace inequality,
with $G_{xxx}|_{x=L}=0$.

For $2\e\partial_{xx}\cL$, integrate by
parts in $Y$, expand $G_{xxY}$, and integrate the $u_{sY}q_{xx}$
term once more in $Y$. This gives

\begin{align} \n
(2\eps \cL G_{xx},\e q_{xx}w^2) &=  (- 2\eps^2 G_{xxY}, q_{xxY}w^2)+(2\eps^2 G_{xxY},q_{xx}w^2)_{Y=0}^{Y=1}\\ \n
  &=  (-2\eps^2[u_{sxxY} q + 2 u_{sxY}q_x + u_{sxx}q_Y + u_s q_{xxY} + 2 u_{sx}q_{xY} + u_{sY} q_{xx}], q_{xxY}w^2)
\\ \n
 &=  -( 2\eps^2 [ u_{sxxY} q + 2 u_{sxY} q_x + u_{sxx} q_Y + u_s q_{xxY} \\ 
&\quad + 2 u_{sx} q_{xY} ], q_{xxY}w^2) + ( \eps^2 \cL u_{s}, q_{xx}^2 w^2)-(\eps^2u_{sY}q_{xx},q_{xx}w^2)_{Y=0}^{Y=1} \label{LGxxqxx}
\end{align}

The fourth term is coercive. For the first two terms, apply
\eqref{eq:usbound} and the weighted Hardy inequality
\eqref{eq:hardy-wall}:

\begin{equation}
\begin{aligned}\label{first separate}
|(\ref{LGxxqxx}.\{1, 2 \})| &\lesssim \|\e^\frac12 Y\frac{\{u_{sxxY},u_{sxY}\}}{\sqrt{u_s}}\|_{L^\infty}\|\e^\frac12\frac{\{q,q_x\}}{Y}w \|\|\e\sqrt{u_s} q_{xxY} w\|\\
      &\lesssim \|\e^\frac12\frac{\{q,q_x\}}{Y}w \|\|\e \sqrt{u_s}q_{xxY} w\|\\
           &\lesssim \|\e^\frac12 q_{xY} w\| \|\e\sqrt{u_s}q_{xxY} w\|\\
      &\lesssim L^{-\beta}\|\e^\frac12u_sq_{xY}w\|^2+ o_L(1)\bigl(\cD_2^2[q]+\cD_{3,x}^2[q]\bigr)  
\end{aligned}
     \end{equation}

      Poincar\'e's inequality controls the third and fifth terms.  
     
\begin{align*}
|(\ref{LGxxqxx}.\{ 3, 5 \})| &\lesssim \sqrt{\eps} \| \sqrt\e q_{xY} w \| \| \e\sqrt{u_s}  q_{xxY} w \| \\
|(\ref{LGxxqxx}.\{6\})| &\lesssim
\|\sqrt\e q_{xx}w\|^2. 
\end{align*}
 
The last boundary term vanishes by \eqref{axis-vanish} at the
axis and by $q=0$ at the wall.

For $\cL^2$, integration by parts gives

\begin{align} \n
(\e^2\cL^2 G, q_{xx}w^2) &=  - (\e^2(\cL G)_Y, q_{xxY}w^2) + (\e^2(\cL G)_Y, q_{xx}w^2)_{Y=1} \\ \n
 &=  (\e^2(\cL G_x)_Y, q_{xY}w^2) - (\e^2(\cL G)_Y, q_{xY}w^2)_{x=L} \\\label{L2qxx1} 
 &=  - (\e^2\cL G_{x}, \cL q_{x}w^2) - (\e^2\cL G_{x}, q_{xY}w^2)_{Y = 0} - (\e^2(\cL G)_Y, q_{xY}w^2)_{x = L}.
\end{align}

The expansion
\begin{align*}
    &\cL G_x = u_{sx}\cL q+ u_s\cL q_x +\cL u_{sx} q+ \cL u_s q_x +2u_{sxY}q_Y+2u_{sY}q_{xY}
\end{align*}
 yields
 \begin{align*}
     &(\ref{L2qxx1}.1) = -(\e^2[u_{sx}\cL q+ u_s\cL q_x +\cL u_{sx} q+ \cL u_s q_x +2u_{sxY}q_Y+2u_{sY}q_{xY}],\cL q_x w^2)
 \end{align*}

The second term is coercive. Poincar\'e's inequality gives
$|(\ref{L2qxx1}.1.1)|\lesssim L\|\e\sqrt{u_s}\cL q_xw\|^2$,
so the first term is absorbed for small $L$. Estimate the third
and fourth terms as in \eqref{first separate}.

\begin{align*}
           |(\ref{L2qxx1}.1.3+4)| &\lesssim \sqrt\e
 \left\|Y\frac{\{\cL u_{sx},\cL u_s\}}{\sqrt{u_s}}\right\|_\infty
 \left\|\sqrt\e\frac{\{q,q_x\}}Yw\right\|
 \|\e\sqrt{u_s}\cL q_xw\|\\
  &\lesssim \|\sqrt\e q_{xY}w\|\|\e\sqrt{u_s}\cL q_xw\|\\
  &\leq C_\beta L^{-\beta}\cD_2^2[q]+o_L(1)\cD_{3,x}^2[q].
\end{align*}

We used \eqref{interp.L}. Integration by parts in $Y$ gives

\begin{align*}
(\ref{L2qxx1}.1.5) &= ( 2\e^2 \bar\cL u_{sx} q_{xY}, q_Y w^2) + (2\e^2 u_{sxY} q_{xY}, \cL q w^2) + (2\e^2 u_{sxY} q_Y, q_{xY}w^2)_{Y= 0}
\\
 &\lesssim  \| \sqrt{\e} q_{xY} w \|\| \e \cL q w \|+ L \|\sqrt{\e}q_{xY} w \|^2  + L \|\e \sqrt{u_{sY}} q_{xY} w \|_{Y = 0}^2\\
 &\lesssim \sqrt L\bigl(\cD_2^2[q]+\cD_3^2[q]\bigr).
\end{align*}

We use \eqref{small.guys} and \eqref{b1}.
Integration by parts in $Y$ gives

\begin{align*} \n
(\ref{L2qxx1}.1.6) &=  (\e^2 q_{xY}^2, \bar\cL u_{s} w^2) +(\e^2 q_{xY}^2, u_{sY} w^2)_{Y = 0} \\  \n
 &\leq  C_\beta L^{-\beta}\cD_2^2[q]+o_L(1)\cD_{3,x}^2[q] + \|\e\sqrt{u_{sY}}q_{xY} w\|_{Y = 0}^2. \\ 
\end{align*}

The wall term $(\ref{L2qxx1}.2)$ satisfies

\begin{align*}
(\ref{L2qxx1}.2) &=  - (\e^2\{u_s \cL q_x + \cL u_{sx}q  + \cL u_s q_x + u_{sx} \cL q + 2u_{sxY}q_Y + 2u_{sY} q_{xY}\} ,q_{xY}w^2)_{Y = 0}
\\
 &\le - (2 - C_0 L) \|\e\sqrt{u_{sY}} q_{xY}w\|^2_{Y = 0}, 
\end{align*}

where $C_0<\infty$ is independent of $L$ and $\e$.
The first four terms vanish because $q=q_x=u_s=u_{sx}=0$ at
the wall. The sixth term equals $-2\|\e\sqrt{u_{sY}}q_{xY}w\|_{Y=0}^2$,
and the wall shear satisfies

\begin{align*}
 &u_{sY}(x,0)
 =\e^{-\frac{1}{2}}\partial_y\widetilde u_p^0(x,0)+O(1)
 \ge c\e^{-\frac{1}{2}}.
\end{align*}

For the fifth term, use Poincar\'e's inequality and
$|\frac{u_{sxY}}{u_{sY}}|_{Y=0}\lesssim1$:

\begin{align*}
&|(\ref{L2qxx1}.2.5)| \lesssim L \sup_{x\in(0,L)}| \frac{u_{sxY}}{u_{sY}}| _{Y = 0} \|\e \sqrt{u_{sY}}q_{xY} w \|_{Y = 0}^2
\end{align*} 
The wall terms thus have a negative coefficient of
$\|\e\sqrt{u_{sY}}q_{xY}w\|_{Y=0}^2$ for small $L$.

Since $G|_{x=0}=0$, the fundamental theorem of calculus in $x$ gives 

\begin{align*}
(\ref{L2qxx1}.3) &=  (\e^2(\cL G)_Y, \p_Y \{ \frac{u_{sx}}{u_s} q \} w^2)_{x = L} \\
 &=  (\e^2(\cL G_x)_Y, \p_Y \{ \frac{u_{sx}}{u_s} q \} w^2) + (\e^2(\cL G)_Y, \p_{xY} \{ \frac{u_{sx}}{u_s} q \} w^2) \\
 &=  - (\e^2\cL G_{x}, \cL \{ \frac{u_{sx}}{u_s} q \} w^2) - (\e^2\cL G_{x}, \frac{u_{sx}}{u_s} q_Y w^2)_{Y = 0} + (\e^2(\cL G)_Y, \p_{xY} \{ \frac{u_{sx}}{u_s} q \} w^2)
\end{align*}

Expanding the first term and using \eqref{eq:profile-quotient-bounds} gives

\begin{align*}
|(\ref{L2qxx1}.3.1)| &=  |- (\e^2\cL G_{x}, [\cL \{ \frac{u_{sx}}{u_s} \}q + 2 \p_Y \{ \frac{u_{sx}}{u_s} \} q_Y + \frac{u_{sx}}{u_s} \cL q] w^2)| \\
  &\lesssim  \| \e\cL G_{x} w \| \Big[ \|\sqrt{\e}\cL \{ \frac{u_{sx}}{u_s}\} Y \|_\infty \|\sqrt{\e} \frac{q}{Y} w \| + L \|\sqrt{\e} \p_Y \{ \frac{u_{sx}}{u_s} \} \|_\infty \|\sqrt{\e} q_{xY} w \| + \|\frac{u_{sx}}{u_s} \|_\infty \| \e \cL q w \| \Big]
\\
 &\lesssim  \|\e \cL G_{x} w\| \Big[  L \| \sqrt{\e}q_{xY} w \| + \|\e \cL q w \| \Big] \\
 &\lesssim  o_L(1) \bigl(\cD_2^2[q]+\cD_3^2[q]\bigr). 
\end{align*}

In the boundary term, expand $\cL G_x$:

\begin{align*}
(\ref{L2qxx1}.3.2) &= - \e^2(u_s \cL q_{x} + u_{sx} \cL q + 2 u_{sxY} q_{Y} + 2 u_{sY} q_{xY} + \cL u_{sx} q + \cL u_{s}q_x, \frac{u_{sx}}{u_s}q_Y w^2)_{Y = 0}.
\end{align*}

The first two terms vanish because $u_s=u_{sx}=0$ at the wall,
and the fifth and sixth because $q=q_x=0$ there. The remaining terms satisfy

\begin{align*}
&|(\ref{L2qxx1}.3.2.3)|+|(\ref{L2qxx1}.3.2.4)|
 \lesssim L\|\e\sqrt{u_{sY}}q_{xY}w\|_{Y=0}^2
 \lesssim L\cD_3^2[q].
\end{align*}

The last term satisfies

\begin{align*}
|(\ref{L2qxx1}.3.3)| &\lesssim\|\e^{\frac{3}{2}}(\cL G)_Yw\|\|\sqrt\e q_{xY}w\|\\
  &\leq C_\beta L^{-\beta}\cD_2^2[q]+o_L(1)\cD_3^2[q].
\end{align*}

\subsubsection*{\normalfont \textit{Step 3: $J$ terms}}

Using $r(\beta_s)_x=w_{sxx}-u_{sxY}$, and $r(\beta_s)_Y = -u_{sxx}-\cL u_s-\frac2r u_{sY}+2\beta_s$, we expand the pairing with $-\e q_{xx}w^2$ as

\begin{align}\label{Jqxx}
&(-u_{sx} G_{xx}-u_{sx} \cL G
+G(\cL u_{sx})+G u_{sxxx}
+w_s(\cL G)_Y+w_sG_{xxY} \n \\*
&\quad +2G \beta_{sx}-w_{sxx}G_Y+u_{sxY}G_Y
-\frac2r u_{sY}G_x+2\beta_sG_x,-\e q_{xx}w^2) .
\end{align}

\begin{align*}
&(\ref{Jqxx}.1,4)\lesssim \|\sqrt\e G_{xx}w\|\|\sqrt\e q_{xx}w\|\lesssim L^{-\beta}\|\sqrt\e u_s q_{xx} w\|^2+o_L(1)\bigl(\cD_2^2[q]+\cD_{3,x}^2[q]\bigr). 
\end{align*}

Integration by parts in $x$ and $Y$ gives

\begin{align*}
    (\ref{Jqxx}.2)+(\ref{Jqxx}.9) &= -(u_{sx}G_Y,\e q_{xxY}w^2)
    -2(u_{sxY}G_Y,\e q_{xx}w^2)\\
     &= (u_{sxx}G_Y,\e q_{xY}w^2)
    +(u_{sx}G_{xY},\e q_{xY}w^2)
    -(u_{sx}G_{Y},\e q_{xY}w^2)_{x=L}\\
    &\quad +2(u_{sxxY}G_Y,\e q_xw^2)+2(u_{sxY}G_{xY},\e q_x w^2)-2(\e u_{sxY}G_Y,q_xw^2)_{x=L}\\
     &\lesssim (1+L) \|\sqrt\e G_{xY}w\|\|\sqrt\e q_{xY}w\|+\|\sqrt\e G_{Y}w\|_{x=L}\|\sqrt\e q_{xY}w\|_{x=L}\\
    &\quad +\|Y\{u_{sxxY},u_{sxY}\}\|_{\infty}\|\sqrt\e \{G_{Y},G_{xY}\}w\|\|\sqrt\e \frac{q_{x}}{Y}w\|\\
    &\quad +\|Y u_{sxY}\|_{\infty}\|\sqrt\e G_Y w\|_{x=L}\|\sqrt\e \frac{q_x}{Y}w\|_{x=L}\\
 &\lesssim L^{-\beta}\|\sqrt\e u_s q_{xY} w\|^2+o_L(1)\bigl(\cD_2^2[q]+\cD_{3,x}^2[q]\bigr).     
\end{align*}

\begin{align*}
    (\ref{Jqxx}.3) &= (u_{sxY}G_Y,\e q_{xx}w^2)+(u_{sxY}G,\e q_{xxY}w^2)
    \\
     &= -(u_{sxxY}G_Y,\e q_{x}w^2)
    -(u_{sxY}G_{xY},\e q_{x}w^2)
    -(u_{sxxY}G_{},\e q_{xY}w^2)\\
    &\quad -(u_{sxY}G_{x},\e q_{xY}w^2)
    +(u_{sxY}G_{Y},\e q_{x}w^2)_{x=L}
    +(u_{sxY}G_{},\e q_{xY}w^2)_{x=L}\\
     &\lesssim (1+L)\|\sqrt\e G_{xY}w\|\|\sqrt\e q_{xY}w\|\\
    &\quad +\|Y u_{sxY}\|_{\infty}\left(\|\sqrt\e G_{Y}w\|_{x=L}\|\sqrt\e \frac{q_{x}}{Y}w\|_{x=L}
    +\|\sqrt\e \frac{G}{Y}w\|_{x=L}\|\sqrt\e q_{xY}w\|_{x=L}\right)\\
 &\lesssim L^{-\beta}\|\sqrt\e u_s q_{xY} w\|^2+o_L(1)\bigl(\cD_2^2[q]+\cD_{3,x}^2[q]\bigr).     
\end{align*}

\begin{align*}
    (\ref{Jqxx}.5) &=\e((w_{sY}-3\frac{w_s}{r})\cL G,q_{xx}w^2)
   +\e(w_s\cL G,q_{xxY}w^2)\\
  &=-\e(\partial_x\{(w_{sY}-3\frac{w_s}{r})\cL G\},q_xw^2)
   -\e(\partial_x\{w_s\cL G\},q_{xY}w^2)\\
 &\quad +\e((w_{sY}-3\frac{w_s}{r})\cL G,q_xw^2)_{x=L}
   +\e(w_s\cL G,q_{xY}w^2)_{x=L}\\
  &\lesssim\|\e\cL Gw\|
        (\|\sqrt\e \frac{q_x}{Y}w\|+\|\sqrt\e q_{xY}w\|)
 +\|\e\cL G_xw\|
        (\|\sqrt\e q_xw\|+\|\sqrt\e q_{xY}w\|)\\
 &\quad +\|\e\cL Gw\|_{x=L}
        (\|\sqrt\e q_{xY}w\|_{x=L}
                    +\|\sqrt\e \frac{q_x}{Y}w\|_{x=L})\\
  &\leq C_\beta L^{-\beta}\cD_2^2[q]+o_L(1)\cD_3^2[q].
\end{align*}

Theorem~\ref{thm:approximate-profiles} bounds the remaining coefficients.

\begin{align*}
(\ref{Jqxx}.7,10,11) &\lesssim \|Y\{\beta_{sx},\beta_s,\frac{u_{sY}}{r}\}\|_{\infty}\|\sqrt\e\frac{\{G,G_x\}}{Y}w\|\|\sqrt\e q_{xx}w\|\\*
 &\lesssim L^{-\beta}\|\sqrt\e u_s q_{xx} w\|^2+o_L(1)\bigl(\cD_2^2[q]+\cD_{3,x}^2[q]\bigr). 
\end{align*}

The remaining terms are bounded by \eqref{systemat1} and \eqref{interp.1}:

\begin{align*}
 |(\ref{Jqxx}.6)|+|(\ref{Jqxx}.8)| &\lesssim \|\e G_{xxY}w\|\|\sqrt\e q_{xx}w\|
  +\sqrt\e\|\sqrt\e G_Yw\|\|\sqrt\e q_{xx}w\|\\
  &\leq C_\beta L^{-\beta}\cD_2^2[q]+o_L(1)\cD_3^2[q].
\end{align*}

\end{proof}

\begin{lemma}
Let \(G\) solve \eqref{eq:main-equation}. Then

\begin{align}\label{Lq-est}
&\e^3
\Big(
\|\sqrt{u_s}(\cL q)_Yw\|^2
+
\|\sqrt{u_s}\cL q_xw\|^2
+
\|\sqrt{u_s}q_{xxY}w\|^2
\Big)
+
\e^2\|u_s\cL q\,w\|_{x=L}^2
+
\e^{\frac{5}{2}}\|\cL q\,w\|_{Y=0}^2
\nonumber\\
 &\lesssim
o_L(1)\bigl(\cD_2^2[q]+\cD_3^2[q]\bigr)
+
|(F_G,-\e^2\cL q\,w^2)|.
\end{align}

\end{lemma}

\begin{proof}
\subsubsection*{\normalfont \textit{Step 1: Rayleigh terms}}

\begin{align*}
&(-\p_x R[q],\e^2\cL q\,w^2)
\le
-c\e^2\|u_s\cL q\,w\|_{x=L}^2
+
o_L(1)\bigl(\cD_2^2[q]+\cD_3^2[q]\bigr).
\end{align*}

Using $\frac1{r^3}\p_Y\{r^3u_s^2q_Y\} = u_s^2\cL q+2u_su_{sY}q_Y$, we expand

\begin{align}\label{RLq}
\begin{aligned}
(-\p_x R[q],\e^2\cL q\,w^2) &=
-(\p_{xx}\{u_s^2q_x\},\e^2\cL q\,w^2)
\\
&\quad -(\p_x\{u_s^2\cL q\},\e^2\cL q\,w^2)
\\
&\quad -(\p_x\{2u_su_{sY}q_Y\},\e^2\cL q\,w^2).
\end{aligned}
\end{align}

The term \((\ref{RLq}.1)\) satisfies

\begin{align*}
\begin{aligned}
(\ref{RLq}.1) &=
-\e^2
\left(
u_s^2q_{xxx}
+
4u_su_{sx}q_{xx}
+
2[u_su_{sxx}+u_{sx}^2]q_x,
\cL q\,w^2
\right)
\\
 &\lesssim
\|\e\cL q\,w\|
\left(
        \|\e\sqrt{u_s}q_{xxx}w\|
        +
        \|\sqrt\e q_{xx}w\|
        +
        \|\sqrt\e q_xw\|
\right)
\\
 &\lesssim
o_L(1)\bigl(\cD_2^2[q]+\cD_3^2[q]\bigr).
\end{aligned}
\end{align*}

For \((\ref{RLq}.2)\), integrate by parts in \(x\):

\begin{align*}
\begin{aligned}
(\ref{RLq}.2) &=
-\e^2
\left(
u_s^2(\cL q)_x
+
2u_su_{sx}\cL q,
\cL q\,w^2
\right)
\\
 &=
-\frac{\e^2}{2}\|u_s\cL q\,w\|_{x=L}^2
+
\frac{\e^2}{2}\|u_s\cL q\,w\|_{x=0}^2
-
\e^2
\left(
u_su_{sx}\cL q,
\cL q\,w^2
\right)
\\
 &\le
-\frac{\e^2}{2}\|u_s\cL q\,w\|_{x=L}^2
+
C\|\e\cL q\,w\|^2
\\
 &\le
-\frac{\e^2}{2}\|u_s\cL q\,w\|_{x=L}^2
+
o_L(1)\bigl(\cD_2^2[q]+\cD_3^2[q]\bigr).
\end{aligned}
\end{align*}

Similarly, \((\ref{RLq}.3)\) satisfies

\begin{align*}
\begin{aligned}
(\ref{RLq}.3) &=
-\e^2
\left(
2\p_x\{u_su_{sY}\}q_Y
+
2u_su_{sY}q_{xY},
\cL q\,w^2
\right)
\\
 &\lesssim
\|\e\cL q\,w\|
\left(
        \|\sqrt\e q_Yw\|
        +
        \|\sqrt\e q_{xY}w\|
\right)
\\
 &\lesssim
o_L(1)\bigl(\cD_2^2[q]+\cD_3^2[q]\bigr).
\end{aligned}
\end{align*}

\subsubsection*{\normalfont \textit{Step 2: \(\e\mathcal A^2\) Terms}}

\begin{align*}
\begin{aligned}
(\e\cA^2G,\e^2\cL q\,w^2) &\le
-c\e^3
\left(
        \|\sqrt{u_s}(\cL q)_Yw\|^2
        +
        \|\sqrt{u_s}\cL q_xw\|^2
        +
        \|\sqrt{u_s}q_{xxY}w\|^2
\right)
\\
&\quad -c\e^{\frac52}\|\cL q\,w\|_{Y=0}^2
+
o_L(1)\bigl(\cD_2^2[q]+\cD_3^2[q]\bigr).
\end{aligned}
\end{align*}

Expanding the viscous term gives

\begin{align}\label{ALq}
\begin{aligned}
(\e\cA^2G,\e^2\cL q\,w^2) &=
(\e^3G_{xxxx},\cL q\,w^2) +
(2\e^3\cL G_{xx},\cL q\,w^2) +
(\e^3\cL^2G,\cL q\,w^2).
\end{aligned}
\end{align}

For \((\ref{ALq}.1)\), integrate twice in \(x\) and once in
\(Y\):

\begin{align}\label{GxxxxLq}
\begin{aligned}
(\ref{ALq}.1) &=
(\e^3G_{xx},\cL q_{xx}w^2)
-
(\e^3G_{xx},\cL q_xw^2)_{x=L}
\\
 &=
-\e^3(G_{xxY},q_{xxY}w^2)
-
(\e^3G_{xx},\cL q_xw^2)_{x=L}
\\
 &=
-\e^3(u_sq_{xxY},q_{xxY}w^2)
-\e^3(u_{sY}q_{xx},q_{xxY}w^2)
\\
&\quad -\e^3(2u_{sx}q_{xY},q_{xxY}w^2)
-\e^3(2u_{sxY}q_x,q_{xxY}w^2)
\\
&\quad -\e^3(u_{sxx}q_Y,q_{xxY}w^2)
-\e^3(u_{sxxY}q,q_{xxY}w^2)
\\
&\quad -
(\e^3G_{xx},\cL q_xw^2)_{x=L}.
\end{aligned}
\end{align}

\((\ref{GxxxxLq}.1)\) is favorable: $(\ref{GxxxxLq}.1) = -\e^3\|\sqrt{u_s}q_{xxY}w\|^2$. For \((\ref{GxxxxLq}.2)\), integrate by parts in \(Y\):

\begin{align*}
\begin{aligned}
(\ref{GxxxxLq}.2) &=
\frac{\e^3}{2}(\cL u_s,q_{xx}^2w^2)
-
\frac{\e^3}{2}(u_{sY}q_{xx}^2w^2)_{Y=0}^{Y=1}
\\
 &\lesssim
\e^3\|\cL u_s\|_{L^\infty}\|q_{xx}w\|^2
\\
 &\lesssim
\e\|\sqrt\e q_{xx}w\|^2
\\
 &\lesssim
\e\bigl(\cD_2^2[q]+\cD_3^2[q]\bigr).
\end{aligned}
\end{align*}

The boundary term vanishes at \(Y=0,1\).

For \((\ref{GxxxxLq}.3)+(\ref{GxxxxLq}.5)\),

\begin{align*}
\begin{aligned}
|(\ref{GxxxxLq}.3)+(\ref{GxxxxLq}.5)| &\lesssim
\e^\frac32
\left(
        \|\sqrt\e q_{xY}w\|
        +
        \|\sqrt\e q_Yw\|
\right)
\|\e\sqrt{u_s}q_{xxY}w\|
\\
 &\lesssim
\e^\frac32\bigl(\cD_2^2[q]+\cD_{3,x}^2[q]\bigr).
\end{aligned}
\end{align*}

For \((\ref{GxxxxLq}.4)+(\ref{GxxxxLq}.6)\),

\begin{align*}
\begin{aligned}
|(\ref{GxxxxLq}.4)+(\ref{GxxxxLq}.6)| &\lesssim
\e^{\frac{3}{2}}\left\|
Y\frac{\{u_{sxY},u_{sxxY}\}}{\sqrt{u_s}}
\right\|_{L^\infty}
\|\sqrt\e\frac{\{q_x,q\}}{Y}w\|
\|\e\sqrt{u_s}q_{xxY}w\|
\\
 &\lesssim
\e^\frac32\bigl(\cD_2^2[q]+\cD_{3,x}^2[q]\bigr).
\end{aligned}
\end{align*}

For \((\ref{GxxxxLq}.7)\), use $G_{xx}=u_sq_{xx}+2u_{sx}q_x+u_{sxx}q$:

\begin{align*}
\begin{aligned}
|(\ref{GxxxxLq}.7)| &\le
\left|
(\e^3u_sq_{xx},\cL q_xw^2)_{x=L}
\right|
+
\left|
(2\e^3u_{sx}q_x,\cL q_xw^2)_{x=L}
\right| +
\left|
(\e^3u_{sxx}q,\cL q_xw^2)_{x=L}
\right|
\\
 &\lesssim
 \e\|\e\sqrt{u_s}\{q_{xx},q_x,q\}w\|_{x=L}
     \|\e\sqrt{u_s}\cL q_xw\|_{x=L}\\
 &\lesssim\e L\bigl(\cD_2^2[q]+\cD_{3,x}^2[q]\bigr).
\end{aligned}
\end{align*}

For \((\ref{ALq}.2)\), integrate by parts in \(x\):

\begin{align}\label{LGxxLq}
\begin{aligned}
(\ref{ALq}.2) &=
-2\e^3(\cL G_x,\cL q_xw^2)
\\
 &=
-2\e^3(u_s\cL q_x,\cL q_xw^2)
-2\e^3(u_{sx}\cL q,\cL q_xw^2)
\\
&\quad -2\e^3(\cL u_{sx}\,q,\cL q_xw^2)
-2\e^3(\cL u_s\,q_x,\cL q_xw^2)
\\
&\quad -2\e^3(2u_{sxY}q_Y,\cL q_xw^2)
-2\e^3(2u_{sY}q_{xY},\cL q_xw^2).
\end{aligned}
\end{align}

\((\ref{LGxxLq}.1)\) is favorable: $(\ref{LGxxLq}.1) = -2\e^3\|\sqrt{u_s}\cL q_xw\|^2$. For \((\ref{LGxxLq}.2)\),

\begin{align*}
\begin{aligned}
|(\ref{LGxxLq}.2)| &\lesssim
\e \|\e\cL q\,w\|
\|\e\sqrt{u_s}\cL q_xw\|
\\
 &\lesssim
\e\bigl(\cD_2^2[q]+\cD_3^2[q]\bigr).
\end{aligned}
\end{align*}

For \((\ref{LGxxLq}.3)+(\ref{LGxxLq}.4)\),

\begin{align*}
\begin{aligned}
|(\ref{LGxxLq}.3)+(\ref{LGxxLq}.4)| &\lesssim
\e\left\|
\e^{\frac12}
Y\frac{\{\cL u_{sx},\cL u_s\}}{\sqrt{u_s}}
\right\|_{L^\infty}
\left\|
\sqrt\e\frac{\{q,q_x\}}{Y}w
\right\|
\left\|
\e\sqrt{u_s}\cL q_xw
\right\|
\\
 &\lesssim
\e\bigl(\cD_2^2[q]+\cD_3^2[q]\bigr).
\end{aligned}
\end{align*}

For \((\ref{LGxxLq}.5)\), integrate by parts in \(Y\):

\begin{align*}
\begin{aligned}
|(\ref{LGxxLq}.5)| &\lesssim
\e
\left|
(2\e^2\bar\cL u_{sx}q_{xY},q_Yw^2)
\right|
+
\e
\left|
(2\e^2u_{sxY}q_{xY},\cL q\,w^2)
\right| +
\left|\e^3(2u_{sxY}q_{xY},q_Y\,w^2)\right|_{Y=0}
\\
 &\lesssim
\e L\|\sqrt\e q_{xY}w\|^2
+
\e\|\sqrt\e q_{xY}w\|\|\e\cL q\,w\|
+
\e^{\frac{5}{2}}L\|q_{xY}w\|_{Y=0}^2
\\
 &\lesssim
\e\bigl(\cD_2^2[q]+\cD_3^2[q]\bigr).
\end{aligned}
\end{align*}

For \((\ref{LGxxLq}.6)\), integrate by parts in \(Y\):

\begin{align*}
\begin{aligned}
(\ref{LGxxLq}.6) &=
2\e^3(\bar\cL u_s,q_{xY}^2w^2)
+
2\e^3(u_{sY}q_{xY}^2w^2)_{Y=0}
\\
 &\lesssim
\e\|\sqrt\e q_{xY}w\|^2
+
\e^{\frac52}\|q_{xY}w\|_{Y=0}^2
\\
 &\lesssim
\e\bigl(\cD_2^2[q]+\cD_3^2[q]\bigr).
\end{aligned}
\end{align*}

For \((\ref{ALq}.3)\), the self-adjoint form of \(\cL\) gives
\begin{align*}
\begin{aligned}
&(\ref{ALq}.3)
=
-\e^3((\cL G)_Y,(\cL q)_Yw^2)
+
(\e^3(\cL G)_Y,\cL q\,w^2)_{Y=0}^{Y=1}.
\end{aligned}
\end{align*}
Expanding $\cL G = u_s\cL q+2u_{sY}q_Y+q\cL u_s$ gives
\begin{align*}
\begin{aligned}
&(\cL G)_Y
=
u_s(\cL q)_Y
+
3u_{sY}\cL q
+
(2u_{sYY}+\bar\cL u_s)q_Y
+
q(\cL u_s)_Y.
\end{aligned}
\end{align*}
Therefore

\begin{align}\label{L2GLq}
\begin{aligned}
(\ref{ALq}.3) &=
-\e^3(u_s(\cL q)_Y,(\cL q)_Yw^2)-
3\e^3(u_{sY}\cL q,(\cL q)_Yw^2)
\\
&\quad -
\e^3
\left(
(2u_{sYY}+\bar\cL u_s)q_Y,
(\cL q)_Yw^2
\right)
-
\e^3
\left(
q(\cL u_s)_Y,
(\cL q)_Yw^2
\right)
\\
&\quad +
(\e^3(\cL G)_Y,\cL q\,w^2)_{Y=0}^{Y=1}.
\end{aligned}
\end{align}

\((\ref{L2GLq}.1)\) is favorable: $(\ref{L2GLq}.1) = -\e^3\|\sqrt{u_s}(\cL q)_Yw\|^2$. For \((\ref{L2GLq}.2)\), integrate by parts in \(Y\):
\begin{align*}
\begin{aligned}
&(\ref{L2GLq}.2)
=
\frac{3}{2}\e^3(\cL u_s,(\cL q)^2w^2)
+
\frac{3}{2}\e^3(u_{sY}(\cL q)^2w^2)_{Y=0}.
\end{aligned}
\end{align*}
Since the boundary term at $Y=1$ vanishes,
\begin{align*}
\begin{aligned}
(\ref{L2GLq}.5)
=
-3\e^3(u_{sY}(\cL q)^2w^2)_{Y=0}
+
\e^3
\left(
(2u_{sYY}+\bar\cL u_s)q_Y,
\cL q\,w^2
\right)_{Y=0}^{Y=1}.
\end{aligned}
\end{align*}
The wall contribution from \((\ref{L2GLq}.2)+(\ref{L2GLq}.5)\) is $-\frac{3}{2}\e^3(u_{sY}(\cL q)^2w^2)_{Y=0}$. The interior coefficient term satisfies

\begin{align*}
\e^3(\cL u_s,(\cL q)^2w^2)
\lesssim
\|\e\cL q\,w\|^2
\lesssim
L\bigl(\cD_2^2[q]+\cD_3^2[q]\bigr).
\end{align*}

For \((\ref{L2GLq}.3)\) and the corresponding boundary contribution in
\((\ref{L2GLq}.5)\), integration by parts in \(Y\) gives

\begin{align*}
\begin{aligned}
&-\e^3
\left(
(2u_{sYY}+\bar\cL u_s)q_Y,
(\cL q)_Yw^2
\right)
+
\e^3
\left(
(2u_{sYY}+\bar\cL u_s)q_Y,
\cL q\,w^2
\right)_{Y=0}^{Y=1}
\\
 &=
\e^3
\left(
(2u_{sYY}+\bar\cL u_s)\cL q
+
(2u_{sYYY}+(\bar\cL u_s)_Y)q_Y,
\cL q\,w^2
\right)
\\
 &\lesssim
\|\e\cL q\,w\|^2
+
\|\sqrt\e q_Yw\|\|\e\cL q\,w\|
\\
 &\lesssim
o_L(1)\bigl(\cD_2^2[q]+\cD_3^2[q]\bigr).
\end{aligned}
\end{align*}

For \((\ref{L2GLq}.4)\),

\begin{align*}
\begin{aligned}
(\ref{L2GLq}.4) &\lesssim
\left\|
\e
Y\frac{(\cL u_s)_Y}{\sqrt{u_s}}
\right\|_{L^\infty}
\left\|
\sqrt\e\frac{q}{Y}w
\right\|
\left\|
\e^{\frac32}\sqrt{u_s}(\cL q)_Yw
\right\|
\\
 &\le
\eta\e^3\|\sqrt{u_s}(\cL q)_Yw\|^2
+
C_\eta o_L(1)\bigl(\cD_2^2[q]+\cD_3^2[q]\bigr).
\end{aligned}
\end{align*}

For sufficiently small \(\eta>0\), \((\ref{ALq}.3)\) satisfies

\begin{align*}
\begin{aligned}
(\ref{ALq}.3) &\le
-c\e^3\|\sqrt{u_s}(\cL q)_Yw\|^2
-\frac{3}{2}\e^3(u_{sY}(\cL q)^2w^2)_{Y=0} +
o_L(1)\bigl(\cD_2^2[q]+\cD_3^2[q]\bigr).
\end{aligned}
\end{align*}

\subsubsection*{\normalfont \textit{Step 3: $J$ terms}}

\begin{align*}
&(-u_{sx} G_{xx}-u_{sx} \cL G
+G(\cL u_{sx})+G u_{sxxx}
+w_s(\cL G)_Y+w_sG_{xxY} \n \\*
&\quad +2G \beta_{sx}-w_{sxx}G_Y+u_{sxY}G_Y
-\frac2r u_{sY}G_x+2\beta_sG_x,-\e^2 \cL q w^2)\lesssim\sqrt L\bigl(\cD_2^2[q]+\cD_3^2[q]\bigr).
\end{align*}

The profile bounds and \eqref{systemat1} give
$\|\varepsilon\mathcal J[G]w\|\lesssim\mathcal D_2[q]+\mathcal D_3[q]$.
Cauchy's inequality and \eqref{small.guys} give the result.

\end{proof}

\begin{lemma}
Let \(G\) solve \eqref{eq:main-equation}, and let
\(w=r^{\frac32}\). Then
\begin{equation}\label{qx-est}
\|\sqrt\e u_s q_{xY}w\|^2
+
\|\sqrt\e u_s q_{xx}w\|^2
\lesssim
\sqrt L\left(
        \bigl(\cD_2^2[q]+\cD_3^2[q]\bigr)+\cD_4^2[G]
\right)
+
|(F_G,\e q_xw^2)|.
\end{equation}
\end{lemma}

\begin{proof}
\subsubsection*{\normalfont \textit{Step 1: Rayleigh terms}}

\begin{align*}
\begin{aligned}
&(-\p_x R[q],\e q_xw^2)
\ge
c\Big(\|\sqrt\e u_sq_{xY}w\|^2
+
\|\sqrt\e u_sq_{xx}w\|^2\Big)
-
CL\bigl(\cD_2^2[q]+\cD_3^2[q]\bigr).
\end{aligned}
\end{align*}

Using $R[q] = \p_x\{u_s^2q_x\} + \frac1{r^3}\p_Y\{r^3u_s^2q_Y\}$, we expand

\begin{align}\label{Rqx}
\begin{aligned}
(-\p_x R[q],\e q_xw^2) &=
-(\p_{xx}\{u_s^2q_x\},\e q_xw^2) -
\left(
        \frac1{r^3}\p_{xY}\{r^3u_s^2q_Y\},
        \e q_xw^2
\right).
\end{aligned}
\end{align}

For \((\ref{Rqx}.1)\), integrate by parts in \(x\):

\begin{align}\label{RqxX}
\begin{aligned}
(\ref{Rqx}.1) &=
\e(\p_x\{u_s^2q_x\},q_{xx}w^2)
-
(\e\p_x\{u_s^2q_x\},q_xw^2)_{x=L}
+
(\e\p_x\{u_s^2q_x\},q_xw^2)_{x=0}
\\
 &=
\e(u_s^2q_{xx},q_{xx}w^2)
+
2\e(u_su_{sx}q_x,q_{xx}w^2)
\\
&\quad -
(\e u_s^2q_{xx},q_xw^2)_{x=L}
-
(2\e u_su_{sx}q_x,q_xw^2)_{x=L}
\\
&\quad +
(\e u_s^2q_{xx},q_xw^2)_{x=0}
+
(2\e u_su_{sx}q_x,q_xw^2)_{x=0}.
\end{aligned}
\end{align}

\((\ref{RqxX}.1)\) is favorable: $(\ref{RqxX}.1) = \|\sqrt\e u_sq_{xx}w\|^2$. For \((\ref{RqxX}.2)\), Poincar\'e's inequality gives

\begin{align*}
\begin{aligned}
|(\ref{RqxX}.2)| &\lesssim
\|\sqrt\e q_xw\|
\|\sqrt\e u_sq_{xx}w\|
\\
 &\lesssim
L\bigl(\cD_2^2[q]+\cD_3^2[q]\bigr).
\end{aligned}
\end{align*}

At $x=0$, $G=u_sq$ and the conditions $G=G_{xx}=0$ in
\eqref{eq:main-equation} give $(\ref{RqxX}.5)+(\ref{RqxX}.6)=0$. At $x=L$, retain the $u_s$ weights and use
$u_sq_x=-u_{sx}q$:

\begin{align*}
 |(\ref{RqxX}.3)+(\ref{RqxX}.4)| &\lesssim\|\sqrt\e u_sq_{xx}w\|_{x=L}
              \|\sqrt\e u_sq_xw\|_{x=L}
          +\|\sqrt\e u_sq_xw\|_{x=L}^2,\\
 \|\sqrt\e u_sq_xw\|_{x=L} &\lesssim\|\sqrt\e qw\|_{x=L}
 \leq\sqrt L\|\sqrt\e q_xw\|
 \lesssim L^{\frac{3}{2}}\bigl(\cD_2[q]+\cD_3[q]\bigr).
\end{align*}

These boundary terms are bounded by
$CL(\cD_2^2[q]+\cD_3^2[q])$.

For \((\ref{Rqx}.2)\), use \(w^2=r^3\) and integrate by parts in \(Y\):

\begin{align}\label{RqxY}
\begin{aligned}
(\ref{Rqx}.2) &=
\e(\p_x\{u_s^2q_Y\},q_{xY}w^2)
-
(\e\p_x\{u_s^2q_Y\},q_xw^2)_{Y=0}^{Y=1}
\\
 &=
\e(u_s^2q_{xY},q_{xY}w^2)
+
2\e(u_su_{sx}q_Y,q_{xY}w^2).
\end{aligned}
\end{align}

The boundary contribution vanishes at \(Y=0,1\). The first term is favorable: $(\ref{RqxY}.1) = \|\sqrt\e u_sq_{xY}w\|^2$. For \((\ref{RqxY}.2)\),

\begin{align*}
\begin{aligned}
|(\ref{RqxY}.2)| &\lesssim
\|\sqrt\e q_Yw\|
\|\sqrt\e q_{xY}w\|
\\
 &\lesssim
L\bigl(\cD_2^2[q]+\cD_3^2[q]\bigr).
\end{aligned}
\end{align*}

\subsubsection*{\normalfont \textit{Step 2: \(\e\cA^2\) Terms}}

\begin{align*}
&\left|
(\e\cA^2G,\e q_xw^2)
\right|
\lesssim
\sqrt L\left(
        \bigl(\cD_2^2[q]+\cD_3^2[q]\bigr)+\cD_4^2[G]
\right).
\end{align*}

Expanding the viscous term gives

\begin{align}\label{Aqx}
\begin{aligned}
(\e\cA^2G,\e q_xw^2) &=
(\e^2G_{xxxx},q_xw^2) +
(2\e^2\cL G_{xx},q_xw^2) +
(\e^2\cL^2G,q_xw^2).
\end{aligned}
\end{align}

For \((\ref{Aqx}.1)\), integrate by parts in \(x\):
\begin{align}\label{Gxxxxqx}
\begin{aligned}
&(\ref{Aqx}.1)
=
-\e^2(G_{xxx},q_{xx}w^2)
+
(\e^2G_{xxx},q_xw^2)_{x=L}
-
(\e^2G_{xxx},q_xw^2)_{x=0}.
\end{aligned}
\end{align}
The \(x=L\) boundary term vanishes. Expanding $G_{xxx} = u_sq_{xxx} + 3u_{sx}q_{xx} + 3u_{sxx}q_x + u_{sxxx}q$ gives

\begin{align*}
\begin{aligned}
|(\ref{Gxxxxqx}.1)| &\lesssim
\sqrt\e\|\e\sqrt{u_s}q_{xxx}w\|
\|\sqrt\e q_{xx}w\|
+\e
\|\sqrt\e q_{xx}w\|^2
+
\|\e q_xw\|\|\e q_{xx}w\|
+\e
\|\sqrt\e qw\|\|\sqrt\e q_{xx}w\|
\\
 &\lesssim
\sqrt\e \,\bigl(\cD_2^2[q]+\cD_3^2[q]\bigr).
\end{aligned}
\end{align*}

The \(x=0\) boundary term satisfies

\begin{align*}
\begin{aligned}
|(\ref{Gxxxxqx}.3)| &\lesssim\|\e^{\frac{3}{2}}G_{xxx}w\|_{x=0}
                   \|\sqrt\e q_xw\|_{x=0}\\
  &\lesssim\sqrt L\bigl(\cD_2[q]+\cD_3[q]\bigr)^{\frac{3}{2}}
                  \cD_4^{\frac{1}{2}}[G]
 \lesssim\sqrt L\bigl(\cD_2^2[q]+\cD_3^2[q]+\cD_4^2[G]\bigr).
\end{aligned}
\end{align*}

Hence

\begin{align*}
&|(\ref{Aqx}.1)|
\lesssim
\sqrt L\left(
        \bigl(\cD_2^2[q]+\cD_3^2[q]\bigr)+\cD_4^2[G]
\right).
\end{align*}

For \((\ref{Aqx}.2)\), the self-adjoint form of \(\cL\) gives
\begin{align}\label{LGxxqx}
\begin{aligned}
&(\ref{Aqx}.2)
=
-2\e^2(G_{xxY},q_{xY}w^2)
+
(2\e^2G_{xxY},q_xw^2)_{Y=0}^{Y=1}.
\end{aligned}
\end{align}
The \(Y\)-boundary terms vanish. Expanding
\[
        G_{xxY}
        =
        u_sq_{xxY}
        +
        u_{sY}q_{xx}
        +
        2u_{sx}q_{xY}
        +
        2u_{sxY}q_x
        +
        u_{sxx}q_Y
        +
        u_{sxxY}q ,
\]
we obtain

\begin{align*}
\begin{aligned}
|(\ref{LGxxqx}.1)| &\lesssim
\sqrt\e\|\e\sqrt{u_s}q_{xxY}w\|
\|\sqrt\e q_{xY}w\|
+
\sqrt\e\|\sqrt\e q_{xx}w\|\|\sqrt\e q_{xY}w\|
\\
&\quad +
\e\|\sqrt\e q_{xY}w\|^2
+
\sqrt\e\|\sqrt\e q_xw\|\|\sqrt\e q_{xY}w\|
+
\e\|\sqrt\e q_Yw\|\|\sqrt\e q_{xY}w\|
\\
&\quad +
\sqrt\e\|\sqrt\e qw\|\|\sqrt\e q_{xY}w\|
\\
 &\lesssim
\sqrt \e\,\bigl(\cD_2^2[q]+\cD_3^2[q]\bigr).
\end{aligned}
\end{align*}

For \((\ref{Aqx}.3)\), integrate by parts twice in \(Y\), then expand
\(\cL G\):

\begin{align}\label{L2Gqx}
\begin{aligned}
(\ref{Aqx}.3) &=
(\e^2\cL G,\cL q_xw^2)
-
(\e^2(\cL G)_Y,q_xw^2)_{Y=0}
+
(\e^2\cL G,q_{xY}w^2)_{Y=0}
\\
  &=
\e^2(u_s\cL q,\cL q_xw^2)
+
2\e^2(u_{sY}q_Y,\cL q_xw^2)
+
\e^2(q\cL u_s,\cL q_xw^2) +
(\e^2\cL G,q_{xY}w^2)_{Y=0}.
\end{aligned}
\end{align}

The boundary term \(-(\e^2(\cL G)_Y,q_xw^2)_{Y=0}\) vanishes since \(q_x|_{Y=0}=0\).
For \((\ref{L2Gqx}.1)\),

\begin{align*}
\begin{aligned}
|(\ref{L2Gqx}.1)| &\lesssim
\|\e\cL q\,w\|
\|\e\sqrt{u_s}\cL q_xw\|
\\
 &\lesssim
\sqrt L\,\bigl(\cD_2^2[q]+\cD_3^2[q]\bigr).
\end{aligned}
\end{align*}

For \((\ref{L2Gqx}.2)\), integrate by parts in \(Y\):

\begin{align*}
\begin{aligned}
(\ref{L2Gqx}.2) &=
-2\e^2(u_{sY}\cL q,q_{xY}w^2)
-2\e^2(\bar\cL u_s q_Y,q_{xY}w^2) -(2\e^2u_{sY}q_Y,q_{xY}w^2)_{Y=0}.
\end{aligned}
\end{align*}

The boundary term in \((\ref{L2Gqx}.2)\) cancels with \((\ref{L2Gqx}.4)\), since
\[
(\e^2\cL G,q_{xY}w^2)_{Y=0}
=
(2\e^2u_{sY}q_Y,q_{xY}w^2)_{Y=0}.
\]
Thus

\begin{align*}
\begin{aligned}
|(\ref{L2Gqx}.2)+(\ref{L2Gqx}.4)| &\lesssim
\|\e\cL q\,w\|\|\sqrt\e q_{xY}w\|
+
L\|\sqrt\e q_{xY}w\|^2
\\
 &\lesssim
\sqrt L\,\bigl(\cD_2^2[q]+\cD_3^2[q]\bigr).
\end{aligned}
\end{align*}

For \((\ref{L2Gqx}.3)\),

\begin{align*}
\begin{aligned}
|(\ref{L2Gqx}.3)| &\lesssim \|\sqrt\e Y\frac{\cL u_s}{\sqrt{u_s}}\|_{\infty} \|\sqrt\e\frac{q}{Y}w\| \|\e\sqrt{u_s}\cL q_x w\|\\
 &\lesssim L\,\bigl(\cD_2^2[q]+\cD_3^2[q]\bigr).
\end{aligned}
\end{align*}

\subsubsection*{\normalfont \textit{Step 3: $J$ terms}}

\begin{align}\label{Jqx}
&(-u_{sx} G_{xx}-u_{sx} \cL G
+G(\cL u_{sx})+G u_{sxxx}
+w_s(\cL G)_Y+w_sG_{xxY} \n \\*
&\quad +2G \beta_{sx}-w_{sxx}G_Y+u_{sxY}G_Y
-\frac2r u_{sY}G_x+2\beta_sG_x,-\e q_{x}w^2) .
\end{align}

Integration by parts in $x$ and $Y$ gives

\begin{align*}
    (\ref{Jqx}.2)+(\ref{Jqx}.9) &= -(u_{sx}G_Y,\e q_{xY}w^2)
    -2(u_{sxY}G_Y,\e q_{x}w^2)\\
    |(\ref{Jqx}.2)+(\ref{Jqx}.9)| &\lesssim L\|\sqrt\e G_{xY}w\|\|\sqrt\e q_{xY}w\|+L\|Yu_{sxY}\|_{\infty}\|\sqrt\e G_{xY} w\|\|\sqrt\e \frac{q_x}{Y} w\|\\
 &\lesssim L\bigl(\cD_2^2[q]+\cD_3^2[q]\bigr).   
\end{align*}

\begin{align*}
    (\ref{Jqx}.3) &= (u_{sxY}G_Y,\e q_{x}w^2)+(u_{sxY}G,\e q_{xY}w^2)
    \\
    |(\ref{Jqx}.3)| &\lesssim L\|\sqrt\e G_{xY}w\|\|\sqrt\e q_{xY}w\| \lesssim L\bigl(\cD_2^2[q]+\cD_3^2[q]\bigr)
\end{align*}

\begin{align*}
    (\ref{Jqx}.5) &=\e((w_{sY}-3\frac{w_s}{r})\cL G,q_xw^2)
   +\e(w_s\cL G,q_{xY}w^2)\\
 |(\ref{Jqx}.5)| &\lesssim\left\|\frac{Y(w_{sY}-3\frac{w_s}{r})}{\sqrt\e}\right\|_\infty
       \|\e\cL Gw\|\|\sqrt\e \frac{q_x}{Y}w\|
 +\left\|\frac{w_s}{\sqrt\e}\right\|_\infty
       \|\e\cL Gw\|\|\sqrt\e q_{xY}w\|\\
  &\lesssim L\bigl(\cD_2^2[q]+\cD_3^2[q]\bigr).
\end{align*}

Theorem~\ref{thm:approximate-profiles} bounds the remaining coefficients.

\begin{align*}
|(\ref{Jqx}.7,10,11)| &\lesssim L\|Y\{\beta_{sx},\beta_s,\frac{u_{sY}}{r}\}\|_{\infty}\|\sqrt\e\frac{\{G,G_x\}}{Y}w\|\|\sqrt\e q_{xx}w\|\\
 &\lesssim L\bigl(\cD_2^2[q]+\cD_3^2[q]\bigr). 
\end{align*}

The remaining terms satisfy

\begin{align*}
 |(\ref{Jqx}.1)|+|(\ref{Jqx}.4)| &\lesssim\|\sqrt\e\{G,G_{xx}\}w\|\|\sqrt\e q_xw\|,\\
 |(\ref{Jqx}.6)|+|(\ref{Jqx}.8)| &\lesssim L\|\e G_{xxY}w\|\|\sqrt\e q_{xx}w\|
       +\sqrt\e\|\sqrt\e G_Yw\|\|\sqrt\e q_xw\|\\
  &\lesssim L\bigl(\cD_2^2[q]+\cD_3^2[q]\bigr),
\end{align*}

and \eqref{small.guys} gives the same bound for the first line.

\end{proof}
\subsection{Side-trace estimates}

Set $\bar B_G(w) := \|\e^{\frac32}u_s\cL G_xw\|_{x=0}^2 + \|\e^{\frac32}u_sG_{xxY}w\|_{x=L}^2$.
\begin{lemma}
Let $G$ solve \eqref{eq:main-equation}, and let
$w=r^{\frac32}$. Then for any fixed $\eta>0$,

\begin{align}
&\bar B_G(w)
+
\big\|
\{\e^2\sqrt{u_s}(\cL G_x)_Y,\e^2\sqrt{u_s}\cL G_{xx},\e^2\sqrt{u_s}G_{xxxY}\}w
\big\|^2
\n\\ \label{trace 1}
 &\lesssim
\eta\|\e^2\cL G_{xx}w\|^2
+
C_\eta\e^{\frac13}\bigl(\cD_2^2[q]+\cD_3^2[q]\bigr)
+
\left|
(F_G,-\e^3u_s\cL G_{xx}w^2)
\right|.
\end{align}

\end{lemma}

\begin{proof}
Using $R[q]=u_s\cA G-G\cA u_s$, we find that $\partial_xR[q]$
and the first three terms of $\cJ[G]$ sum to $u_s(\cA G)_x$.

\subsubsection*{\normalfont \textit{Step 1: Rayleigh terms and the canceling part of $J$}}

\begin{align}\label{RtraceMix}
\begin{aligned}
(u_s(\cA G)_x,-\e^3u_s\cL G_{xx}w^2) &=
-(u_sG_{xxx},\e^3u_s\cL G_{xx}w^2)
-(u_s(\cL G)_x,\e^3u_s\cL G_{xx}w^2).
\end{aligned}
\end{align}

For $(\ref{RtraceMix}.1)$, integrate by parts first in $Y$ and then in $x$ to obtain

\begin{align}\label{RtraceMixI}
\begin{aligned}
(\ref{RtraceMix}.1) &=
\e^3(u_s^2G_{xxxY},G_{xxY}w^2)
+
\e^3(2u_su_{sY}G_{xxx},G_{xxY}w^2)
\\
 &=
\frac12\|\e^{\frac32}u_sG_{xxY}w\|_{x=L}^2
-
\frac12\|\e^{\frac32}u_sG_{xxY}w\|_{x=0}^2
-
\e^3(u_su_{sx}G_{xxY},G_{xxY}w^2)
\\
&\quad +
\e^3(2u_su_{sY}G_{xxx},G_{xxY}w^2).
\end{aligned}
\end{align}

The $x=0$ boundary term vanishes since $G_{xx}|_{x=0}=0$, and

\begin{align*}
\begin{aligned}
|(\ref{RtraceMixI}.3)| &\lesssim\e^3\|G_{xxY}w\|^2,
\\
|(\ref{RtraceMixI}.4)| &\lesssim
\sqrt\e\,\|\e G_{xxx}w\|\,\|\e G_{xxY}w\|
\lesssim
\sqrt\e\,\bigl(\cD_2^2[q]+\cD_3^2[q]\bigr).
\end{aligned}
\end{align*}

For $(\ref{RtraceMix}.2)$, integrate by parts in $x$:

\begin{align}\label{RtraceMixII}
\begin{aligned}
(\ref{RtraceMix}.2) &=
-\e^3(u_s(\cL G)_x,u_s\cL G_{xx}w^2)
\\
 &=
\frac12\|\e^{\frac32}u_s\cL G_xw\|_{x=0}^2
-
\frac12\|\e^{\frac32}u_s\cL G_xw\|_{x=L}^2
+
\e^3(u_su_{sx}\cL G_x,\cL G_xw^2).
\end{aligned}
\end{align}

Since $G_x|_{x=L}=0$, $\cL G_x|_{x=L}=0$. The coefficient term satisfies

\begin{align*}
 &|(\ref{RtraceMixII}.3)|\lesssim\e^3\|\cL G_xw\|^2.
\end{align*}

\subsubsection*{\normalfont \textit{Step 2: $\e\cA^2$ Terms}}

Use the finite-mode approximation from the proof of
Lemma~\ref{lem:main-principal-solvability} with $F_G=\e\cA^2G$.
The integration by parts used there to prove \eqref{sum.Gxxyy.lap} gives

\begin{align*}
\begin{aligned}
(\e\cA^2G,\e^3u_s\cL G_{xx}w^2) &\ge
c\big\|
\{\e^2\sqrt{u_s}(\cL G_x)_Y,\e^2\sqrt{u_s}\cL G_{xx},\e^2\sqrt{u_s}G_{xxxY}\}w
\big\|^2 -
C\e\bigl(\cD_2^2[q]+\cD_3^2[q]\bigr).
\end{aligned}
\end{align*}

\subsubsection*{\normalfont \textit{Step 3: The remaining $J$ terms}}

Write

\begin{align}\label{JtraceMix}
\begin{aligned}
&(\bigl(w_s(\cA G)_Y+(2G-rG_Y)(\beta_s)_x+rG_x(\beta_s)_Y\bigr),-\e^3u_s\cL G_{xx}w^2)\\
 &=
(w_s(\cA G)_Y,-\e^3u_s\cL G_{xx}w^2)
\\
&\quad +((2G-rG_Y)(\beta_s)_x+rG_x(\beta_s)_Y,
        -\e^3u_s\cL G_{xx}w^2).
\end{aligned}
\end{align}

For $(\ref{JtraceMix}.1)$, the $G_{xxY}$-part satisfies

\begin{align*}
\begin{aligned}
\left|
-\e^3(w_sG_{xxY},u_s\cL G_{xx}w^2)
\right| &\lesssim
\sqrt\e\,\|\e G_{xxY}w\|\,
\|\e^2\sqrt{u_s}\cL G_{xx}w\|
\\
 &\le
\eta\|\e^2\sqrt{u_s}\cL G_{xx}w\|^2
+
C_\eta\e\bigl(\cD_2^2[q]+\cD_3^2[q]\bigr).
\end{aligned}
\end{align*}

For the \((\cL G)_Y\)-part, integrate by parts in \(x\):

\begin{align}\label{JtraceMixI}
\begin{aligned}
-\e^3(w_s(\cL G)_Y,u_s\cL G_{xx}w^2) &=
\e^3(\p_x\{w_su_s(\cL G)_Y\},\cL G_xw^2) -
(\e^3w_su_s(\cL G)_Y,\cL G_xw^2)_{x=0}^{x=L}.
\end{aligned}
\end{align}

The boundary term vanishes by $\cL G_x|_{x=L}=0$, $(\cL G)_Y|_{x=0}=0$. Expanding the \(x\)-derivative,
\begin{align*}
\begin{aligned}
&(\ref{JtraceMixI}.1)
=
\e^3(w_su_s(\cL G_x)_Y,\cL G_xw^2)
+
\e^3(\p_x\{w_su_s\}(\cL G)_Y,\cL G_xw^2).
\end{aligned}
\end{align*}
Hence

\begin{align*}
\begin{aligned}
&|(\ref{JtraceMix}.1)|
\le
\eta
\big\|
\{\e^2\sqrt{u_s}(\cL G_x)_Y,\e^2\sqrt{u_s}\cL G_{xx}\}w
\big\|^2
+
C_\eta\sqrt\e\,\bigl(\cD_2^2[q]+\cD_3^2[q]\bigr).
\end{aligned}
\end{align*}

In $(\ref{JtraceMix}.2)$, substitute $r(\beta_s)_x=w_{sxx}-u_{sxY}$ and $r(\beta_s)_Y = -u_{sxx}-\cL u_s-\frac2r u_{sY}+2\beta_s$:

\begin{align}\label{JtraceMixBeta}
\begin{aligned}
(\ref{JtraceMix}.2) &=
-(2\e^3G\beta_{sx},u_s\cL G_{xx}w^2)
+
(\e^3w_{sxx}G_Y,u_s\cL G_{xx}w^2)
\\
&\quad -
(\e^3u_{sxY}G_Y,u_s\cL G_{xx}w^2)
+
(\e^3u_{sxx}G_x,u_s\cL G_{xx}w^2)
\\
&\quad +
(\e^3\cL u_s\,G_x,u_s\cL G_{xx}w^2)
+
\left(\frac{2\e^3}{r}u_{sY}G_x,u_s\cL G_{xx}w^2\right)
\\
&\quad -
(2\e^3\beta_sG_x,u_s\cL G_{xx}w^2).
\end{aligned}
\end{align}

For \((\ref{JtraceMixBeta}.1)\), multiply and divide by \(Y\):

\begin{align*}
\begin{aligned}
|(\ref{JtraceMixBeta}.1)| &\lesssim
\sqrt\e
\left\|\sqrt\e\frac{G}{Y}w\right\|
\|\e^2\sqrt{u_s}\cL G_{xx}w\|
\\
 &\le
\eta\|\e^2\sqrt{u_s}\cL G_{xx}w\|^2
+
C_\eta\sqrt\e\,\bigl(\cD_2^2[q]+\cD_3^2[q]\bigr).
\end{aligned}
\end{align*}

The \(w_{sxx}G_Y\)-term satisfies

\begin{align*}
\begin{aligned}
|(\ref{JtraceMixBeta}.2)| &\lesssim
\e\,\|\sqrt\e G_Yw\|\,
\|\e^2\sqrt{u_s}\cL G_{xx}w\|
\\
 &\le
\eta\|\e^2\sqrt{u_s}\cL G_{xx}w\|^2
+
C_\eta\e^2\bigl(\cD_2^2[q]+\cD_3^2[q]\bigr).
\end{aligned}
\end{align*}

In the \(u_{sxY}G_Y\)-term, integrate by parts in \(x\):

\begin{align*}
\begin{aligned}
(\ref{JtraceMixBeta}.3) &=
-\e^3(u_{sxY}G_Y,u_s\cL G_{xx}w^2)
\\
 &=
\e^3(\p_x\{u_su_{sxY}G_Y\},\cL G_xw^2)
-
(\e^3u_su_{sxY}G_Y,\cL G_xw^2)_{x=0}^{x=L}.
\end{aligned}
\end{align*}

The boundary term vanishes at \(x=L\) since \(\cL G_x|_{x=L}=0\), and at
\(x=0\) since \(G_Y|_{x=0}=0\).

\begin{align*}
\begin{aligned}
|(\ref{JtraceMixBeta}.3)| &\lesssim
\e\left(
        \|\sqrt\e G_{xY}w\|
        +
        \|\sqrt\e G_Yw\|
\right)
\|\e\cL G_xw\|
\\
 &\lesssim
\e\bigl(\cD_2^2[q]+\cD_3^2[q]\bigr).
\end{aligned}
\end{align*}

The \(u_{sxx}G_x\)-term satisfies

\begin{align*}
\begin{aligned}
|(\ref{JtraceMixBeta}.4)| &\lesssim
\sqrt\e\,
\left\|\sqrt\e G_xw\right\|
\|\e^2\sqrt{u_s}\cL G_{xx}w\|
\\
 &\le
\eta\|\e^2\sqrt{u_s}\cL G_{xx}w\|^2
+
C_\eta\e\bigl(\cD_2^2[q]+\cD_3^2[q]\bigr).
\end{aligned}
\end{align*}

In the \(\cL u_sG_x\)-term, integrate by parts in \(x\) and use
Young's inequality with exponents \(\frac43\) and \(4\):

\begin{align*}
\begin{aligned}
(\ref{JtraceMixBeta}.5) &=
-\e^3(u_s\cL u_sG_{xx},\cL G_xw^2)
-\e^3((u_s\cL u_s)_xG_x,\cL G_xw^2) -(\e^3u_s\cL u_sG_x,\cL G_xw^2)_{x=0}
\\
|(\ref{JtraceMixBeta}.5)| &\lesssim
\sqrt\e\,\|\e\cL u_s\|_\infty
\|\sqrt\e G_{xx}w\|
\|\e\cL G_xw\|
\\
&\quad +
\sqrt\e
\left\|\sqrt\e\frac{G_x}{Y}w\right\|
\|\e\cL G_xw\|
\\
&\quad +
\e^{\frac14}\|\sqrt\e Y\cL u_s\|_\infty
\left\|\sqrt\e\frac{G_x}{Y}w\right\|^{\frac12}
\left\|\e\frac{G_{xx}}{Y}w\right\|^{\frac12}
\|\e\cL G_xw\|^{\frac12}
\|\e^2\cL G_{xx}w\|^{\frac12}
\\
 &\lesssim
\sqrt\e\,\bigl(\cD_2^2[q]+\cD_{3,x}^2[q]\bigr)
+
\widetilde C_\eta\e^{\frac13}\bigl(\cD_2^2[q]+\cD_{3,x}^2[q]\bigr)
+
\eta\|\e^2\cL G_{xx}w\|^2.
\end{aligned}
\end{align*}

For the \(\frac{u_{sY}G_x}{r}\) and \(\beta_sG_x\) terms, multiply and divide by \(Y\):

\begin{align*}
\begin{aligned}
|(\ref{JtraceMixBeta}.6)+(\ref{JtraceMixBeta}.7)| &\lesssim
\sqrt\e\,
\left\|\sqrt\e\frac{G_x}{Y}w\right\|
\|\e^2\sqrt{u_s}\cL G_{xx}w\|
\\
 &\le
\eta\|\e^2\sqrt{u_s}\cL G_{xx}w\|^2
+
C_\eta\sqrt\e\,\bigl(\cD_2^2[q]+\cD_3^2[q]\bigr).
\end{aligned}
\end{align*}

\end{proof}

\begin{lemma}
Let $G$ solve \eqref{eq:main-equation}, and let
$w=r^{\frac32}$. Then

\begin{align}
&\|\e^{\frac32}\sqrt{u_s}G_{xxx}w\|_{x=0}^2
+
\big\|\{\e^2G_{xxxx},\e^2G_{xxxY},\e^2\cL G_{xx}\}w\big\|^2
\n\\ \label{trace 2}
 &\lesssim
\cD_2^2[q]+\cD_3^2[q]
+
\left|
(F_G,-\e^3G_{xxxx}w^2)
\right|.
\end{align}

\end{lemma}

\begin{proof}

\subsubsection*{\normalfont \textit{Step 1: Rayleigh terms and the canceling part of $J$}}

\begin{align*}
\begin{aligned}
&(u_s(\cA G)_x,-\e^3G_{xxxx}w^2)
\\
 &\ge
c\|\e^{\frac32}\sqrt{u_s}G_{xxx}w\|_{x=0}^2
-
\eta\big\|\{\e^2G_{xxxx},\e^2G_{xxxY},\e^2\cL G_{xx}\}w\big\|^2
-
C_\eta\bigl(\cD_2^2[q]+\cD_3^2[q]\bigr).
\end{aligned}
\end{align*}

\begin{align}\label{RtraceX}
\begin{aligned}
(u_s(\cA G)_x,-\e^3G_{xxxx}w^2) &=
-(u_sG_{xxx},\e^3G_{xxxx}w^2)
-(u_s\cL G_x,\e^3G_{xxxx}w^2).
\end{aligned}
\end{align}

For $(\ref{RtraceX}.1)$, integrate by parts in $x$:

\begin{align}\label{RtraceXI}
\begin{aligned}
(\ref{RtraceX}.1) &=
-\e^3(u_sG_{xxx},G_{xxxx}w^2)
\\
 &=
\frac12\|\e^{\frac32}\sqrt{u_s}G_{xxx}w\|_{x=0}^2
-
\frac12\|\e^{\frac32}\sqrt{u_s}G_{xxx}w\|_{x=L}^2
+
\frac12\e^3(u_{sx}G_{xxx},G_{xxx}w^2).
\end{aligned}
\end{align}

The $x=L$ boundary term vanishes by $G_{xxx}|_{x=L}=0$, and
$|(\ref{RtraceXI}.3)|\lesssim\e^3\|G_{xxx}w\|^2$.
For $(\ref{RtraceX}.2)$, integrate by parts in $x$:

\begin{align}\label{RtraceXII}
\begin{aligned}
(\ref{RtraceX}.2) &=
-\e^3(u_s\cL G_x,G_{xxxx}w^2)
\\
 &=
\e^3(u_s\cL G_{xx}+u_{sx}\cL G_x,G_{xxx}w^2)
-
(\e^3u_s\cL G_x,G_{xxx}w^2)^{x=L}_{x=0}.
\end{aligned}
\end{align}

The $x=L$ boundary term vanishes. For the $x=0$ boundary term,

\begin{align*}
\begin{aligned}
&\left|(\e^3u_s\cL G_x,G_{xxx}w^2)_{x=0}\right|
\le
\eta\|\e^{\frac32}\sqrt{u_s}G_{xxx}w\|_{x=0}^2
+
C_\eta\|\e^{\frac32}\sqrt{u_s}\cL G_xw\|_{x=0}^2.
\end{aligned}
\end{align*}

For the first interior term in $(\ref{RtraceXII}.1)$, integrate by parts in $Y$ and then in $x$:

\begin{align}\label{RtraceXIIa}
\begin{aligned}
\e^3(u_s\cL G_{xx},G_{xxx}w^2) &=
-\e^3(u_sG_{xxY},G_{xxxY}w^2)
-
\e^3(u_{sY}G_{xxY},G_{xxx}w^2)
+
(\e^3u_sG_{xxY},G_{xxx}w^2)_{Y=0}^{Y=1}
\\
 &=
-
\frac12\|\e^{\frac32}\sqrt{u_s}G_{xxY}w\|_{x=L}^2
+
\frac12\|\e^{\frac32}\sqrt{u_s}G_{xxY}w\|_{x=0}^2
+
\frac12\e^3(u_{sx}G_{xxY},G_{xxY}w^2)
\\
&\quad -
\e^3(u_{sY}G_{xxY},G_{xxx}w^2).
\end{aligned}
\end{align}

The $x=0$ boundary term vanishes since $G_{xx}|_{x=0}=0$, and

\begin{align*}
\begin{aligned}
|(\ref{RtraceXIIa}.3)| &\lesssim\e^3\|G_{xxY}w\|^2,
\\
|(\ref{RtraceXIIa}.4)| &\lesssim
\sqrt\e\,\|\e G_{xxY}w\|\,\|\e G_{xxx}w\|
\lesssim
\sqrt\e\,\bigl(\cD_2^2[q]+\cD_3^2[q]\bigr).
\end{aligned}
\end{align*}

The second interior term in $(\ref{RtraceXII}.1)$ satisfies

\begin{align*}
&\left|\e^3(u_{sx}\cL G_x,G_{xxx}w^2)\right|
 \lesssim\e^3\|\cL G_xw\|\|G_{xxx}w\|.
\end{align*}

For the two remaining boundary terms, use
$\cL G_x|_{x=L}=0$ and $G_{xxY}|_{x=0}=0$:

\begin{align*}
 \|\e^{\frac{3}{2}}\sqrt{u_s}\cL G_xw\|_{x=0}^2 &\lesssim\|\e\cL G_xw\|\|\e^2\cL G_{xx}w\|,\\
 \|\e^{\frac{3}{2}}\sqrt{u_s}G_{xxY}w\|_{x=L}^2 &\lesssim\|\e G_{xxY}w\|\|\e^2G_{xxxY}w\|.
\end{align*}

By Young's inequality, their sum is at most $\eta\|\{\e^2\cL G_{xx},\e^2G_{xxxY}\}w\|^2 +C_\eta\bigl(\cD_2^2[q]+\cD_3^2[q]\bigr)$.

\subsubsection*{\normalfont \textit{Step 2: $\e\cA^2$ Terms}}

Integration by parts as in \eqref{sum.Gxxxx.lap} gives

\begin{align*}
\begin{aligned}
&(\e\cA^2G,\e^3G_{xxxx}w^2)
 \ge
c\big\|\{\e^2\cL G_{xx},\e^2G_{xxxY},\e^2G_{xxxx}\}w\big\|^2
-
C\e\bigl(\cD_2^2[q]+\cD_{3,x}^2[q]\bigr).
\end{aligned}
\end{align*}

\subsubsection*{\normalfont \textit{Step 3: The remaining $J$ terms}}

\begin{align*}
\begin{aligned}
&\left|(\bigl(w_s(\cA G)_Y+(2G-rG_Y)(\beta_s)_x+rG_x(\beta_s)_Y\bigr),-\e^3G_{xxxx}w^2)\right|\\
 &\le
\eta
\big\|\{\e^2G_{xxxx},\e^2G_{xxxY},\e^2\cL G_{xx}\}w\big\|^2
\\
&\quad +
\eta\|\e^{\frac32}\sqrt{u_s}G_{xxx}w\|_{x=0}^2
+
C_\eta\e^{\frac13}\bigl(\cD_2^2[q]+\cD_3^2[q]\bigr).
\end{aligned}
\end{align*}

Write

\begin{align}\label{JtraceX}
\begin{aligned}
&(\bigl(w_s(\cA G)_Y+(2G-rG_Y)(\beta_s)_x+rG_x(\beta_s)_Y\bigr),-\e^3G_{xxxx}w^2)\\
 &=
(w_s(\cA G)_Y,-\e^3G_{xxxx}w^2)
\\
&\quad +
((2G-rG_Y)(\beta_s)_x+rG_x(\beta_s)_Y,
        -\e^3G_{xxxx}w^2).
\end{aligned}
\end{align}

For $(\ref{JtraceX}.1)$, the $G_{xxY}$-part satisfies

\begin{align*}
\begin{aligned}
\left|
-\e^3(w_sG_{xxY},G_{xxxx}w^2)
\right| &\lesssim
\sqrt\e\,\|\e G_{xxY}w\|\,
\|\e^2G_{xxxx}w\|
\\
 &\le
\eta\|\e^2G_{xxxx}w\|^2
+
C_\eta\e\bigl(\cD_2^2[q]+\cD_3^2[q]\bigr).
\end{aligned}
\end{align*}

For the \((\cL G)_Y\)-part, integrate by parts first in \(x\):
\begin{align}\label{JtraceXLG}
\begin{aligned}
&-\e^3(w_s(\cL G)_Y,G_{xxxx}w^2)
=
\e^3(w_{sx}(\cL G)_Y,G_{xxx}w^2)
+
\e^3(w_s(\cL G)_{xY},G_{xxx}w^2).
\end{aligned}
\end{align}
The boundary terms vanish by \(G_{xxx}|_{x=L}=0\) and \(G|_{x=0}=0\).
For the first term in \((\ref{JtraceXLG})\),

\begin{align*}
\begin{aligned}
&|(\ref{JtraceXLG}.1)|
\lesssim
\e\|\e^{\frac32}(\cL G)_Yw\|
\|\e G_{xxx}w\|
\lesssim
\e\bigl(\cD_2^2[q]+\cD_3^2[q]\bigr).
\end{aligned}
\end{align*}

In the second term of \((\ref{JtraceXLG})\), integrate by parts in $Y$.
The wall and axis terms vanish:

\begin{align*}
\begin{aligned}
(\ref{JtraceXLG}.2) &=
\e^3(w_s(\cL G)_{xY},G_{xxx}w^2)
\\
 &=
-\e^3(w_s\cL G_x,G_{xxxY}w^2)
-
\e^3\left(
\left(w_{sY}-\frac3r w_s\right)\cL G_x,
G_{xxx}w^2
\right)\\
|(\ref{JtraceXLG}.2)| &\lesssim
\sqrt\e\,\|\e\cL G_xw\|\|\e^2G_{xxxY}w\|
+
\e\,\|\e\cL G_xw\|\|\e G_{xxx}w\|
\\
 &\le
\eta\|\e^2G_{xxxY}w\|^2
+
C_\eta\sqrt\e\,\bigl(\cD_2^2[q]+\cD_3^2[q]\bigr).
\end{aligned}
\end{align*}

Hence

\begin{align*}
&|(\ref{JtraceX}.1)|
\le
\eta\big\|\{\e^2G_{xxxx},\e^2G_{xxxY}\}w\big\|^2
+
C_\eta\sqrt\e\,\bigl(\cD_2^2[q]+\cD_3^2[q]\bigr).
\end{align*}

In $(\ref{JtraceX}.2)$, substitute $r(\beta_s)_x=w_{sxx}-u_{sxY}$ and $r(\beta_s)_Y = -u_{sxx}-\cL u_s-\frac2r u_{sY}+2\beta_s$:

\begin{align}\label{JtraceXBeta}
\begin{aligned}
(\ref{JtraceX}.2) &=
-(2\e^3G\beta_{sx},G_{xxxx}w^2)
+
(\e^3w_{sxx}G_Y,G_{xxxx}w^2)
\\
&\quad -
(\e^3u_{sxY}G_Y,G_{xxxx}w^2)
+
(\e^3u_{sxx}G_x,G_{xxxx}w^2)
\\
&\quad +
(\e^3\cL u_s\,G_x,G_{xxxx}w^2)
+
\left(\frac{2\e^3}{r}u_{sY}G_x,G_{xxxx}w^2\right)
\\
&\quad -
(2\e^3\beta_sG_x,G_{xxxx}w^2).
\end{aligned}
\end{align}

For \((\ref{JtraceXBeta}.1)\), multiply and divide by \(Y\):

\begin{align*}
\begin{aligned}
|(\ref{JtraceXBeta}.1)| &\lesssim
\sqrt\e
\left\|\sqrt\e\frac{G}{Y}w\right\|
\|\e^2G_{xxxx}w\|
\\
 &\le
\eta\|\e^2G_{xxxx}w\|^2
+
C_\eta\sqrt\e\,\bigl(\cD_2^2[q]+\cD_3^2[q]\bigr).
\end{aligned}
\end{align*}

The \(w_{sxx}G_Y\)-term satisfies

\begin{align*}
\begin{aligned}
|(\ref{JtraceXBeta}.2)| &\lesssim
\e\,\|\sqrt\e G_Yw\|\,
\|\e^2G_{xxxx}w\|
\\
 &\le
\eta\|\e^2G_{xxxx}w\|^2
+
C_\eta\e^2\bigl(\cD_2^2[q]+\cD_3^2[q]\bigr).
\end{aligned}
\end{align*}

In the \(u_{sxY}G_Y\)-term, integrate by parts in \(x\):

\begin{align*}
\begin{aligned}
(\ref{JtraceXBeta}.3) &=
-\e^3(u_{sxY}G_Y,G_{xxxx}w^2)
\\
 &=
\e^3(\p_x\{u_{sxY}G_Y\},G_{xxx}w^2)
-
(\e^3u_{sxY}G_Y,G_{xxx}w^2)_{x=0}^{x=L}.
\end{aligned}
\end{align*}

The boundary term vanishes at \(x=L\) since \(G_{xxx}|_{x=L}=0\), and at
\(x=0\) since \(G_Y|_{x=0}=0\).

\begin{align*}
\begin{aligned}
|(\ref{JtraceXBeta}.3)| &\lesssim
\e
\left(
        \|\sqrt\e G_{xY}w\|
        +
        \|\sqrt\e G_Yw\|
\right)
\|\e G_{xxx}w\|
\\
 &\lesssim
\e\bigl(\cD_2^2[q]+\cD_3^2[q]\bigr).
\end{aligned}
\end{align*}

The \(u_{sxx}G_x\)-term satisfies

\begin{align*}
\begin{aligned}
|(\ref{JtraceXBeta}.4)| &\lesssim
\sqrt\e\,
\left\|\sqrt\e G_xw\right\|
\|\e^2G_{xxxx}w\|
\\
 &\le
\eta\|\e^2G_{xxxx}w\|^2
+
C_\eta\e\bigl(\cD_2^2[q]+\cD_3^2[q]\bigr).
\end{aligned}
\end{align*}

In the \(\cL u_sG_x\)-term, integrate by parts in \(x\) and use
Young's inequality with exponents \(\frac43\) and \(4\):

\begin{align*}
\begin{aligned}
(\ref{JtraceXBeta}.5) &=
-\e^3(\cL u_sG_{xx},G_{xxx}w^2)
-\e^3(\cL u_{sx}G_x,G_{xxx}w^2) -(\e^3\cL u_sG_x,G_{xxx}w^2)_{x=0}
\\
|(\ref{JtraceXBeta}.5)| &\lesssim
\sqrt\e\,\|\sqrt\e G_{xx}w\|
\|\e G_{xxx}w\|
+
\sqrt\e\left\|\sqrt\e\frac{G_x}{Y}w\right\|
\|\e G_{xxx}w\|
\\
&\quad +
\e^{\frac14}\|\sqrt\e Y\cL u_s\|_\infty
\left\|\sqrt\e\frac{G_x}{Y}w\right\|^{\frac12}
\left\|\e\frac{G_{xx}}{Y}w\right\|^{\frac12}
\|\e G_{xxx}w\|^{\frac12}
\|\e^2G_{xxxx}w\|^{\frac12}
\\
 &\le
\eta\|\e^2G_{xxxx}w\|^2
+
\eta\|\e^{\frac32}\sqrt{u_s}G_{xxx}w\|_{x=0}^2
+
C_\eta\e^{\frac13}\bigl(\cD_2^2[q]+\cD_{3,x}^2[q]\bigr).
\end{aligned}
\end{align*}

For the \(\frac{u_{sY}G_x}{r}\) and \(\beta_sG_x\) terms, multiply and divide by \(Y\):

\begin{align*}
\begin{aligned}
|(\ref{JtraceXBeta}.6)+(\ref{JtraceXBeta}.7)| &\lesssim
\sqrt\e\,
\left\|\sqrt\e\frac{G_x}{Y}w\right\|
\|\e^2G_{xxxx}w\|
\\
 &\le
\eta\|\e^2G_{xxxx}w\|^2
+
C_\eta\sqrt\e\,\bigl(\cD_2^2[q]+\cD_3^2[q]\bigr).
\end{aligned}
\end{align*}

\end{proof}

\begin{proposition}[Closure of the a priori estimate]
\label{prop:apriori-closure}
Let $0<\e\le L^{\frac{1}{8}}$ and let $G\in\mathcal H_0^4(w)$ solve
\eqref{eq:main-equation}. Then
\begin{equation}\label{eq:main-gain}
 \cD[q]\lesssim\sqrt\e L^{\frac1{16}}\|F_Gw\|
             \lesssim\|F_Gw\|.
\end{equation}
\end{proposition}

\begin{proof}
The estimates \eqref{trace 2}, \eqref{trace 1}, and Cauchy's inequality give
\[
 |G|_{\p,3}^2+
 \|\{\e^2G_{xxxx},\e^2G_{xxxY},\e^2\cL G_{xx},
       \e^2\sqrt{u_s}(\cL G_x)_Y\}w\|^2
 \lesssim\cD_2^2[q]+\cD_3^2[q]+\e^2\|F_Gw\|^2.
\]
Equation~\eqref{eq:main-equation} gives
\[
 \e^2\cL^2G
 =\e\{\partial_xR[q]+\cJ[G]\}-\e F_G
  -\e^2G_{xxxx}-2\e^2\cL G_{xx}.
\]
By $R[q]=u_s\cA G-G\cA u_s$, \eqref{systemat1}, and the profile bounds,
$\|\e\{\partial_xR[q]+\cJ[G]\}w\|\lesssim\cD_2[q]+\cD_3[q]$.
Thus
\begin{equation}\label{eq:main-fourth-recovery}
 \cD_4^2[G]\lesssim\cD_2^2[q]+\cD_3^2[q]+\e^2\|F_Gw\|^2.
\end{equation}

Fix $\beta=\frac{3}{8}$ in \eqref{qxx}. Multiply \eqref{qx-est} by a
sufficiently large constant times $L^{-\frac38}$ and add
\eqref{qxx} and \eqref{Lq-est}. After substituting
\eqref{eq:main-fourth-recovery}, absorb the error terms for
$L\ll1$ to obtain
\begin{align}
 &\cD_3^2[q]+L^{-\frac38}\cD_2^2[q]
 \le C\|F_Gw\|\bigl(\e\|q_{xx}w\|
                  +L^{-\frac38}\e\|q_xw\|
                  +\e^2\|\cL q\,w\|\bigr)
       +C\e^2L^{\frac18}\|F_Gw\|^2.
 \label{eq:main-forcing-before}
\end{align}
Use \eqref{interp.1} with $\eta=L^{\frac1{16}}$. Then

\begin{align*}
 \|\sqrt\e q_{xx}w\| &\le L^{\frac1{16}}\cD_{3,x}[q]
                         +CL^{-\frac18}\cD_2[q]\\
  &\le CL^{\frac1{16}}
       \bigl(\cD_3^2[q]+L^{-\frac38}\cD_2^2[q]\bigr)^{\frac12}.
\end{align*}

Also, \eqref{Poincare.2} and \eqref{small.guys} give
\[
 \|\sqrt\e q_xw\|\lesssim L\|\sqrt\e q_{xx}w\|,
 \qquad
 \|\e\cL q\,w\|\lesssim\sqrt L\cD_3[q].
\]
Insert these in \eqref{eq:main-forcing-before}. Since $L\le1$ and
$\e\le1$, Young's inequality gives
\[
 \cD_3^2[q]+L^{-\frac38}\cD_2^2[q]
 \lesssim(\e+\e^2)L^{\frac18}\|F_Gw\|^2.
\]
Hence $\cD_2[q]+\cD_3[q]\lesssim\sqrt\e L^{\frac1{16}}\|F_Gw\|$.
Equation~\eqref{eq:main-fourth-recovery} gives
$\cD_4[G]\lesssim(\sqrt\e L^{\frac1{16}}+\e)\|F_Gw\|$.
Since $0<\e\le L^{\frac{1}{8}}$,
$\e\le\sqrt\e L^{\frac1{16}}$, proving \eqref{eq:main-gain}.
This also verifies the a priori estimate assumed in
Proposition~\ref{prop:main-solvability}.
\end{proof}

\subsection{Scalar swirl estimates}
\label{subsec:swirl-energy}

\begin{lemma}[Scalar swirl estimate]\label{prop:sw-scalar-estimates}
Let $s^\e\in\mathscr V_s^2$ solve \eqref{eq:swirl-equation} with
$F_sw\in L^2(\Omega)$. Then
\begin{equation}\label{eq:sw-S-mixed-source}
 \mathcal S_\e[s^\e]^2\lesssim\|F_sw\|^2.
\end{equation}
\end{lemma}

\begin{proof}
Multiplying \eqref{eq:swirl-equation} by $s^\e w^2$ and using
\eqref{eq:usbound} and \eqref{eq:profile-quotient-bounds} gives

\begin{equation}\label{eq:sw-H-preestimate}
\begin{aligned}
 \e\|\nabla s^\e w\|^2
 +\tfrac12\|\sqrt{u_s}s^\e w\|_{x=L}^2 &=(F_s,s^\e w^2)
       +\left(\frac{w_s}{r}s^\e,s^\e w^2\right)\\
  &\lesssim\|F_sw\|\bigl(
       \eta\sqrt\e\|s_Y^\e w\|
       +\eta^{-2}L\|\sqrt{u_s}s_x^\e w\|\bigr)
       +\sqrt\e L^2\|\sqrt{u_s}s_x^\e w\|^2\\
  &\le C_\eta\|F_sw\|^2+\eta^2\e\|s_Y^\e w\|^2
       +CL^2\|\sqrt{u_s}s_x^\e w\|^2\\
  &\lesssim\|F_sw\|^2+L^2\|\sqrt{u_s}s_x^\e w\|^2,
\end{aligned}
\end{equation}

where we used the interpolation argument of
Lemma~\ref{lem:basic-recovery} and Lemma~\ref{lem:weighted-poincare} to get

\begin{equation}\label{eq:sw-scalar-interpolation}
\begin{aligned}
 \|s^\e w\| &\le\eta\sqrt\e\|s_Y^\e w\|
                         +C\eta^{-2}\|\sqrt{u_s}s^\e w\|\\
  &\le\eta\sqrt\e\|s_Y^\e w\|
                         +C\eta^{-2}L\|\sqrt{u_s}s_x^\e w\|.
\end{aligned}
\end{equation}

Here $\eta>0$ is fixed sufficiently small, and the
$\eta^2\e\|s_Y^\e w\|^2$ term is absorbed on the left.
Multiplying \eqref{eq:swirl-equation} by $\e s_x^\e w^2$ gives

\begin{equation}\label{eq:sw-Hx-preestimate}
\begin{aligned}
 \e\|\sqrt{u_s}s_x^\e w\|^2
 +\tfrac{\e^2}{2}\|s_x^\e w\|_{x=0}^2
 +\tfrac{\e^2}{2}\|s_Y^\e w\|_{x=L}^2 &=(F_s,\e s_x^\e w^2)-(w_ss_Y^\e,\e s_x^\e w^2)
       +2\left(\frac{w_s}{r}s^\e,\e s_x^\e w^2\right)
\\
  &\le C\|F_sw\|^2+C\e^2\|s_x^\e w\|^2
       +C\eta^{-1}\e^2\|s_Y^\e w\|^2\\
 &\quad +\eta\e\|\sqrt{u_s}s_x^\e w\|^2
       +C\e^{\frac{3}{2}}L\|\sqrt{u_s}s_x^\e w\|^2\\
  &\lesssim\|F_sw\|^2+\e^2\|\nabla s^\e w\|^2.
\end{aligned}
\end{equation}

The last two terms are absorbed since $\eta\ll1$ and $L\ll1$.
Multiply \eqref{eq:sw-H-preestimate} by a sufficiently large
constant times $\e$ and add \eqref{eq:sw-Hx-preestimate}.
Using $L\ll1$, we get

\[
\begin{aligned}
 \e^2\|\nabla s^\e w\|^2
 +\e\|\sqrt{u_s}s_x^\e w\|^2
 +\e\|\sqrt{u_s}s^\e w\|_{x=L}^2 +\e^2\|s_x^\e w\|_{x=0}^2
 +\e^2\|s_Y^\e w\|_{x=L}^2 &\lesssim\|F_sw\|^2.
\end{aligned}
\]

By \eqref{eq:swirl-equation} and
Lemma~\ref{lem:swirl-principal-solvability},

\begin{equation}\label{eq:sw-L2-elliptic}
\begin{aligned}
 \e^3\bigl(\|s_{xx}^\e w\|^2+2\|s_{xY}^\e w\|^2
                         +\|\cL s^\e\,w\|^2\bigr) &\lesssim\e\|\sqrt{u_s}s_x^\e w\|^2
       +\e^2\|\nabla s^\e w\|^2+\e\|F_sw\|^2
 \lesssim\|F_sw\|^2.
\end{aligned}
\end{equation}

\end{proof}

\subsection{Coupled linear estimates}
\label{subsec:swirl-coupling}

Consider \eqref{eq:main-equation} and \eqref{eq:swirl-equation} with

\begin{equation}\label{eq:coupled-frozen-forcing}
\begin{aligned}
 F_G &=-\e^{-\frac72}\mathfrak R_\beta
       +F_G^{\mathrm L}[s^\e]
       +F_G^{\mathrm{NL}},\\
 F_s &=-\e^{-\frac72}\mathfrak R_S
       +F_s^{\mathrm L}[G]
       +F_s^{\mathrm{NL}}.
\end{aligned}
\end{equation}

Here $F_G^{\mathrm{NL}}$ and $F_s^{\mathrm{NL}}$ are treated as prescribed source terms.
The linear coupling terms satisfy

\begin{align}\label{eq:sw-C-bound}
 \left|\left(F_s^{\mathrm L}[G],hw^2\right)\right|\nonumber &\le C\bigl(\|Gw\|+\|G_Yw\|+\|\frac{G_x}{Y}w\|\bigr)
                            \|\sqrt{u_s}hw\|\nonumber
\\
  &\le C\e^{-\frac{1}{2}}\cD[q]\|\sqrt{u_s}hw\|,\\
 &\sqrt\e\|F_G^{\mathrm L}[s^\e]w\|
 \le C\sqrt\e\bigl(\|\sqrt{u_s}s_x^\e w\|+\|s^\e w\|\bigr)
 \le C\mathcal S_\e[s^\e].\label{eq:sw-feedback-bound}
\end{align}

Here \eqref{eq:sw-feedback-bound} follows from
\eqref{eq:sw-scalar-interpolation}. For \eqref{eq:sw-C-bound}, use
\[
 |(s_s)_x|+|Y(r(s_s)_Y-2s_s)|\le C\sqrt{u_s}
\]
from Theorem~\ref{thm:approximate-profiles}, followed by
\eqref{systemat1} and Hardy's inequality.
Since $u_s$ is bounded, \eqref{eq:sw-C-bound} also gives
\[
 \|F_s^{\mathrm L}[G]w\|\le C\e^{-\frac{1}{2}}\cD[q].
\]

\begin{proposition}[Coupled linear estimate]
\label{prop:sw-coupled-estimates}
Given source terms $F_G^{\mathrm{NL}}$ and $F_s^{\mathrm{NL}}$ with
$F_G^{\mathrm{NL}}w,F_s^{\mathrm{NL}}w\in L^2(\Omega)$, the coupled
linear system with the forcing in \eqref{eq:coupled-frozen-forcing}
has a unique solution $(G,s^\e)\in\mathscr Z_\e$ satisfying

\begin{align}\label{eq:sw-linear-final}
 \cD[q]+\mathcal S_\e[s^\e] &\le C_K\Bigl[
 \e^{-\frac72}\bigl(\sqrt\e\|\mathfrak R_\beta w\|
                                      +\|\mathfrak R_Sw\|\bigr) +\sqrt\e\|F_G^{\mathrm{NL}}w\|
                       +\|F_s^{\mathrm{NL}}w\|\Bigr].
\end{align}

For two choices of the prescribed sources, the same estimate holds for
the difference of the solutions, with $F_G^{\mathrm{NL}}$ and
$F_s^{\mathrm{NL}}$ replaced by their differences and the residual terms
omitted.
\end{proposition}
\begin{proof}
Proposition~\ref{prop:apriori-closure} and
\eqref{eq:sw-feedback-bound} give
\begin{align}&\label{eq:main-feedback}
 \cD[q]\le C\sqrt\e L^{\frac{1}{16}}
       \bigl(\e^{-\frac72}\|\mathfrak R_\beta w\|
                   +\|F_G^{\mathrm{NL}}w\|\bigr)
                  +CL^{\frac{1}{16}}\mathcal S_\e[s^\e].
\end{align}
By \eqref{eq:sw-C-bound} and \eqref{eq:sw-scalar-interpolation},

\begin{align*}
 |(F_s^{\mathrm L}[G],\e s^\e w^2)| &\le CL\sqrt\e\,\cD[q]\|\sqrt{u_s}s_x^\e w\|,\\
 |(F_s^{\mathrm L}[G],\e s_x^\e w^2)| &\le C\sqrt\e\,\cD[q]\|\sqrt{u_s}s_x^\e w\|,\\
 \e\|F_s^{\mathrm L}[G]w\|^2 &\lesssim\cD[q]^2.
\end{align*}

Using these bounds in the proof of
Lemma~\ref{prop:sw-scalar-estimates} gives
\begin{align}&\label{eq:sw-scalar-feedback}
 \mathcal S_\e[s^\e]
 \le C\bigl(\e^{-\frac72}\|\mathfrak R_Sw\|
                   +\|F_s^{\mathrm{NL}}w\|\bigr)
       +C\cD[q].
\end{align}
Combining \eqref{eq:main-feedback} and
\eqref{eq:sw-scalar-feedback}, we absorb
$CL^{\frac{1}{16}}\cD[q]$ and obtain \eqref{eq:sw-linear-final}.

Set $G^{(0)}=0$ and write $G^{(j)}=u_s q^{(j)}$.
Given $G^{(j)}$, solve \eqref{eq:swirl-equation} for $s^{(j)}$, then
use $s^{(j)}$ to solve \eqref{eq:main-equation} for $G^{(j+1)}$.
Propositions~\ref{prop:swirl-solvability} and \ref{prop:main-solvability}
justify both steps with the sources in \eqref{eq:coupled-frozen-forcing}.
The differences satisfy

\begin{align*}
 \cD[q^{(j+1)}-q^{(j)}] &\le CL^{\frac{1}{16}}\mathcal S_\e[s^{(j)}-s^{(j-1)}]\\
  &\le CL^{\frac{1}{16}}
       \cD[q^{(j)}-q^{(j-1)}],
       \qquad j\ge1.
\end{align*}

Since $L\ll1$ and $\mathscr Z_\e$ is complete by
Lemma~\ref{lem:energy-norm-equivalence},
the iteration converges. The same estimate gives uniqueness.
For different prescribed sources, subtract the equations and repeat
the estimate. The residuals cancel.
\end{proof}

\section{Proof of Theorem~\ref{thm:main}}\label{sec:nonlinear}
\label{subsec:swirl-nonlinear}

\begin{proof}[Proof of Theorem~\ref{thm:main}]
For $(G,s^\e)\in\mathscr Z_\e$, \eqref{systemat1} and
\eqref{eq:grad-G-Linfty-H4eps} give
\[
\begin{aligned}
 &\|\bigl((2G-rG_Y)(\cA G)_x+rG_x(\cA G)_Y\bigr)w\|\\
 &\le C\|G\|_{W^{1,\infty}}
       \bigl(\|G_{xxx}w\|+\|\cL G_xw\|
                     +\|G_{xxY}w\|+\|(\cL G)_Yw\|\bigr)\\
 &\le C_L\e^{-\frac74}\cD[q]\,\e^{-\frac32}\cD[q]
   =C_L\e^{-\frac{13}{4}}\cD[q]^2.
\end{aligned}
\]
By $H^2(\Omega^\sharp)\hookrightarrow L^{10}(\Omega^\sharp)$ and
$H^1(\Omega^\sharp)\hookrightarrow L^{\frac{5}{2}}(\Omega^\sharp)$,
\begin{align*}
&\|\partial_x((s^\e)^2)w\|
 \le C\|(s^\e)^\sharp\|_{L^{10}}
          \|(s_x^\e)^\sharp\|_{L^{\frac{5}{2}}}\\
 &\le C_L\|s^\e\|_{H^2(w)}^2
 \le C_L\e^{-3}\mathcal S_\e[s^\e]^2.
\end{align*}
By \eqref{eq:grad-G-Linfty-H4eps},
\begin{align}&\label{eq:sw-N-transport}
 \|[(2G-rG_Y)s_x^\e+rG_xs_Y^\e-2G_xs^\e]w\|\nonumber\\
 &\quad\le C\|G\|_{W^{1,\infty}}
            \bigl(\|\nabla s^\e w\|+\|s^\e w\|\bigr)\nonumber\\
 &\quad\le C_L\e^{-\frac{11}{4}}\cD[q]\mathcal S_\e[s^\e].
\end{align}
Bilinearity gives the corresponding difference estimates.

Let $\mathcal T(\bar G,\bar s^\e)$ be the solution given by
Proposition~\ref{prop:sw-coupled-estimates}, with the nonlinear terms
in \eqref{eq:coupled-frozen-forcing} evaluated at $(\bar G,\bar s^\e)$. By
\eqref{eq:sw-residual-target},
$\e^{-\frac{7}{2}}\bigl(\sqrt\e\|\mathfrak R_\beta w\|+\|\mathfrak R_Sw\|\bigr)\le C\e$.
The product estimates above give
\begin{align*}
 &\sqrt\e\|F_G^{\mathrm{NL}}[G,s^\e]w\|
 \le C_{L,K}\bigl(\e^{\frac34}\cD[q]^2
                              +\e\mathcal S_\e[s^\e]^2\bigr),\\
 &\|F_s^{\mathrm{NL}}[G,s^\e]w\|
 \le C_{L,K}\e^{\frac{3}{4}}\cD[q]\mathcal S_\e[s^\e].
\end{align*}
On the closed unit ball, Proposition~\ref{prop:sw-coupled-estimates} gives
\begin{align*}
 &\|\mathcal T(G,s^\e)\|_{\mathscr Z_\e}
 \le C_{L,K}\bigl(\e+\e^{\frac{3}{4}}\|(G,s^\e)\|_{\mathscr Z_\e}^2\bigr),\\
 &\|\mathcal T(G_1,s_1^\e)-\mathcal T(G_2,s_2^\e)\|_{\mathscr Z_\e}
 \le C_{L,K}\e^{\frac{3}{4}}
      \|(G_1-G_2,s_1^\e-s_2^\e)\|_{\mathscr Z_\e}.
\end{align*}
For small $\e$, $C_{L,K}(\e+\e^{\frac{3}{4}})\le1$ and
$C_{L,K}\e^{\frac{3}{4}}\le\frac12$. Banach's theorem gives a unique
solution of \eqref{eq:sw-coupled-remainder} in the closed unit ball.
At this fixed point,
\[
 \|(G,s^\e)\|_{\mathscr Z_\e}
 \le C_{L,K}\e+\frac12\|(G,s^\e)\|_{\mathscr Z_\e},
\]
so $\cD[q]+\mathcal S_\e[s^\e]\le2C_{L,K}\e$.
By \eqref{eq:grad-G-Linfty-H4eps},
\[
 \|U-u_s\|_\infty+\|W-w_s\|_\infty
 \le C_{L,K}\e^{\frac{7}{2}}\|G\|_{W^{1,\infty}}
 \le C_{L,K}\e^{\frac{7}{4}}\cD[q]\le C_{L,K}\e^{\frac{11}{4}}.
\]

Cartesian differentiation gives
\[
 \|rs^\e\mathbf e_\theta\|_{H^2((0,L)\times B_2)}
 \lesssim
 \sum_{|\alpha|\le2}\|\partial_{x,Y}^{\alpha}s^\e w\|
 +\|r^{\frac{1}{2}}\{s^\e,s_x^\e,s_Y^\e\}\|.
\]
Since $s^\e|_{Y=0}=s_x^\e|_{Y=0}=0$, Hardy's inequality gives
\[
 \int_0^1 r|s^\e|^2\,dY\le C\int_0^1r^3|s_Y^\e|^2\,dY,
 \qquad
 \int_0^1 r|s_x^\e|^2\,dY\le C\int_0^1r^3|s_{xY}^\e|^2\,dY.
\]
By \eqref{equiv-norm},
\[
 \|r^{\frac{1}{2}}s_Y^\e\|
 =\left\|\frac{s_Y^\e}{r}w\right\|
 \le C\|\cL s^\e\,w\|
 \le C\|s^\e\|_{H^2(w)}.
\]
Thus, by \eqref{eq:radial-sobolev-norms}, $\|rs^\e\,\mathbf e_\theta\|_{H^2((0,L)\times B_2)} \le C\|s^\e\|_{H^2(w)}$. Since $H^2((0,L)\times B_2)\hookrightarrow L^\infty$,
\[
 \|\Theta-rs_s\|_\infty
 \le C_L\e^{\frac72}\|s^\e\|_{H^2(w)}
 \le C_L\e^2\mathcal S_\e[s^\e]\le C_{L,K}\e^3.
\]
This proves \eqref{M2}.
\end{proof}

\bibliographystyle{plain}
\bibliography{references}
\end{document}